\documentclass[10pt,a4paper,makeidx]{amsart}%
\RequirePackage{ifthen}
\newboolean{omarfont}
\setboolean{omarfont}{false}
\newboolean{usesvjour}
\setboolean{usesvjour}{false}%
    \RequirePackage{ifthen}
    \provideboolean{useutopia}
    \setboolean{useutopia}{false}%
    \provideboolean{usesvjour}
    \setboolean{usesvjour}{false}%
    \ifthenelse{\boolean{useutopia}}{%
      \RequirePackage{ifthen}
      \RequirePackage{iftex}
      \provideboolean{usemathrsfs}%
      \provideboolean{useutopia}
      \setboolean{useutopia}{false}%
      \ifXeTeX
        \RequirePackage[american,british]{babel}
        \RequirePackage{xltxtra}%
        \DeclareSymbolFont{usualmathcal}{OMS}{cmsy}{m}{n}
        \DeclareSymbolFontAlphabet{\mathcalbf}{usualmathcal}
        \ifthenelse{\boolean{useutopia}}{
          \RequirePackage[utopia,euro]{mathdesign}
          \RequirePackage[OMLmathbf,OMLmathsfit]{isomath}
        }{
          \RequirePackage{amssymb}    %
          \RequirePackage{bbold}    %
        }
        \RequirePackage{mathrsfs} %
        \setboolean{usemathrsfs}{true}%
      \else%
        \RequirePackage[utf8]{inputenc}%
        \ifthenelse{\boolean{useutopia}}{%
          \RequirePackage[utopia,euro]{mathdesign}%
          \RequirePackage[OMLmathrm,OMLmathbf,OMLmathsfit]{isomath}
          \def\nobreakspace{\nobreak\space\nobreak}%
          \DeclareSymbolFont{usualmathcal}{OMS}{cmsy}{m}{n}
          \DeclareSymbolFontAlphabet{\mathcalbf}{usualmathcal}
          \RequirePackage{stmaryrd} %
          \providecommand{\diracdelta}[1][]{\ensuremath{\deltaup_{#1}}}
          
          \providecommand{\lap}{\ensuremath{\Deltaup}}
          \providecommand{\measure}[1]{\ensuremath{\mathcalbf{\uppercase{#1}}}}
        }{%
          \RequirePackage{amssymb}    %
          \RequirePackage{mathrsfs} %
          \RequirePackage{bbold}    %
          \RequirePackage{stmaryrd} %
          \setboolean{usemathrsfs}{true}%
          \providecommand{\mathcalbf}{\mathcal}
        }%
        \RequirePackage[american,british]{babel}
      \fi%
      \RequirePackage{mathtools}%
      \RequirePackage{nicefrac}
      \RequirePackage{stmaryrd} %
      \RequirePackage{amsthm}%
      \RequirePackage[shortlabels]{enumitem}%
      \RequirePackage{xspace}%
      \RequirePackage{verbatim}%
      \RequirePackage{fancyvrb}
      \RequirePackage{listings}%
      \RequirePackage{xcolor}
      \RequirePackage{epsfig}
      \RequirePackage{graphicx}
      \RequirePackage[space]{grffile}%
      \RequirePackage{booktabs}%
      \RequirePackage{tikz}%
      \usetikzlibrary{calc}%
      \provideboolean{usebeamer}%
      \provideboolean{isthesis}%
      \provideboolean{isamsltex}%
      \provideboolean{issiamltex}%
      \providecommand{\nobreakspace}{~}
    }{%
      \RequirePackage{ifthen}
      \RequirePackage{iftex}
      \provideboolean{usemathrsfs}%
      \provideboolean{useutopia}
      \setboolean{useutopia}{false}%
      \ifXeTeX
        \RequirePackage[american,british]{babel}
        \RequirePackage{xltxtra}%
        \DeclareSymbolFont{usualmathcal}{OMS}{cmsy}{m}{n}
        \DeclareSymbolFontAlphabet{\mathcalbf}{usualmathcal}
        \ifthenelse{\boolean{useutopia}}{
          \RequirePackage[utopia,euro]{mathdesign}
          \RequirePackage[OMLmathbf,OMLmathsfit]{isomath}
        }{
          \RequirePackage{amssymb}    %
          \RequirePackage{bbold}    %
        }
        \RequirePackage{mathrsfs} %
        \setboolean{usemathrsfs}{true}%
      \else%
        \RequirePackage[utf8]{inputenc}%
        \ifthenelse{\boolean{useutopia}}{%
          \RequirePackage[utopia,euro]{mathdesign}%
          \RequirePackage[OMLmathrm,OMLmathbf,OMLmathsfit]{isomath}
          \def\nobreakspace{\nobreak\space\nobreak}%
          \DeclareSymbolFont{usualmathcal}{OMS}{cmsy}{m}{n}
          \DeclareSymbolFontAlphabet{\mathcalbf}{usualmathcal}
          \RequirePackage{stmaryrd} %
          \providecommand{\diracdelta}[1][]{\ensuremath{\deltaup_{#1}}}
          
          \providecommand{\lap}{\ensuremath{\Deltaup}}
          \providecommand{\measure}[1]{\ensuremath{\mathcalbf{\uppercase{#1}}}}
        }{%
          \RequirePackage{amssymb}    %
          \RequirePackage{mathrsfs} %
          \RequirePackage{bbold}    %
          \RequirePackage{stmaryrd} %
          \setboolean{usemathrsfs}{true}%
          \providecommand{\mathcalbf}{\mathcal}
        }%
        \RequirePackage[american,british]{babel}
      \fi%
      \RequirePackage{mathtools}%
      \RequirePackage{nicefrac}
      \RequirePackage{stmaryrd} %
      \RequirePackage{amsthm}%
      \RequirePackage[shortlabels]{enumitem}%
      \RequirePackage{xspace}%
      \RequirePackage{verbatim}%
      \RequirePackage{fancyvrb}
      \RequirePackage{listings}%
      \RequirePackage{xcolor}
      \RequirePackage{epsfig}
      \RequirePackage{graphicx}
      \RequirePackage[space]{grffile}%
      \RequirePackage{booktabs}%
      \RequirePackage{tikz}%
      \usetikzlibrary{calc}%
      \provideboolean{usebeamer}%
      \provideboolean{isthesis}%
      \provideboolean{isamsltex}%
      \provideboolean{issiamltex}%
      \providecommand{\nobreakspace}{~}
    }
    \colorlet{a}{magenta}
    \colorlet{b}{green!75!blue}
    \colorlet{c}{yellow!87.5!red}
    \colorlet{d}{cyan}
    \colorlet{e}{red}
    \colorlet{f}{blue}
    \colorlet{g}{white}
    \colorlet{h}{black!50}
    \colorlet{i}{black}
    \colorlet{j}{black!75}

    \RequirePackage[obeyspaces]{url}
    \RequirePackage{hyperref}
    \providecommand{\linkedurl}[1]{\url{#1}}
    \providecommand{\linkedemail}[1]{\href{mailto:#1}{#1}}
    
    \providecommand{\email}[1]{{\linkedemail{#1}}}

    \providecommand{\Ignore}[1]{}
    \providecommand{\ignore}[1]{}
    \providecommand{\freeze}[1]{}%
    \providecommand{\crossout}[1]{{\textcolor{i!20}{#1}}}
    \providecommand{\highlightcolor}{a}
    \providecommand{\highlight}[1]{{\color{\highlightcolor}#1}}

    \providecommand{\memo}[1]{%
      \ensuremath{%
        \framebox{\tiny\textbf{\kern-2pt\textsf{#1}}\kern-2pt}%
      }%
      \xspace}

    \RequirePackage{alphalph}
    \provideboolean{shownotes}
    \setboolean{shownotes}{true}%
    \newcounter{margnote}[page]
    \providecommand{\mgcolor}{a}
    \providecommand{\mgcolorset}[1]{\renewcommand{\mgcolor}{\alphalph{#1}}}
    \providecommand{\mgcolorsetbycounter}[1]{%
      \ifthenelse{\value{#1}<11}{%
        \renewcommand{\mgcolor}{\alph{#1}}%
      }{%
        \renewcommand{\mgcolor}{a}}%
    }
    \providecommand{\mgcolormake}{\mgcolorsetbycounter{margnote}}
    \providecommand{\mgcolorstepby}[1]{
      \setcounter{tmpcounter}{\value{margnote}}%
      \addtocounter{tmpcounter}{#1}%
      \mgcolorsetbycounter{tmpcounter}%
    }%
    \providecommand{\margnotecolor}{%
      \ifthenelse{\value{margnote}=0}{%
        \mgcolorset{10}
      }{%
        \ifthenelse{\value{margnote}<7}{%
          \mgcolormake%
        }{%
          \ifthenelse{\value{margnote}=7}{\mgcolorset{10}}{%
            \ifthenelse{\value{margnote}<11}{\mgcolormake}{%
              \ifthenelse{\value{margnote}<17}{\mgcolorstepby{-10}}{%
                \mgcolorset{10}%
              }%
            }%
          }%
        }%
      }%
    }%
    \providecommand{\margnotemark}{{\colorbox{\mgcolor}{\tiny\color{g}\upshape\texttt{\arabic{page}.\arabic{margnote}}}}}
    \providecommand{\margnote}[2][]{%
      \ifthenelse{%
        \boolean{shownotes}%
      }{%
        \stepcounter{margnote}%
        \margnotecolor%
        \margnotemark %
        \marginpar{%
          \color{\mgcolor}%
          \texttt{%
            \begin{minipage}{2cm}%
              \raggedright\tiny%
              \margnotemark%
              #2%
              \\
              {\ifx|#1|{}\else{ - #1}\fi}%
            \end{minipage}%
          }%
        }%
      }{%
      }%
    }%
    \providecommand{\mathnote}[2][]{%
      \ifthenelse{%
        \boolean{shownotes}%
      }{%
        \stepcounter{margnote}%
        \margnotecolor%
        \text{%
          \colorbox{\mgcolor}{%
            \color{g}%
            \texttt{%
              \tiny%
                  \margnotemark: %
                  \ifx|#1|{}\else{#1:}\fi%
                  #2%
            }%
          }%
        }%
      }{%
      }%
    }%
    \providecommand{\textnote}[2][]{%
      \ifthenelse{%
        \boolean{shownotes}%
      }{%
        \stepcounter{margnote}%
        \margnotecolor%
        \ \\
        \text{%
          \colorbox{\mgcolor}{%
            \begin{minipage}{.9\textwidth}
            \color{g}%
            \texttt{%
              \margnotemark: %
              \ifx|#1|{}\else{#1: }\fi%
              #2%
            }%
            \end{minipage}
          }%
        }%
      }{%
      }%
    }%

    \providecommand{\Todo}[1]{
      \ifthenelse{\boolean{shownotes}}{
        \begin{center}
        \begin{tikzpicture}
         \node[fill=a!17]{
           \begin{minipage}{\textwidth}
             \texttt{To do:}
             \\
             \texttt{\bfseries{\small #1}}
           \end{minipage}
         };
        \end{tikzpicture}
        \end{center}
      }{}}
    \provideboolean{showrevisions}
    \setboolean{showrevisions}{true}
    \newcommand{\revisionsheader}{***\newline\Warning{the following part is under revision}}
    \newcommand{\revisionsfooter}{***\newline\Warning{end of part under revision}}

    \providecommand{\Warning}[1]{    
      \begin{tikzpicture}
        \node[fill=a!27]{
          \begin{minipage}{\textwidth}
            \texttt{\bfseries{\small Warning: #1}}
          \end{minipage}
        };
      \end{tikzpicture}
    }

    \provideboolean{showcomments}
    \providecommand{\margincomment}[1]{
    \ifthenelse{\boolean{showcomments}}{\marginpar{\tiny #1}}{}
    }
    \provideboolean{showchanges}
    \setboolean{showchanges}{false}
    \providecommand{\changes}[2][]{%
      \ifthenelse{\boolean{showchanges}}{{\ifx|#1|{}\else\margnote{#1}\fi\highlight{#2}}}{#2}}
    \providecommand{\mathchanges}[2][]{%
      \ifthenelse{\boolean{showchanges}}{{\ifx|#1|{}\else\mathnote{#1}\fi\highlight{#2}}}{#2}}

    \providecommand{\changefromto}[3][replace with]{%
      \ifthenelse{\boolean{showchanges}}{%
        {\crossout{#2}\margnote{#1}}{\highlight{#3}}}{%
        #3\xspace}%
    }
    \providecommand{\ChangePar}[3][]{%
      \ifthenelse{\boolean{showchanges}}{
        {\par\textcolor{i!20}{#2}\ifx|#1|\else\margnote{#1}\fi}{\par\textcolor{a}{#3}}
      }{%
        \par #3%
      }%
    }
    \providecommand{\InsertPar}[1]{
      \ifthenelse{\boolean{showchanges}}
      {{\par$\mapsto$ \textcolor{blue}{#1}}}
      {\par #1}
    }
    
    \providecommand{\mathchangefromto}[3][]{\crossout{#2}\ifx|#1|\else\mathnote{#1}\fi\highlight{#3}}

    \providecommand{\mathscript}
    	   {\mathscr}
     \providecommand{\bbbold}{\mathbb}

     \providecommand{\rN}{\ensuremath{\bbbold N}\xspace}
     
     \providecommand{\rP}{\ensuremath{\bbbold P}\xspace}
     
     \providecommand{\rR}{\ensuremath{\bbbold R}\xspace}
     
     \providecommand{\rT}{\ensuremath{\bbbold T}\xspace}

    \providecommand{\Ae}[1][]{\ensuremath{\ifx|#1|{\ }\else{\:#1\text{-}}\fi\text{almost everywhere }}\xspace}
    \providecommand{\Aa}[1][]{\ensuremath{\text{ for }\ifx|#1|{}\else{\:#1\text{-}}\fi\text{almost all }}}
    \providecommand{\as}[1][]{\ensuremath{\ifx|#1|{\ }\else{#1\text{-}}\fi\text{almost surely}}\xspace}
    \providecommand{\aposteriori}{aposteriori\xspace}
    
    \providecommand{\apriori}{{apriori}\xspace}

     \providecommand{\naturals}{\ensuremath{\rN}}
     
     \providecommand{\NO}[1][]{\ensuremath{\naturals_0\ifx|#1|{}\else^{#1}\fi}}

     \providecommand{\reals}{\rR}

     \providecommand{\R}[1]{\reals^{#1}}
     
     \providecommand{\fieldmats}[3][F]{\csname#1\endcsname{#2\times#3}}
     \providecommand{\realmats}[2]{\fieldmats[R]{#1}{#2}}
     
     \providecommand{\RO}[1][]{{\reals_{0+}\ifx|#1|{}\else^{#1}\fi}}
     \providecommand{\RP}[1][]{{\reals_+\ifx|#1|\else^{#1}\fi}}

     \providecommand{\torus}[1]{\rT\ifthenelse{\equal{#1}1}{}{^#1}}

     \providecommand{\diracdelta}[1][]{\ensuremath{{\mathrm{\delta}}\ifx|#1|{}\else_{#1}\fi}}

     \providecommand{\pic}{\ensuremath{\mathrm\pi}}
     
     \providecommand{\pifracl}[2][]{\fracl{\ifx|#1|\else#1\fi\pic}{#2}}
     \providecommand{\pifrac}[2][]{\frac{\ifx|#1|\else#1\fi\pic}{#2}}

     \providecommand{\closure}[1]{\overline{#1}}
     \providecommand{\inner}{\cdot}
     \providecommand{\outerp}{\wedge}

     \providecommand{\frobinner}{\!:\!}

     \providecommand{\W}{\ensuremath{\varOmega}\xspace}

     \providecommand{\w}{\ensuremath{\omega}\xspace}

     \providecommand{\qgroup}[1]{{#1}}%

     \providecommand{\qp}[2][]{\ensuremath{\ifx|#1|\left(\else\csname#1\endcsname(\fi{#2}\ifx|#1|\right)\else\csname#1\endcsname)\fi}}

     \providecommand{\qpreg}[1]{\ensuremath{(#1)}}
     \providecommand{\qpbig}[1]{\qp[big]{#1}}%
     \providecommand{\qpBig}[1]{\ensuremath{\Big(#1\Big)}}
     \providecommand{\qpbigg}[1]{\ensuremath{\bigg(\!#1\!\bigg)}}
     \providecommand{\qpBigg}[1]{\ensuremath{\Bigg(\!#1\!\Bigg)}}
     \providecommand{\qb}[2][]{\ifx|#1|\left[\else\csname#1\endcsname[\fi{#2}\ifx|#1|\right]\else\csname#1\endcsname]\fi}
     \providecommand{\qc}[2][]{\ensuremath{\ifx|#1|\left\{\else\csname#1\endcsname\{\fi{#2}\ifx|#1|\right\}\else\csname#1\endcsname\}\fi}}

     \providecommand{\qa}[1]{\ensuremath{\left\langle{#1}\right\rangle}}
     \providecommand{\qareg}[1]{\ensuremath{\langle#1\rangle}}
     \providecommand{\qabig}[1]{\ensuremath{\big\langle#1\big\rangle}}
     \providecommand{\qaBig}[1]{\ensuremath{\Big\langle#1\Big\rangle}}
     \providecommand{\qabigg}[1]{\ensuremath{\bigg\langle#1\bigg\rangle}}
     \providecommand{\qaBigg}[1]{\ensuremath{\Bigg\langle#1\Bigg\rangle}}
     \providecommand{\opinter}[2]{\ensuremath{\left(#1,#2\right)}\xspace}

     \providecommand{\expp}[1]{\ensuremath{\e^{#1}}}
     
     \providecommand{\compowqp}[2]{\ensuremath{\qp{\!#2\!\!}^{\kern -.4em #1}\!}}
     
     \providecommand{\powqpreg}[2]{\ensuremath{%
         \qpreg{#2}^{\kern 0em\lower .1ex\hbox{\scriptsize $#1$}}\kern-.3em}}
     \providecommand{\powqpbig}[2]{\ensuremath{%
         \qpbig{#2}^{\kern -.2em\lower .3ex\hbox{\scriptsize $#1$}}\kern-.3em}}
     \providecommand{\powqpBig}[2]{\ensuremath{%
         \qpBig{#2}^{\kern -.2em\lower .3ex\hbox{\scriptsize $#1$}}\kern-.3em}}
     \providecommand{\powqpbigg}[2]{\ensuremath{%
         \qpbigg{#2}^{\kern -.2em\lower .3ex\hbox{\scriptsize $#1$}}\kern-.3em}}
     \providecommand{\powqpBigg}[2]{\ensuremath{%
         \qpBigg{#2}^{\kern -.2em\lower .3ex\hbox{\scriptsize $#1$}}}}
     \providecommand{\powp}[3][]{{#3}\ifx|#1|^{#2}\else{#1}^{#2}\fi}%
     \providecommand{\pow}[2][]{\ifx|#1|\operatorname{pow}^{#2}\else\powp{#2}{#1}\fi}%
     \providecommand{\ppow}[3][]{\powp[#1]{#3}{#2}}

     \providecommand{\powqp}[3][]{\powp[#1]{#2}{\qp{#3}}}%
     \providecommand{\qppow}[3][]{\ppow[#1]{\qp{#2}}{#3}}%

     \providecommand{\pownorm}[2]{\powp{#1}{\norm{#2}}}

     \providecommand{\qpsqrt}[1]{\powqp{\fracl12}{#1}}

     \providecommand{\qpsqrt}[2][2]{\powqpsqrt[#1]{#2}}

     \providecommand{\norm}[2][]{\ifx|#1|\left|\else\csname#1\endcsname|\fi#2\ifx|#1|\right|\else\csname#1\endcsname|\fi}
     \providecommand{\normon}[2]{\norm{#1}_{#2}}

     \providecommand{\normreg}[1]{\ensuremath{|#1|}}
     \providecommand{\abs}[2][]{\ensuremath{\ifx|#1|{\left|#2\right|}\else{\csname#1\endcsname|{#2}\csname#1\endcsname|}\fi}}

     \providecommand{\Norm}[2][]{\ifx|#1|\left\|\else\csname#1\endcsname\|\fi{#2}\ifx|#1|\right\|\else\csname#1\endcsname\|\fi}
     
     \providecommand{\Normbig}[1]{\ensuremath{\big\|#1\big\|}}
     
     \providecommand{\Normon}[2]{\Norm{#1}_{#2}}

     \providecommand{\normonsob}[4][]{\normon{#2}{\sob{#3}{#4}\if|#1|{}\else(#1)\fi}}
     \providecommand{\Normonsob}[4][]{\Normon{#2}{\sob{#3}{#4}\if|#1|{}\else(#1)\fi}}
     \providecommand{\Normonleb}[3][]{\Normon{#2}{\leb{#3}\if|#1|{}\else(#1)\fi}}

     \providecommand{\ltwop}[3][]{\ensuremath{\qa{#2,#3}\ifx|#1|\else_{#1}\fi}}
     \providecommand{\ltwopreg}[2]{\ensuremath{\qareg{#1,#2}\ifx|#1|\else_{#1}\fi}}
     \providecommand{\ltwopbig}[2]{\ensuremath{\qabig{#1,#2}\ifx|#1|\else_{#1}\fi}}
     \providecommand{\ltwopBig}[2]{\ensuremath{\qaBig{#1,#2}\ifx|#1|\else_{#1}\fi}}
     \providecommand{\ltwopbigg}[2]{\ensuremath{\qabigg{#1,#2}\ifx|#1|\else_{#1}\fi}}
     \providecommand{\ltwopBigg}[2]{\ensuremath{\qaBigg{#1,#2}\ifx|#1|\else_{#1}\fi}}

     \providecommand{\average}[2][]{{\qa{#2}\ifx|#1|\else_{#1}\fi}}

     \providecommand{\ensemble}[2]{\ensuremath{\left\{ #1:\;#2 \right\}}}
     \providecommand{\setofsuch}{\ensemble}%

     \providecommand{\ceil}[1]{\ensuremath{\left\lceil{#1}\right\rceil}}

     \providecommand{\setof}[1]{{\qc{#1}}}
     \providecommand{\pair}[2]{\qp{#1,#2}}
     
     \providecommand{\triple}[3]{\qp{#1,#2,#3}}

     \providecommand{\conditionalto}[1]{{\left|{#1}\right.}}

    \providecommand{\measure}[1]{\ensuremath{\mathcalbf{\MakeUppercase{#1}}}}
    
    \providecommand{\probmeasure}[2][]{{\measure{#2}}\ifx|#1|\else_{#1}\fi}
    \providecommand{\Prob}{}
    \renewcommand{\Prob}[1][]{\probmeasure[{#1}]{p}}

    \providecommand{\randvars}[1][\Prob]{\operatorname{RV}\ifx|#1|{}\else{(#1)}\fi}
    \providecommand{\discrandvars}[1][\Prob]{\operatorname{DRV}\ifx|#1|{}\else{({#1)}\fi}} 
    \providecommand{\contrandvars}[1][\Prob]{\ensuremath{\operatorname{CDRV}\ifx|#1|{}\else(#1)\fi}} 
     \def\env@matrix{\hskip -\arraycolsep
      \let\@ifnextchar\new@ifnextchar
      \array{*\c@MaxMatrixCols c}}
     \makeatletter
     \renewcommand*\env@matrix[1][c]{\hskip -\arraycolsep
       \let\@ifnextchar\new@ifnextchar
       \array{*\c@MaxMatrixCols #1}}
     \makeatother
     \providecommand{\irow}[2]{#1_{#2}}%
     \providecommand{\icol}[2]{#1^{#2}}%

     \providecommand{\ijrowcol}[3]{\icol{\irow{#1}{#2}}{#3}}
     \providecommand{\entry}[1]{\qb{#1}}
     \providecommand{\vecentry}[2]{\irow{#1}{#2}}

     \providecommand{\rowof}[1]{\qb{#1}}

     \providecommand{\getentryi}[2]{\irow{\entry{#1}}{#2}}
     
     \providecommand{\getvecentry}[2]{\getentryi{\vec #1}{#2}}

     \providecommand{\dismatof}[2][r]{\begin{bmatrix}[#1]#2\end{bmatrix}}

     \providecommand{\matentry}[3]{\ijrowcol{#1}{#2}{#3}}

     \providecommand{\block}[5]{\ijrowcol{#1}{\ifx#2#3{\rowof{#2}}\else\rowof{{#2}\dotsc{#3}}\fi}{\ifx#4#5{\rowof{#4}}\else\rowof{{#4}\dotsc{#5}}\fi}}
     \providecommand{\colblock}[3]{\getvecentry{#1}{\ifx#2#3{#2}\else\fromto{#2}{#3}\fi}}

     \providecommand{\dismatskeldots}[4]{
       \dismatof[c]{
         #1&\dotsc&#3
         \\
         \vdots & \ddots &\vdots
         \\
         #2&\dotsc&#4
       }
     }
     \providecommand{\dismatcommfromtofromto}[5]{
       \dismatskeldots{#1#2#4}{#1#3#4}{#1#2#5}{#1#3#5}
     }
     \providecommand{\dismatcustfromtofromto}[6][matentry]{
       \dismatcommfromtofromto{\csname#1\endcsname{#2}}#3#4#5#6
     }
     \providecommand{\dismatcustfromtofromto}[6][matentry]{
       \dismatskeldots{%
         \csname#1\endcsname{#2}{#3}{#4}%
       }{%
         \csname#1\endcsname{#2}{#3}{#6}%
       }{%
         \csname#1\endcsname{#2}{#5}{#4}%
       }{%
         \csname#1\endcsname{#2}{#5}{#6}%
       }%
     }%
     \providecommand{\dismatcustfromtofromto}[6][matentry]{
       \dismatof{
         \csname#1\endcsname{#2}{#3}{#4}&\dotsc&\csname#1\endcsname{#2}{#3}{#6}
         \\
         \vdots & \ddots &\vdots
         \\
         \csname#1\endcsname{#2}{#5}{#4}&\dotsc&\csname#1\endcsname{#2}{#5}{#6}
       }
     }

     \providecommand{\dissysaxbdotsnm}[5]{\begin{matrix}[r]%
         \matentry{#1}11\vecentry{#2}1&+\dotsb&+\matentry{#1}1{#5}\vecentry{#2}{#5}
         &
         =
         \ifx|#3|0\else{\vecentry {#3}1}\fi
         \\
         \dotsb
         \\
         \matentry{#1}{#4}1\vecentry{#2}1&+\dotsb&+\matentry{#1}{#4}{#5}\vecentry{#2}{#5}
         &
         =
         \ifx|#3|0\else{\vecentry {#3}{#4}}\fi
     \end{matrix}}
     \providecommand{\seqof}[1]{\qp{#1}}%
     \providecommand{\seq}[1]{\seqof{#1}}%
     \providecommand{\seqs}[2]{\seqof{#1}_{#2}}
     \providecommand{\sets}[2]{\setof{#1}_{#2}}%
     \providecommand{\seqi}[3][]{\seqs{#2_{#3}}{\ifx|#1|{#3}\else{{#3}\in{#1}}\fi}}%
     \providecommand{\subseqi}[4][]{\seqs{#2_{{#3}_{#4}}}{\ifx|#1|{#4}\else{{#4}\in{#1}}\fi}}%
     
     \providecommand{\seqsinat}[2]{\seqsi{#1}{#2}{\naturals}}

      \providecommand{\seti}[3][]{\sets{#2_{#3}}{\ifx|#1|_{#3}\else_{{#3}\in{#1}}\fi}}%
     \providecommand{\sequ}[3][]{\seqs{#2^{#3}}{\ifx|#1|{#3}\else{{#3}\in{#1}}\fi}}%
     \providecommand{\setu}[3][]{\sets{#2^{#3}}{\ifx|#1|{#3}\else{{#3}\in{#1}}\fi}}%
     \providecommand{\seqsi}[3]{\seqi[#3]{#1}{#2}}

     \providecommand{\limofat}[3][]{\ensuremath{\lim_{\ifx|#1|{}\else{#1\ni}\fi#3}{#2}}}
     \providecommand{\limsupofat}[3][]{\ensuremath{\limsup_{\ifx|#1|{}\else{#1\ni}\fi#3}{#2}}}
     \providecommand{\liminfofat}[3][]{\ensuremath{\liminf_{\ifx|#1|{}\else{#1\ni}\fi#3}{#2}}}

     \providecommand{\listdotsfrom}[3][]{\ensuremath{#2\ifx|#1|\else#1\fi,#3\ifx|#1|\else#1\fi,\dotsc}}
     \providecommand{\listdotsfromto}[3][]{\ensuremath{#2\ifx|#1|\else#1\fi,\dotsc,#3\ifx|#1|\else#1\fi}}
     \providecommand{\listifromto}[5][]{\ensuremath{{#2}_{#3}\ifx|#1|\else#1\fi},\text{ for }\ensuremath{\rangefromto{#3}{#4}{#5}}\xspace}
     \providecommand{\listufromto}[5][]{\ensuremath{{#2}^{#3}\ifx|#1|\else#1\fi},\text{ for }\ensuremath{\rangefromto{#3}{#4}{#5}}\xspace}
     \providecommand{\listitwo}[2][]{\ensuremath{#2_1\ifx|#1|\else#1\fi,#2_2\ifx|#1|\else#1\fi}}
     \providecommand{\listutwo}[2][]{\ensuremath{#2^1\ifx|#1|\else#1\fi,#2^2\ifx|#1|\else#1\fi}}
     \providecommand{\listithree}[2][]{\ensuremath{#2_1\ifx|#1|\else#1\fi,#2_2\ifx|#1|\else#1\fi,#2_3\ifx|#1|\else#1\fi}}
     \providecommand{\listuthree}[2][]{\ensuremath{#2^1\ifx|#1|\else#1\fi,#2^2\ifx|#1|\else#1\fi,#2^3\ifx|#1|\else#1\fi}}

     \providecommand{\listidotsfromto}[4][]{\listdotsfromto[#1]{#2_{#3}}{#2_{#4}}}

     \providecommand{\seqidotsfromto}[3]{\seq{\listidotsfromto{#1}{#2}{#3}}}

     \providecommand{\sums}[2]{\ensuremath{\sum_{#1\in #2}}}

     \providecommand{\sumifromto}[3]{\ensuremath{\sum_{#1=#2}^{#3}}}
     \providecommand{\sumsifromto}[4]{\ensuremath{\sumifromto{#2}{#3}{#4}{#1}_{#2}}}

     \providecommand{\jump}[2][]{\ensuremath{\left\llbracket #2\right\rrbracket\ifx|#1|{}\else_{#1}\fi}}

     \providecommand{\fromto}[2]{\ensuremath{\setof{#1\dotsc#2}}}%

     \providecommand{\integerbetween}[2]{\ensuremath{={#1},\dotsc,{#2}}}
     
     \providecommand{\rangefromto}[3]{\ensuremath{#1\integerbetween{#2}{#3}}}

     \providecommand{\e}{\ensuremath{\operatorname{e}\!}\xspace}

     \providecommand{\d}{}
     \renewcommand{\d}[1][]{\ensuremath{\operatorname{d}\!\ifx|#1|\else{_{#1}}\fi}}
     
     \providecommand{\ds}[1][]{\d{\measure S}}
     \providecommand{\D}[1][]{\ensuremath{\operatorname{D}\!\ifx|#1|\else{_{#1}}\fi}}

    \providecommand{\registered}%
    {\ensuremath{^\text{\textregistered}}}

    \providecommand{\tand}{\ensuremath{\text{ and }}}

    \providecommand{\card}{\ensuremath{\#}}

    \providecommand{\constant}[1]{\ensuremath{C_{#1}}}
    \providecommand{\constext}[2][]{\constant{\textup{#2}{\ifx|#1|{}\else{,\ensuremath{#1}}\fi}}}            %
    \providecommand{\constref}[2][]{\ensuremath{\constant{\textup{\ref{#2}{\ifx|#1|{}\else{,\ensuremath{#1}}\fi}}}}}
    \providecommand{\constdef}[2][]{\label{#2}\ensuremath{\constant{\textup{\ref{#2}{\ifx|#1|{}\else{,\ensuremath{#1}}\fi}}}}}

    \providecommand{\funkref}[3][]{\ensuremath{{#3}_{\textup{\ref{#2}{\ifx|#1|{}\else{,\ensuremath{#1}}\fi}}}}}

    \providecommand{\diam}{\operatorname{diam}}
    
    \providecommand{\curl}{\operatorname{curl}}
    \renewcommand{\div}[1][]{\nabla\ifx|#1|{}\else\kern-2pt_{#1}\fi\kern-2pt\inner}
    \providecommand{\divof}[2][]{\div[#1]\ifx|#2|{}\else\qb{#2}\fi}
    \providecommand{\divideabyb}[2]{\operatorname{div}(a,b)}
    \providecommand{\grad}{}
    \renewcommand{\grad}[1][]{\nabla\ifx|#1|\else_{#1}\fi}

    \providecommand{\rot}[1][]{\nabla\ifx|#1|\else_{#1}\fi\outerp}
    
    \providecommand{\rowdiv}[1][]{\D\ifx|#1|{}\else\kern-1pt_{#1}\kern-2pt\fi\cdot}
    \providecommand{\rowdivof}[2][]{\rowdiv[#1]\ifx|#2|{}\else\qb{#2}\fi}

    \providecommand{\inv}[1][]{\operatorname{inv}\ifx|#1|\else^{#1}\fi}
    \providecommand{\ivt}[1]{\operatorname{ivt}\ifx|#1|\else^{#1}\fi}
    \providecommand\tensorinvariant\ivt

    \providecommand{\mod}{}
    \renewcommand{\mod}[1][]{\operatorname{mod}\ifx|#1|\else\kern-1pt_{#1}\fi}
    
    \let\oldfrac\frac
    \renewcommand{\frac}[3][]{\ifx|#1|\oldfrac{#2}{#3}\else\begin{array}{#1}{#2}\\\hline{#3}\end{array}\fi}
    \providecommand{\fracl}[3][]{\ifx|#1|\nicefrac{#2}{#3}\else{#2}#1/{#3}\fi}

    \providecommand{\qpfracl}[3][]{\qp{\ifx|#1|\fracl{#2}{#3}\else{#2}#1/{#3}\fi}}
    \providecommand{\qpfrac}[3][]{\qp{\ifx|#1|\frac{#2}{#3}\else{#2}#1/{#3}\fi}}
    \providecommand{\absfracl}[3][]{\abs{\ifx|#1|\fracl{#2}{#3}\else{#2}#1/{#3}\fi}}
    \providecommand{\absfrac}[3][]{\abs{\ifx|#1|\frac{#2}{#3}\else{#2}#1/{#3}\fi}}
    \providecommand{\fraclff}[3][]{\ifx|#1|{#2}/{#3}\else{#2}#1/{#3}\fi}

    \providecommand{\eye}[1][]{\vec{\mathrm I}\ifx|#1|{}\else_{#1}\fi}%
    \providecommand{\numeye}[1][]{\boldsymbol{\mathsf{I}}\ifx|#1|{}\else_{#1}\fi}%
    \providecommand{\Eye}[1]{
      \begin{bmatrix}
      \ifthenelse{#1>1}{
        \ifthenelse{#1>2}{
          \ifthenelse{#1>3}{
            \ifthenelse{#1>4}{
              1&\zeroentry&\dotso&\zeroentry
              \\
              \zeroentry&1&\dotso&\zeroentry
              \\
              \vdots&\vdots&\ddots&\vdots
              \\
              \zeroentry&\zeroentry&\dotso&1
            }{        
              1&\zeroentry&\zeroentry&\zeroentry
              \\
              \zeroentry&1&\zeroentry&\zeroentry
              \\
              \zeroentry&\zeroentry&1&\zeroentry
              \\
              \zeroentry&\zeroentry&\zeroentry&1
            }
          }{
            1&\zeroentry&\zeroentry
            \\
            \zeroentry&1&\zeroentry
            \\
            \zeroentry&\zeroentry&1
          }
        }{
          1&\zeroentry
          \\
          \zeroentry&1
        }
      }{
        1
      }
      \end{bmatrix}
    }

    \providecommand{\lebmeas}[1][]{\measure L^{#1}}     %
    \providecommand{\lebmeasof}[2][]{\ifx|#1|\left|#2\right|\else\lebmeas[#1]\qp{#2}\fi}         %
    \providecommand{\meshsize}[1][]{h\ifx|#1|\else_{#1}\fi}
    \providecommand{\maxi}[2]{#1\vee#2}                       %
    \providecommand{\qpmaxi}[2]{\qp{\maxi{#1}{#2}}}
    \providecommand{\mini}[2]{#1\wedge#2}                     %
    \providecommand{\qpmini}[2]{\qp{\mini{#1}{#2}}}

    \providecommand{\argmin}{\operatorname{argmin}}

    \providecommand{\Oh}{\operatorname{O}}                   %
    \providecommand{\dash}[1][']{\ifthenelse{\equal{#1}{'}\OR\equal{#1}{''}}{#1}{^{(#1)}}}
    
    \providecommand{\pdfrac}[2][]{\ensuremath{\frac{\partial\ifx|#1|\phantom{#2}\else{#1}\fi}{\partial{#2}}}} %
    \providecommand{\pdfracpow}[3][]{\ensuremath{\frac{\partial^{#3}\ifx|#1|\phantom{#2}\else{#1}\fi}{\partial{#2}^{#3}}}} %

    \providecommand{\pd}[2][]{\ensuremath{\partial_{#2}}{\ifx|#1|{}\else{\qb{#1}}\fi}} %

    \renewcommand{\Im}{\operatorname{im}}                 %
    \renewcommand{\Re}{\operatorname{re}}                 %
    \providecommand{\imaginpart}[1][]{\Im{\ifx|#1|{}\else\qp{#1}\fi}} %
    \providecommand{\realpart}[1][]{\Re{\ifx|#1|{}\else\qp{#1}\fi}} %
    \providecommand{\sign}{\operatorname{sign}}               %

    \providecommand\determinant\det
    
    \providecommand{\trace}{\operatorname{tra}}             %
    \providecommand{\traceof}[1]{\trace\qp{#1}}             %

    \providecommand{\transpose}{\intercal}%

    \providecommand{\transposeof}[1]{\ensuremath{\qp{#1}^\transpose}}
    \providecommand{\transposed}{{}^\transpose}

    \providecommand{\orthogonalto}[1][]{\ensuremath{\perp\ifx|#1|{}\else{_{#1}}\fi}}
    
    \providecommand{\rowof}[1]{\ensuremath{\vecof{#1}}}

    \providecommand{\discolvec}[2][r]{\ensuremath{\begin{bmatrix}[#1]#2\end{bmatrix}}}

    \providecommand{\discolvecthree}[4][r]{\ensuremath{\discolvec[#1]{#2\\#3\\#4}}}

    \providecommand{\discolvecithree}[1]{\discolvecthree{\vecentry{#1}1}{\vecentry{#1}2}{\vecentry{#1}3}}

    \providecommand{\disrowof}[1]{\rowof{\begin{matrix}[r]#1\end{matrix}}}
    \providecommand{\disrowvec}[1]{\disrowof{#1}}
    \providecommand{\zeroentry}{\phantom0}%
    
    \providecommand{\disrowvectwo}[2]{\ensuremath{\disrowvec{#1 &#2}}}

    \providecommand{\dismatrix}[2][r]{\ensuremath{\begin{bmatrix}[#1]#2\end{bmatrix}}}
    \providecommand{\dismattwo}[5][r]{\dismatrix[#1]{#2&#3\\ #4&#5}}

    \providecommand{\dismatthree}[9]{\dismatrix[r]{#1&#2&#3\\ #4&#5&#6\\ #7&#8&#9}}%

    \providecommand{\smint}{\ensuremath{{\text{\textbf{/}}}\kern-.75em\smallint}}
    \renewcommand{\smint}[1][]{\lower12.3pt\hbox{\begin{tikzpicture}\draw[line width=.75pt] (-3pt,-0.5)--(1pt,-0.5) node[pos=0.6]{$\int$};\path (3pt,-24pt)node {\scriptsize $#1$};\end{tikzpicture}}}

    \providecommand{\lap}{\ensuremath{\Delta}}
    \providecommand{\lapin}[1][]{\lap\ifx|#1|\else_{#1}\fi}
    \providecommand{\normalsymbol}{\operatorname{n}}
    \providecommand{\normal}[1][]{\normalsymbol\ifx|#1|\else_{#1}\fi}%
    \providecommand{\normalto}[2][]{\ensuremath{\normal[#2]\ifx|#1|\else\qp{#1}\fi}}
    \providecommand{\normalder}[1][]{\ensuremath{\normal\ifx|#1|\else\qp{#1}\fi{\inner\grad}}}
    \providecommand{\normalderto}[2][]{\ensuremath{\normalto[#1]{#2}{\inner\grad}}}

    \providecommand{\tangentialsymbol}{\operatorname{t}}
    \providecommand{\tangentialto}[2][]{\tangentialsymbol\ifx|#1|\else^{#1}\fi\ifx|#2|\else_{#2}\fi}
    \providecommand{\explicittangentialto}[2][]{\restrictionqf{\vec{#2}-\normalto{#1}\normalto{#1}\inner{\vec{#2}}}{\boundary{#1}}}

    \providecommand{\intersected}{\ensuremath{\cap}}
    
    \providecommand{\meet}{\intersected}

    \providecommand{\union}[1]{\ensuremath{\bigcup\nolimits_{#1}}}

    \providecommand{\unions}[3][]{\union{#2\in{#3}\ifx|#1|\else:#1\fi}}

    \ifthenelse{\boolean{usesvjour}}{}{
      \let\vec\undefined
      \providecommand{\vec}[1]{\ensuremath{\boldsymbol{#1}}}
      \renewcommand{\vec}[1]{\ensuremath{\boldsymbol{#1}}}
    }

    \providecommand{\hatmat}[1]{\hat{\mat{#1}}}

    \providecommand{\geomat}[1]{\vec{\MakeUppercase{#1}}}

    \providecommand{\mat}[1]{\geomat{#1}} %
    \providecommand{\Prob}[1][]{\ensuremath{\operatorname{Prob}\ifx|#1|{}\else_{#1}\fi}}

    \providecommand{\pdf}[2][]{\ensuremath{\operatorname{pdf}_{#2\ifx|#1|{}\else{\conditionalto{#1}}\fi}}\xspace}

    \providecommand{\expectation}{\ensuremath{\operatorname{E}}}
    \providecommand{\EX}[1][]{\ensuremath{\expectation\ifx|#1|{}\else_{#1}\fi}}

    \providecommand{\gausskernel}[3][x]{%
      \ensuremath{
        \exp\frac{-\if#20{#1}\else(#1-\mu)\fi^2}{%
          2\if#31{}\else\powp2{#3}\fi}%
      }%
    }
    \providecommand{\gaussdistribution}[3][x]{%
      \ensuremath{\frac1{\sqrt{2\pic}\if#31{}\else#3\fi}%
        \gausskernel[#1]{#2}{#3}
      }%
    }%

    \providecommand{\boundary}{\partial}
    \providecommand{\SPD}{\operatorname{SPD}}
    \providecommand{\spdmats}[2][F]{\SPD(\csname#1\endcsname{#2})}
     \providecommand{\Continuous}{\ensuremath{\operatorname C}\xspace}%
     \providecommand{\Hspace}{\ensuremath{\operatorname H}\xspace}
     \providecommand{\Lebesgue}{\ensuremath{\operatorname L}\xspace}
     \providecommand{\Besovspace}{\ensuremath{\operatorname B}\xspace}

     \providecommand{\Weaklyder}{\ensuremath{\operatorname W}\xspace}
     
     \providecommand{\dual}[1]{\ensuremath{{#1}'}}
     
     \providecommand{\dualspace}[2][]{\dual{\linspace{#2}\ifx|#1|\else{_{#1}}\fi}}
     \providecommand{\bidual}[1]{\ensuremath{{#1}''}}
     \providecommand{\bidualspace}[2][]{\bidual{\linspace{#2}\ifx|#1|\else{_{#1}}\fi}}

     \providecommand{\cont}[1]{\ensuremath{\Continuous^{#1}}}

     \providecommand{\BV}[1]{\ensuremath{\operatorname{BV}}}

     \providecommand{\leb}[1]{\ensuremath{\Lebesgue_{#1}}}
     \providecommand{\lebloc}[1]{\ensuremath{{{\Lebesgue}^{\kern-.20em\lower .1ex\hbox{\tiny\textrm{\textup{loc}}}}_{#1}}}}
     \providecommand{\lebnorm}[3][]{\ensuremath{\Norm{#2}_{\leb{#3}\ifx|#1|{}\else(#1)\fi}}}
     \providecommand{\bes}[3][]{\ensuremath{\Besovspace^{#2}_{#3\ifx|#1|\else,#1\fi}}}

     \providecommand{\sob}[2]{\ensuremath{{\smash\Weaklyder}^{#1}_{#2}}}

     \providecommand{\sobh}[1]{\ensuremath{\Hspace^{#1}}}
     \providecommand{\vecsobh}[1]{\ensuremath{\vec\Hspace^{#1}}}
     \providecommand{\hdiv}[1][]{\vecsobh{\operatorname{div}}\ifx|#1|\else(#1)\fi}
     \providecommand{\hcurl}[1][]{\vecsobh{\operatorname{curl}}\ifx|#1|\else(#1)\fi}
     
     \providecommand{\sobhz}[2][]{\sobh{#2}_{0\ifx+#1+\else|#1\fi}}

     \providecommand{\Lip}[1][]{\ensuremath{\operatorname{Lip}}\ifx|#1|{}\else{\qp{#1}}\fi}

     \providecommand{\holder}[2]{\cont{#1,#2}}
     \providecommand{\poly}[1]{\ensuremath{\rP}^{#1}}

     \providecommand{\Symmatrices}[2][R]{\ensuremath{\operatorname{Sym}{(\csname#1\endcsname{#2})}}}
     
     \providecommand{\SAmatrices}[2][F]{\ensuremath{\operatorname{SA}{(\csname#1\endcsname{#2})}}}
     \providecommand{\mesh}[2][]{{\ensuremath{\mathcalbf{\MakeUppercase{#2}}\ifx|#1|\else_{#1}\fi}}}

    \providecommand{\crouzeixraviart}[1][1]{\operatorname{CR}\ifx|#1|{}\else{^{#1}}\fi}

    \providecommand{\linspace}[1]{\mathscript{\MakeUppercase{#1}}}
    \providecommand{\fatlinspace}[1]{\mathcalbf{\MakeUppercase{#1}}}
    
    \providecommand{\linop}[1]{\mathcalbf{\MakeUppercase{#1}}}

    \providecommand{\clinopss}[2]{\clinopss{\linspace{#1}}{\linspace{#2}}}
    \providecommand{\fepartition}[2][]{\mathscript{\MakeUppercase{#2}}\ifx|#1|{}\else_{#1}\fi}
    \providecommand{\fespace}[2][]{\mathbb{\uppercase{#2}}\ifx|#1|{}\else_{#1}\fi}
    
    \providecommand{\vespace}[1][]{\fespace v\ifx|#1|\else_{#1}\fi}

    \providecommand{\fe}[2][]{\ensuremath{\uppercase{#2}\ifx|#1|\else_{#1}\fi}}%

    \providecommand{\fespacefun}[2]{\fe[\fespace{#1}]{#2}}
    \providecommand{\vecfe}[2][]{\ensuremath{\vec{\uppercase{#2}}\ifx|#1|{}\else{_{#1}}\fi}}%
    \providecommand{\vecfespacefun}[2]{\vecfe[\fespace{#1}]{#2}}
    \providecommand{\matfe}[2][]{\ensuremath{\mat{\uppercase{#2}}\ifx|#1|{}\else{_{#1}}\fi}}%
    \providecommand{\matfespacefun}[2]{\matfe[\fespace{#1}]{#2}}
    \providecommand{\hatmatfe}[2][]{\ensuremath{\hatmat{\uppercase{#2}}\ifx|#1|{}\else{_{#1}}\fi}}%
    
    \providecommand{\EOC}{\ensuremath{\operatorname{EOC}}\xspace}
    \providecommand{\tol}{\ensuremath{\operatorname{tol}}\xspace}

    \providecommand{\Foreach}{\text{ for each }}%
    \providecommand{\Forsome}{\text{ for some }}

    \RequirePackage{amscd}

    \providecommand{\funk}[3]{\ensuremath{#1:#2\to#3}}

    \providecommand{\dfunkmapsto}[6][]{\ensuremath{
        \begin{array}{rrcl}
          {#2}: & {#4} &  \to   & {#6}
          \\
                & {#3} &\mapsto & {#5\text{\ #1}}
        \end{array}\quad}}

    \providecommand{\restriction}[2]{\left.#1\right|_{#2}}
    \renewcommand{\restriction}[2]{\left.#1\right|_{#2}}

    \providecommand{\restrictionqf}[2]{\qb{#1}_{#2}}

    \providecommand{\evalat}[3][]{\qb{#2}_{\ifx|#1|{}\else#1=\fi#3}}
    \providecommand{\evaldiff}[4][]{\qb{#2}^{\ifx|#1|{}\else#1=\fi#3}_{\ifx|#1|{}\else#1=\fi#4}}
    \providecommand{\boundarytraceof}[2]{\restriction{#2}{\boundary{#1}}}

    \providecommand{\aka}[1]{(also known as {#1})\xspace}
    
    \providecommand{\akaindexemph}[2][]{\aka{\indexemph[#1]{#2}}}
    \providecommand{\CBS}{Cauchy--Bunyakovsky--Schwarz inequality\xspace}
    \providecommand\bs{\char '134}   %
    \providecommand\caret^%

    \providecommand{\codelinecolor}{a}
    \providecommand{\codeline}[2][]{%
      \begin{center}
        \begin{minipage}{\linewidth}
          \begin{tikzpicture}
            \node[fill=\codelinecolor!12.5!g]{
              \tiny
              \begin{minipage}{\linewidth}
                \nolinkurl{#1}\ \nolinkurl{#2}%
              \end{minipage}
            };
          \end{tikzpicture}
        \end{minipage}
      \end{center}
    }

    \providecommand{\texcommand}[1]{\texttt{\bs{\nolinkurl{#1}}}\xspace}
    
    \providecommand{\codename}[1]{\nolinkurl{#1}\xspace}
    \providecommand{\colorvarname}[2][a]{\colorbox{#1!6.25}{\nolinkurl{#2}}}%
    \providecommand{\codevarname}[1]{\colorvarname[a]{#1}}

    \ifthenelse{\boolean{useutopia}}{%
      
    }{%
      
    }

    \providecommand{\codeprint}[2][.]{
      \par
      \begin{center}
        \framebox{Printout of file %
          \ifthenelse{\isundefined\pickuppath}{%
            \codename{#2}
          }{%
            \providecommand{\fullpickuppath}{}
            \renewcommand{\fullpickuppath}{\pickuppath/\ifx|#1|\else#1/\fi#2}
            \href{\fullpickuppath}{\codename{#2}}
        }}
        \\
        \lstinputlisting{#1/#2}
      \end{center}
      \par
    }

    \providecommand{\indexen}[2][]{{\ifthenelse{\boolean{shownotes}}{\color b}{}#2\ifx|#1|\index{#2}\else\index{#1}\fi}}
    \providecommand{\indexemph}[2][]{\emph{\indexen[#1]{#2}}}

    \providecommand{\ListParameters}{}
    \renewcommand{\ListParameters}%
    {
    	 \setlength{\topsep}{0pt}
    	 \setlength{\leftmargin}{0pt}
             \setlength{\itemsep}{0pt}
    	 \setlength{\parsep}{0pt}
    	 \setlength{\parskip}{0pt}
             \setlength{\labelsep}{0pt}
    	 \setlength{\itemindent}{0pt}
    }
    {%
      \begin{list}%
        {}%
        {\ListParameters%
        
    }}%
    {\end{list}}
    \newcounter{tmpcounter}
    \newcounter{LetterListItem}
    \renewcommand{\theLetterListItem}{(\alph{LetterListItem})}

    \newcounter{CapitalListItem}
    \renewcommand{\theCapitalListItem}{\Alph{CapitalListItem}.}

    \newcounter{NumberListItem}
    \renewcommand{\theNumberListItem}{\arabic{NumberListItem}}
    {
    	\begin{list}%
    	{\theNumberListItem.\ }%
    	{\usecounter{NumberListItem}%
    	 \ListParameters
    	}
    }%
    {\end{list}}
    \newcounter{QuestionListItem}
    \renewcommand{\theQuestionListItem}{\textbf{Question \arabic{QuestionListItem}}}
    {
    	\begin{list}%
    	{\theQuestionListItem.\ }%
    	{\usecounter{QuestionListItem}%
    	 \ListParameters
    	}
    }%
    {\end{list}}
    \newcounter{RomanListItem}
    \renewcommand{\theRomanListItem}{(\roman{RomanListItem})}
    {
    	\begin{list}%
    	{\theRomanListItem\ }%
    	{\usecounter{RomanListItem}
    	 \ListParameters
    	}
    }%
    {\end{list}}
    \newcounter{StepsItem}
    {
    	\begin{list}%
    	{Step \theStepsItem.\ }%
    	{\usecounter{StepsItem}%
    	 \ListParameters
    	}
    }%
    {\end{list}}
    \newcounter{CasesListItem}
    \renewcommand{\theCasesListItem}{\Alph{CasesListItem}}
    \newenvironment{CasesList}%
    {
    	\begin{list}%
    	{\emph{Case \theCasesListItem.}\ }%
    	{\usecounter{CasesListItem}%
    	 \ListParameters
    	}
    }%
    {\end{list}}
    \newcounter{QAListItem}
    \renewcommand{\theQAListItem}{Q\arabic{QAListItem}:}
    {
    	\begin{list}%
    	{\theQAListItem}%
    	{\usecounter{QAListItem}
    	 \ListParameters
    	}
    }%
    {\end{list}}

    \ifthenelse{\boolean{isthesis}}{%
      \setcounter{secnumdepth}{1}%
    }{
      \setcounter{secnumdepth}{2} %
    }
    \providecommand{\ListParameters}{}
    \renewcommand{\ListParameters}
    {
    	 \setlength{\topsep}{0em}
    	 \setlength{\leftmargin}{0em}
             \setlength{\itemsep}{0ex}
    	 \setlength{\parsep}{.5ex}
    	 \setlength{\itemindent}{\labelsep}
    	 \addtolength{\itemindent}{\labelwidth}
    }

      \providecommand{\ObsName}{Remark}%
      \providecommand{\RemName}{Remark}%
      \providecommand{\NotName}{Notation}%
      \providecommand{\BFNName}{Big~Fat~Note}%
      \providecommand{\DefName}{Definition}%
      \providecommand{\ExaName}{Example}%
      \providecommand{\TheName}{Theorem}%
      \providecommand{\LemName}{Lemma}%
      \providecommand{\ProName}{Proposition}%
      \providecommand{\CorName}{Corollary}%
      \providecommand{\PbmName}{Problem}%
      \providecommand{\HypName}{Hypothesis}%
      \providecommand{\AlgName}{Algorithm}%
      \providecommand{\ExeName}{Exercise}%
      \providecommand{\SolName}{Solution}%
      \providecommand{\ClaName}{Claim}%
      \providecommand{\EsyName}{Essay}%
      \providecommand{\Proofname}{Proof}%
      \providecommand{\Derivename}{Derivation}%
    
    \ifthenelse{\boolean{isthesis}}{%
      \providecommand{\Thecounter}{The}
    }{%
      \providecommand{\Thecounter}{subsection}
    }
    \makeatletter
    \newcommand{\oltikzgetxy}[3]{%
      \tikz@scan@one@point\pgfutil@firstofone#1\relax
      \edef#2{\the\pgf@x}%
      \edef#3{\the\pgf@y}%
    }
    \makeatother
    \providecommand{\pdfformat}[1]{
       \provideboolean{pdfoutput}
       \setboolean{pdfoutput}{#1}%
      \ifthenelse{\boolean{pdfoutput}}{
        \typeout{using pdf}
\usepackage{pdfsync}
        \providecommand{\graphext}{pdf}
        \renewcommand{\graphext}{pdf}
        \providecommand{\graphextex}{pdf_t}
        \renewcommand{\graphextex}{pdf_t}
      }{
        \typeout{using eps}
        \RequirePackage[dvips]{graphicx,xcolor}
        \providecommand{\graphext}{eps}
        \renewcommand{\graphext}{eps}
        \providecommand{\graphextex}{eps_t}
        \renewcommand{\graphextex}{eps_t}
      }
      \RequirePackage{epsfig}
      \RequirePackage{tikz}
      \RequirePackage{rotating}
      \RequirePackage{graphicx}
      \RequirePackage{xcolor}
      \provideboolean{darkcolortheme}
      \definecolor{SussexFlint}{rgb}{.00,.19,.21}
      \definecolor{SussexGrey}{rgb}{.51,.58,.49}
      \definecolor{SussexOrange}{rgb}{.94,.29,.00}
      \definecolor{SussexYellow}{rgb}{1.00,.73,.00}
      \definecolor{SussexRed}{rgb}{.94,.01,.49}
      \definecolor{SussexPurple}{rgb}{.48,.06,.44}
      \definecolor{SussexGreen}{rgb}{.00,.58,.46}
      \definecolor{OmarGreen}{rgb}{.00,.68,.36}
      \definecolor{SussexBlue}{rgb}{.00,.58,.65}
      \definecolor{OmarBlue}{rgb}{.00,.38,.65}
      \colorlet{a}{OmarBlue}%
      \colorlet{b}{SussexOrange}
      \colorlet{c}{SussexGreen}
      \colorlet{d}{SussexPurple}%
      \colorlet{e}{SussexRed}
      \colorlet{f}{SussexYellow}
      \colorlet{g}{white}%
      \colorlet{h}{SussexGrey}%
      \colorlet{i}{black}%
      \colorlet{j}{SussexFlint}
      \colorlet{colora}{a}
      \colorlet{colorb}{b}
      \colorlet{colorc}{c}
      \colorlet{colord}{d}
      \colorlet{colore}{e}
      \colorlet{colorf}{f}
      \colorlet{colorg}{g}
      \colorlet{colorh}{h}
      \colorlet{colori}{i}
      \colorlet{colorj}{j}
      \newcommand{\mausDarkColorTheme}{
        \colorlet{a}{SussexYellow!50!yellow}
        \colorlet{b}{SussexBlue}%
        \colorlet{c}{SussexRed!50!red}
        \colorlet{d}{SussexOrange!50!yellow}
        \colorlet{e}{SussexGreen!50!green}
        \colorlet{f}{SussexPurple!50!magenta}
        \colorlet{g}{black}%
        \colorlet{h}{SussexFlint!50!black}
        \colorlet{i}{white}%
        \colorlet{j}{SussexGrey}
      }
      \ifthenelse{\boolean{darkcolortheme}}{\mausDarkColorTheme}{}
    }
    \providecommand{\solution}{\textbf{\SolName.}\xspace}

     \newcounter{phantombox}[enumi]%
     \provideboolean{showphantoms}
     \renewcommand{\thephantombox}{\Alph{phantombox}}%
     \newcommand{\phantombox}[1]{\stepcounter{phantombox}%
       \ensuremath{\boxed{%
           {\ifthenelse{\boolean{showphantoms}}{#1}{\phantom{#1}}}%
           {\texttt{\tiny\ \colorbox{i!50}{\color g\thephantombox}}
           }%
         }%
       }%
     }

     \provideboolean{hidesolution}
     \newcommand{\consolution}[2][]{
       \ifthenelse{\boolean{hidesolution}}{#1\setboolean{showphantoms}{false}}{%
         {\setboolean{showphantoms}{true}\color{i!50}\par \small {\solution}\ #2\par\ \\[5pt]}}
     }
     \provideboolean{showmarks}
     \providecommand{\showmarks}[1]{%
       \ifthenelse{%
         \boolean{showmarks}}{%
         \marginpar{%
           \tiny [$#1$ mark\ifthenelse{\equal{#1}1}{\phantom{s}}s]}%
       }{}}%

     \newcommand{\condibreak}{\ifthenelse{\boolean{hidesolution}}{\newpage}{}}

     \providecommand{\qeyword}[1]{\index{#1}\ifthenelse{\boolean{shownotes}}{{\tiny\color e\colorbox{e!6.25}{#1}}}{}}
     \providecommand{\pathword}[1]{\ifthenelse{\boolean{shownotes}}{\ \\\index{#1@\tiny\codevarname{#1}}{\tiny\color f\colorvarname[f]{#1}}}{}}
     \providecommand{\sourcecite}[2][]{\ifthenelse{\boolean{shownotes}}{{\ \\\tiny\colorbox{d!6.25}{\color d\texttt{source: \citet[#1]{#2}}}}}{\nocite{#2}}}
     \providecommand{\conword}[1]{\ifthenelse{\boolean{shownotes}}{#1}{}}

     \RequirePackage{hyphenat}
    \RequirePackage{lineno}
    \ifthenelse{\boolean{showchanges}}{
      \newcommand{\llabel}[1]{\hypertarget{llineno:#1}{\linelabel{#1}}}
      \newcommand{\lref}[1]{\hyperlink{llineno:#1}{\ref*{#1}}}
    }{
      \newcommand\llabel[1]{}
      \newcommand\lref[1]{}
    }
    \provideboolean{includeresponses}
    \setboolean{includeresponses}{false}
    \providecommand{\mailto}[1]{\href{mailto:#1}{\nolinkurl{#1}}}

\newtheoremstyle{plain}%
  {}%
  {}%
  {\mdseries\slshape}%
  {\parindent}%
  {\bfseries}%
  {.}%
  {.5em}%
  {}%

\newtheoremstyle{note}%
  {}%
  {}%
  {}%
  {\parindent}%
  {\bfseries}%
  {.}%
  {.5em}%
  {}%

\newtheoremstyle{claim}%
  {}%
  {}%
  {\mdseries\slshape}%
  {}%
  {\bfseries}%
  {}%
  {.5em}%
  {}%

\newtheoremstyle{exercise}%
  {}%
  {}%
  {}%
  {}%
  {\bfseries}%
  {.}%
  {1em}%
  {}%

\newtheoremstyle{break}%
  {}%
  {}%
  {}%
  {}%
  {\bfseries}%
  {.}%
  {\newline}%
  {}%

\swapnumbers{
  \theoremstyle{plain}
  \ifthenelse
      {\boolean{isthesis}}
      {\newtheorem{The}{\TheName}[section]}%
      {
      }%
{
   \theoremstyle{plain}

   \renewcommand{\Thecounter}{subsection}

   \newtheorem*{The*}{\TheName}
   \newtheorem*{Lem*}{\LemName}
   \newtheorem*{Pro*}{\ProName}
   \newtheorem*{Cor*}{\CorName}
   \newtheorem*{Pbm*}{\PbmName}
   \newtheorem*{Hyp*}{\HypName}
   \newtheorem*{Exe*}{\ExeName}
   \newtheorem*{Txx*}{\ExeName} %
   \newtheorem*{Con*}{Conclusion}
   \newtheorem*{Sum*}{Summary}
 }
 {
   \theoremstyle{claim}

 }
 {
   \theoremstyle{note}

   \newtheorem*{Obs*}{\ObsName}

   \newtheorem*{Def*}{\DefName}
   \newtheorem*{Exa*}{\ExaName}
   \newtheorem*{Alg*}{\AlgName}
 }

 {
   \theoremstyle{break}
 }
}

\newenvironment{The}[1][]{%
  \ifx&#1&%
  \subsection{\TheName\xspace}%
  \else%
  \subsection[#1 theorem]{\TheName\ (#1)}%
  \fi%
  \slshape}{%
  \upshape}

\newenvironment{Lem}[1][]{\subsection{\LemName\xspace{\ifx&#1&{}\else{ (#1)}\fi}}\slshape}{\upshape}

\newenvironment{Def}[1][]{\subsection{\DefName\xspace{\ifx&#1&{}\else{ of #1}\fi}}}{}
\newenvironment{Obs}[1][]{\subsection{\ObsName\xspace{\ifx&#1&{}\else{ (#1)}\fi}}}{}

\newenvironment{Alg}[1][]{\subsection{\AlgName\xspace{\ifx&#1&{}\else{ (#1)}\fi}}}{}
\providecommand{\qed}{\vrule height 5pt depth 0pt width 3pt}
\providecommand{\qqed}{{\raggedright{\ \hfill\qed}}}

\newcounter{passo}

\newenvironment{Proof}[1][]%
{\par\noindent{\bf \Proofname\ifx|#1|.\ \else\ #1.\ \fi}\setcounter{passo}{0}}%
{\qqed\par}
{\par\noindent{\bf \Derivename\ #1}\setcounter{passo}{0}}%
{\qqed\par}
\newenvironment{Proof*}[1][{}]%
{\subsection{\Proofname\ #1}\setcounter{passo}{0}}
{\qqed\par}

\newboolean{useimanum}
\makeatletter%
\@ifclassloaded{imanum}{
  \setboolean{useimanum}{true}
}{
  \setboolean{useimanum}{false}
}%
\makeatother%
\setboolean{showchanges}{true}%
\setboolean{shownotes}{false}%
\setboolean{includeresponses}{false}%
\makeindex
\usepackage{natbib}
\usepackage{subfig}
\usepackage{graphics}
\usepackage{algorithm}%
\usepackage{algcompatible}
\algnewcommand\algorithmicprocedure{\textbf{procedure}}
\algnewcommand\PROCEDURE{\item[\algorithmicprocedure]}%
\algnewcommand\algorithmicendprocedure{\textbf{end procedure}}
\algnewcommand\ENDPROCEDURE{\item[\algorithmicendprocedure]}%
\renewcommand{\leq}{\leqslant}
\renewcommand{\geq}{\geqslant}
\renewcommand{\rot}{\nabla\!\times\!}%
\providecommand{\codevarname}[1]{\colorvarname{#1}}
\renewcommand{\codevarname}[1]{\colorvarname{#1}}
\providecommand{\tol}{\codevarname{tol}}
\renewcommand{\tol}{\codevarname{tol}}
\providecommand{\maxiter}{\codevarname{maxiter}}
\renewcommand{\linop}[1]{\mathcal{\MakeUppercase{#1}}}
\providecommand{\linoptheta}[1][\theta]{\linop M_{#1}}
\providecommand{\abil}[3][a]{#1\qp{#2\,;\,#3}}
\providecommand{\aqua}[2][a]{\abil[#1]{#2}{#2}}
\ifthenelse{\boolean{useutopia}}{
  \providecommand{\fatlinspace}[1]{\mathcalbf{\MakeUppercase{#1}}}
}{
  \providecommand{\fatlinspace}[1]{\mathfrak{\MakeUppercase{#1}}}
}
\providecommand{\feop}[2][h]{\mathcal{\MakeUppercase{#2}}\ifx|#1|\else_{#1}\fi}
\providecommand{\pointinter}[1]{\feop[\fespace{#1}]I}
\providecommand{\pinteron}[1]{\pointinter{\fespace{#1}}}
\ifthenelse{\boolean{useutopia}}{
  \renewcommand{\fe}[2][]{\ensuremath{\mathsfit{#2}\ifx|#1|\else_{#1}\fi}}%
  \renewcommand{\vecfe}[2][]{\ensuremath{\vec{\mathsfit{#2}}\ifx|#1|{}\else{_{#1}}\fi}}%
  \renewcommand{\matfe}[2][]{\ensuremath{\mat{\mathsfit{\MakeUppercase{#2}}}\ifx|#1|{}\else{_{#1}}\fi}}%
}{
  \renewcommand{\fe}[2][]{\ensuremath{\mathsf{#2}\ifx|#1|\else_{#1}\fi}}%
  \renewcommand{\vecfe}[2][]{\ensuremath{\vec{\mathsf{#2}}\ifx|#1|{}\else{_{#1}}\fi}}%
  \renewcommand{\matfe}[2][]{\ensuremath{\mat{\mathsf{\MakeUppercase{#2}}}\ifx|#1|{}\else{_{#1}}\fi}}%
}
\providecommand{\aposteriori}{a posteriori\xspace}

\providecommand{\apriori}{{a priori}\xspace}

\renewcommand{\aposteriori}{a posteriori\xspace}

\renewcommand{\apriori}{{a priori}\xspace}

\providecommand{\mesht}{{\mesh t}}

\provideboolean{includeresponses}
\setboolean{includeresponses}{false}
\renewcommand{\colorvarname}[2][i]{\texttt{#2}}%
\numberwithin{equation}{section}

\providecommand{\ourtitle}{A least-squares Galerkin approach to gradient and Hessian recovery for nondivergence-form elliptic equations}
\providecommand{\ourshorttitle}{Least squares Galerkin for nondivergence elliptic equations}

\providecommand{\authoramireh}{Amireh Mousavi\xspace}
\providecommand{\aaddressamireh}{Isfahan University of Technology\xspace}
\providecommand{\baddressamireh}{Isfahan}
\providecommand{\caddressamireh}{Iran}
\providecommand{\emailamireh}{amireh.mousavi@math.iut.ac.ir}

\providecommand{\authoromar}{Omar Lakkis\xspace}
\providecommand{\aaddressomar}{University of Sussex\xspace}
\providecommand{\baddressomar}{Brighton\xspace}
\providecommand{\caddressomar}{England UK\xspace}
\providecommand{\emailomar}{lakkis.o.maths@gmail.com}
\providecommand{\ourthanks}{
    This work was supported by the ModCompShock Marie Skłodowska--Curie
    International Training Network and was possible thanks to DISIM of
    the University of L'Aquila, Italy, where most of the reported research
    took place in 2017--18.}
\providecommand{\ourabstract}{%
  We propose a least-squares method involving the recovery of the
  gradient and possibly the Hessian for elliptic equation in
  nondivergence form. As our approach is based on the Lax--Milgram
  theorem with the curl-free constraint built into the target (or
  cost) functional, the discrete spaces require no inf-sup
  stabilization. We show that standard conforming finite elements can
  be used yielding \apriori and \aposteriori convergence results. We
  illustrate our findings with numerical experiments with uniform or
  adaptive mesh
  refinement.\xspace%
}
\providecommand{\ourkeywords}{nondivergence form, elliptic equations, finite element method, convergence}
\ifthenelse{\boolean{useimanum}}{
  \AtBeginDocument{
    \label{CorrectFirstPageLabel}
    \def\fpage{\pageref{CorrectFirstPageLabel}}
  }
  \renewcommand{\vec}[1]{\ensuremath{\mathbf{#1}}}
  \received{1 September 2019}
  \revised{6 April 2021}
}{}%
\begin{document}
\ifthenelse{\boolean{includeresponses}}{
  \renewcommand{\thepage}{Responses to the Referees (\roman{page})}
  \clearpage
  \renewcommand{\thepage}{\arabic{page}}
  \setcounter{page}1
}{}
\ifthenelse{\boolean{useimanum}}{
  \title{\ourtitle}
  \shorttitle{\ourshorttitle}
  \author{
    \textsc{
      \authoromar\thanks{Corresponding author. Email: \email{\emailomar}
      \aaddressomar, \baddressomar, \caddressomar}
    }
    \\[6pt]
    and
    \\[6pt]
    \textsc{
      \authoramireh\thanks{Email: \email{\emailamireh}
      \aaddressamireh, \baddressamireh, \caddressamireh}
    }
  }
  \shortauthorlist{O. Lakkis \& A. Mousavi}
}{
  \title[\ourshorttitle]{\ourtitle}
  \author{\authoromar}
  \address{%
    \authoromar,
    \aaddressomar,
    \baddressomar,
    \caddressomar
  }
  \email{\emailomar}
  \author{\authoramireh}
  \address{%
    \authoramireh,
    \aaddressamireh,
    \baddressamireh,
    \caddressamireh
  }
  \email{\emailamireh}
  \thanks{\ourthanks}
}
\maketitle
\ifthenelse{\boolean{useimanum}}{
  \ \\
  \begin{abstract}
    {\ourabstract}{\ourkeywords}
  \end{abstract}
}{%
  \begin{abstract}
    \ourabstract
  \end{abstract}
}
\section{Introduction}
Elliptic equations in nondivergence form play an important role in
many domains of pure and applied mathematics ranging from nonlinear
PDEs
\citep{CaffarelliCabre:95:book:Fully,ArmstrongSmart:10:article:An-easy}
to Probability Theory
\citep{Evans:85:article:Some,FabesStroock:83:article:The-Lp-intergrability},
continuum Game Theory,
homogenization~\citep{CapdeboscqSprekelerSuli:20:article:Finite} and
wave
propagation~\citep{ArjmandKreiss:17:techreport:An-Equation-Free}. The
numerical approximation of such equations (references to be given
below) %
plays thus an important %
role. Here we propose a least-squares based gradient- or
Hessian-recovery Galerkin finite element method for the numerical
approximating of a function $\funk{u}\W\reals$, $\W$ convex, solving
the following \emph{linear elliptic Dirichlet boundary value problem
  in nondivergence form}\index{$\linop L$}\index{$u$}
\begin{equation}
  \label{eq:nondivergence0-inhomogeneous}
  \linop Lu
  :=
  \mat A\frobinner\D^2 u
  +
  \vec b\inner \grad u
  -
  c u
  =
  f
  \tand
  \restriction u{\boundary\W}=r
\end{equation}
where $f \in \leb2(\W)$, $r \in \sobh{3/2}(\boundary\W)$, all
coefficients are measurable, $\mat A$ is a uniformly elliptic
tensor-valued,\index{$\lambda_\flat$}\index{$\lambda_\sharp$}
\begin{equation}
  \label{def:uniformly-elliptic}
  \lambda_\flat
  \eye
  \leq
  \mat A
  \leq
  \lambda_\sharp
  \eye,
  \text{\Ae in }\W
  ,
  \Forsome\lambda_\sharp\geq\lambda_\flat>0,
\end{equation}
$c$ is non-negative on $\W$ and $\mat A$, $\vec b$, $c$ satisfy either
of the \indexemph{Cordes condition} (\ref{def:general-Cordes-condition}) or
(\ref{def:special-Cordes-condition}) (to be discussed in
\S~\ref{sec:Cordes-conditions}). Roughly speaking, the Cordes condition
allows us to reformulate the operator $\linop L$ so that it is close
enough to an invertible operator in \indexemph{divergence
  form} thereby ensuring the
elliptic problem with discontinuous coefficients is well-posed
(see \S\ref{sec:Cordes-conditions} for more details).

A main difficulty in the study of elliptic PDEs in nondivergence
form is the lack of a natural variational structure which precludes
a straightforward use of weak solutions in $\sobh1(\W)$, say, and
their numerical approximation using the bilinear form given by the
exact problem. One is thus forced to find some suitable
approximation of the Hessian more or less directly. The appropriate
concept of generalized solution for nondivergence form equations is
that of viscosity solution, which relies on the maximum principle.
In this respect, finite difference methods have the advantage over
Galerkin methods, in that they replicate more easily the maximum
principle, which is very useful when aiming at the approximation of
viscosity solutions. On the other hand finite difference methods,
besides lacking the geometric flexibility and the higher order
approximation power of Galerkin methods, must be modified to take
into account coefficients that are more singular than Lipschitz
\citep{FroeseOberman:09:article:Numerical}. Dealing with the
boundary is also not that straightforward as with Galerkin methods
Which we deal with in this article.

Galerkin methods for general elliptic PDEs in nondivergence form were
studied by \citet{Bohmer:10:book:Numerical}, but $\cont1(\W)$ finite
elements are required for their practical
implementations~\citep{DavydovSaeed:13:article:Numerical}. A
\indexemph{recovered Hessian} finite element method for approximating
the solution of nondivergence form elliptic equation was introduced by
\citet{LakkisPryer:11:article:A-finite}; this method was later
generalized and fully analyzed by
\citet{Neilan:17:article:Convergence}. \indexemph{Discontinuous 
Galerkin} approaches have been proposed by
\citet{SmearsSuli:13:article:Discontinuous},
\citet{FengHenningsNeilan:17:article:Finite} and
\citet{FengNeilanSchnake:18:article:Interior}. Further Galerkin
approaches for nondivergence form equation do exist such as the
\indexemph{two-scale Galerkin method} which is based on an
integro-differential scheme by
\citet{NochettoZhang:18:article:Discrete} and the somewhat related
method of \citet{FengJensen:17:article:Convergent}, which draws on the
\indexemph{semi-Lagrangian methods} and the celebrated convergence
theorem of \citet{BarlesSouganidis:91:article:Convergence}, the
\indexemph{primal-dual weak Galerkin} method
\citet{WangWang:18:article:A-primal-dual} and the variational
formulation of elliptic problems in nondivergence form of
\citet{Gallistl:17:article:Variational}.

In this paper, we propose a least-squares approach combined with a
gradient and Hessian recovery. Our approach is related to the method
of \citet{SmearsSuli:13:article:Discontinuous} in that the test
function is the elliptic operator (or an approximation thereof)
applied to the ``variable function'', but, unlike them, we use
conforming finite elements. Our work is also connected to that of
\citet{Gallistl:17:article:Variational} with the key departure that
our least-squares approach allows a cost-functional enforcement of the
curl-free requirement rather than imposing this on the function space
and having to enforce it discretely via inf-sup stable
discretizations. Indeed, a feature of the method we will propose is
that it is coercive and based on the idea of gradient- or
Hessian-recovery combined with Lax--Milgram theorem, which as noted by
\citet{BrambleLazarovPasciak:97:article:A-least-squares} it is one of
the two main approaches of least squares Galerkin methods (the other
is a weighted-residual approach based on the Agmon--Douglis--Nirenberg
theory). An obviously non-exhaustive list of references to the
least-squares based Galerkin methods for linear and nonlinear we came
across is further complemented by
\citet{AzizKelloggStephens:85:article:Least} (based on the ADN theory)
\citet{BochevGunzburger:06:incollection:Least-squares} (which gives a
thorough survey at writing time)
\citet{DeanGlowinski:06:article:Numerical} (which uses least-squares
to solve the Monge--Ampère equation, related to nondivergence PDEs)
and its further refinement
\citet{CaboussatGlowinskiSorensen:13:article:A-least-squares}. We
deem it worth noting that an early attempt at least-squares FEMs for
elliptic equation in nondivergence form by
\citet{BrambleSchatz:70:article:Rayleigh-Ritz-Galerkin} is quite
inspiring, despite the difficulties in practical implementations of
methods there proposed (they require the same $\sobh2$-conformity as
\cite{Bohmer:10:book:Numerical}).

The rest of this article is structured as follows: in
\S~\ref{sec:least-squares-approach-to-EPDENDF} we introduce the main
background material, the cost (or energy) functional $E_\theta$ (where
\indexen{$\theta$} is a parameter) to be minimized and the associated
bilinear forms; we give some technical remarks. In
\S~\ref{sec:coercivity-continuity-of-cost-functional} we show that the
bilinear forms associated with $E_\theta$ satisfy the Lax--Milgram
theorem's assumptions thereby guaranteeing the least-squares problem
and the equivalent exact PDE are well-posed. In
\S~\ref{sec:a-conforming-Galerkin-FEM} we introduce the Galerkin
discretization, which, thanks to
\S~\ref{sec:coercivity-continuity-of-cost-functional}, enjoys
quasi-optimality and convergence properties on general finite element
spaces without the need to enforce inf-sup; we also derive via a
residual--error \aposteriori estimate, indicators and an adaptive
algorithm. Finally in \S~\ref{sec:numerical-experiments} we
illustrate the theoretical findings with numerical experiments in both
uniform and adaptive mesh refinement frameworks, before giving some
conclusions and outlook in \S~\ref{sec:conclusion}.
\section{Least-squares approach to elliptic problems in nondivergence form}
\label{sec:least-squares-approach-to-EPDENDF}
We now provide the main technical ideas for our approach. After some
preliminaries, function spaces in
\S~\ref{sec:basic-notation-and-function-spaces}, we discuss the Cordes
conditions in \S~\ref{sec:Cordes-conditions} and the nonhomogenous
Dirichlet problem in \ref{sec:Dirichlet-boundary-conditions}. We
introduce in \ref{sec:least-squares-problem} the least-squares
formulation with cost (or energy) functional $E_\theta$ of
problem~(\ref{eq:nondivergence0-inhomogeneous}) with $r=0$ and show
the equivalence between solving this and the Euler--Lagrange equations
in \S~\ref{sec:Euler-Lagrange-equations} and briefly discussing a
Hessian-less variant of our method in
Remark~\ref{obs:Hessian-less-approach}. We close this section by
introducing further the bilinear forms in \S~\ref{sec:bilinear-forms}
and recalling a useful Maxwell-type estimate of
\citet{CostabelDauge:99:article:Maxwell} in
Lemma~\ref{lem:Maxwell-estimate}.
\subsection{Basic notation and function spaces}
\label{sec:basic-notation-and-function-spaces}
For two vectors
$\vec{x}=\seqidotsfromto{x}1m,\vec{y}=\seqidotsfromto{y}1m$ (displayed
as columns with row transposes) in $\R m$ we write
$\vec{x}\inner\vec{y}:=\vec{x}\transposed\vec{y}:=\sumsifromto{x}i1m\vecentry{y}i$.
For a matrix, $\mat M\in\realmats mm$, $\trace\mat M$ denotes the
trace of the matrix $\mat M$, defined as the sum of its eigenvalues
(or, equivalently, its diagonal entries) and $\det\mat M$ denotes the
determinant of $\mat M$ defined as the product of its eigenvalues.
For two matrices $\mat M,\mat N\in\realmats ml$, their
\indexemph{Frobenius inner product} is defined by $\mat
M\frobinner\mat N:=\traceof{\mat{M}\transposed\mat{N}}$ and by $\norm{
  \mat M }$ we mean the \indexemph{Frobenius norm} of the matrix $\mat
M$, defined as \(\norm{\mat{M}}:=\qpsqrt{\mat M \frobinner\mat M}\),
which coincides with the Euclidean norm of $\mat M$'s spectrum.

Throughout the paper, including the above we denote, for a function
(or distribution) $\funk{\vec\phi}\W{\R m}$, $m\in\naturals$, by
$\D\vec\phi$ its first \indexemph{derivative},
$\grad\vec\phi:=\transposeof{\D\vec\phi}$ its \indexemph{gradient}
and, when $m=1$ (with a slight abuse of notation) by $\D^2\vec \phi$
its \indexemph{Hessian} (matrix or tensor). We shall also denote the
\indexemph{divergence} by $\div$, the \indexemph{curl}
\akaindexemph{rotation} by $\rot$ and the \indexemph{Laplace operator}
by $\lap:=\div\!\grad$.
The smallest and largest of two numbers $a,b$ 
are respectively denoted $\mini ab$ and $\maxi ab$.

We help the reader interested in tracking constants by
labeling them in accordance to the display where they are defined or
first appear; to lighten notation their dependence on other constants
or parameters is silent outside the definition, except when strictly
necessary, e.g., the parameters are variables in the given context.
For example, defining
\begin{equation}
  \constdef[\alpha,\beta]{const:example}:=\frac{\maxi\alpha\beta}{\beta},
\end{equation}
would be used as follows: $X\leq\constref{const:example}Y$ for each fixed $\alpha,\beta$,
or $B(\beta)\leq\sum_\alpha\constref[\alpha]{const:example}A(\alpha,\beta)$
for each fixed $\beta$ (but variable $\alpha$).

Consider a real number $p\geq1$ and a non-negative integer $s\in\NO$,
given a normed vector space $\pair X{\norm\cdot}$,
denote by {$\sob sp(\W;X)$}\index{$\sob sp(\W;X)$} the
\indexemph{Sobolev space} of $X$-valued functions $f$ in
\indexen{$\leb p(\W;X)$} whose (generalized/distributional/weak)
derivatives up to order {$s$} are in $\leb p(\W;Y)$ (for the
appropriate $Y$); $\leb p(\W;X)$ is the space of $X$-valued functions
whose norm has $p$-integrable/summable power. Similar definitions hold
with $p=\infty$ where the integrability requirement is replaced by
essential boundedness. When $p=2$ we denote this space by
{$\sobh{s}(\W;X)$}. The $\leb2(\W)$ and $\leb2(\boundary\W)$
inner products of two scalar, vector, or tensor-valued functions
$\varphi$ and $\psi$ is indicated with the brackets
\begin{equation}
  \label{eq:def:L2-inner-product}
  \ltwop{\varphi}{\psi}:=
  \int_\W\varphi(\vec x)\star\psi(\vec x)\d\vec x,
  \qquad
  \ltwop{\varphi}{\psi}_{\boundary\W}:=
  \int_{\boundary\W}\varphi(\vec x)\star\psi(\vec x)
  {\ds(\vec x)}
\end{equation}
where $\star$ stands for one of the arithmetic, Euclidean-scalar, or Frobenius
inner product in $\reals$, $\R d$, or $\realmats dd$ respectively and
$\ds$ is the $(d-1)$-dimensional measure element.

We refer to standard texts, e.g., \cite{Evans:10:book:Partial}, for
details about Sobolev spaces.

The boundary trace of a function $f\in\sob sp(\W;X)$ whenever it
exists, is denoted by $\restriction{f}{\partial\W}$ or just $f$ when
the trace is understood by the context. Since the domain $\W$ is
assumed of class $\holder01$, traces of functions in $\sobh1(\W;X)$
exist on $\partial\W$ and
the outward unit normal vector to $\W$ is denoted by
\index{$\normalto\W$} $\normalto\W(\vec{x})$ for $\measure S$-almost
every $\vec x$ on $\partial\W$.
If $\vec\psi\in\sobh1(\W;\R d)$, denoting by
  $\boundarytraceof\W{\vec\psi}$ the trace we respectively define
  $\vec\psi$'s \indexemph{normal trace}, and \indexemph{tangential
    trace} as
  \begin{equation}
  \normalto\W\normalto\W\inner\restriction{\vec\psi}{\boundary\W}
  ,\tand
  \explicittangentialto[\W]\psi
  .
\end{equation}
Our notation for some of the function spaces 
\index{$\linspace v$}\index{$\fatlinspace w$}\index{$\fatlinspace y$}
\begin{gather}
  \label{eq:def:space-of-zero-tangential-trace-H1-fields}
  \linspace v:= \setofsuch{
    \vec{\psi}\in\sobh1(\W;\R d)
  }{
    \explicittangentialto[\W]\psi=0
  }
  ,
  \\
  \label{eq:def:space-Y-equals-H1-H1d-L2dd}
  \fatlinspace{Y}:=
  \sobh1(\W) 
  \times
  \sobh1\qp{\W;\R d}
  \times
  \leb2\qp{\W;\Symmatrices d}, 
  \\
  \label{eq:def:space-V-zero-trace-H1-zero-tangential-trace-H1dd-L2dd}
  \fatlinspace{w}
  :=
  \sobhz1(\W) \times  \linspace v \times 
  \leb2(\W,\Symmatrices d)
  \subseteq
  \fatlinspace y
  ,
\end{gather}
endowed with the $\sobh1(\W;\R d)${-norm} for $\linspace v$ and
\begin{equation}
  \Norm{(\varphi , \vec\psi , \mat\Xi) }_{\fatlinspace{Y}}^2
  :=
  \Norm{ \varphi  }_{\sobh1(\W)}^2
  +
  \Norm{ \vec\psi  }_{\sobh1(\W)}^2
  +
  \Norm{ \mat\Xi }_{\leb2(\W)}^2
  \Foreach
  (\varphi ,\vec\psi ,\mat\Xi)
  \in
  \fatlinspace{Y}
  \supseteq
  \fatlinspace w
 .
\end{equation}
\subsection{Cordes conditions}
\label{sec:Cordes-conditions}
Let $d\in\naturals$ (typically $d=2,3$), \index{$d$}
\indexen{$\W$} be a bounded convex domain in $\R d$ of class $\holder01$
$\mat{A}\in\leb\infty(\W;\Symmatrices{d})$\index{$\mat A$} a
symmetric-matrix-valued function, $\vec b \in\leb\infty(\W;\R
d)$\index{$\vec b$} a vector field and $c \in
\leb\infty(\W)$\index{$c$} a scalar function which satisfy the
following \indexemph{Cordes condition}
\index{$\lambda$}\index{$\varepsilon$}
\begin{equation}
  \label{def:general-Cordes-condition}
  \frac{\norm{\mat A}^2
    + \fracl{\norm{\vec b}^2}{2\lambda}
    + (\fracl c\lambda)^2}{
    (\trace\mat A
    + \fracl c\lambda ) ^2
  }
  \leq
  \frac{1}{d+\varepsilon}
  \text{\Ae in } \W
\end{equation}
for some $\lambda>0$ and $\varepsilon\in (0,1)$.

In the special case $\vec b=0$ and $c=0$, we may take $\lambda=0$
and the Cordes condition (\ref{def:general-Cordes-condition}) is then replaced by
\begin{equation}
  \label{def:special-Cordes-condition}
  \frac{\norm{ \mat A }^2}{(\trace \mat A)^2}
  \leq
  \frac{1}{d-1+\varepsilon}
  \text{\Ae in }\W
\end{equation}
for some $\varepsilon\in (0,1)$. Since the right hand side of
(\ref{def:general-Cordes-condition}) and 
(\ref{def:special-Cordes-condition}) are decreasing with respect 
to $\varepsilon$, it suffices to find some $\bar{\varepsilon} >0$ 
which satisfies them and then considering $\varepsilon \in (0,\bar
{\varepsilon}]$ small enough. By the same argument, as the 
dimension increases, (\ref{def:general-Cordes-condition}) and 
(\ref{def:special-Cordes-condition}) become more stringent.
It is easy to show that in two dimensions, all symmetric positive 
definite matrices satisfy (\ref{def:special-Cordes-condition}), whereas 
this is not true in three (and higher) dimensions. For instance, taking
\begin{equation}
  \mat{A} :=
  \begin{bmatrix}
    1 & 0 & 0\\
    0 & 1 & 0\\
    0 & 0 & 5
  \end{bmatrix}
\end{equation}
in (\ref{def:special-Cordes-condition}) violates it. Nonetheless the
Cordes conditions (\ref{def:general-Cordes-condition}) or
(\ref{def:special-Cordes-condition}) cover a wide range of
applications including some nonlinear Hamilton--Jacobi--Bellman
equations
\citep[e.g.]{%
  Talenti:65:article:Sopra,%
  SmearsSuli:14:article:Discontinuous,%
  GallistlSuli:19:article:Mixed}.

If the boundary value of~(\ref{eq:nondivergence0-inhomogeneous}) be 
zero ($r = 0$) and the coefficients
satisfy 
$\vec b=0$, $c=0$ and (\ref{def:special-Cordes-condition}),
existence, uniqueness and stability of the strong solution in
$\sobh2(\W)\cap\sobhz1(\W)$ is proved by
\citet[Thm.~1]{Talenti:65:article:Sopra} for $\cont3$ smooth domains,
while a more general version for convex domains based on the
\indexemph{Miranda-Talenti regularity estimate}, is proved by
\citet[Thm.~3]{SmearsSuli:13:article:Discontinuous} while
\citet[Thm.~3]{SmearsSuli:14:article:Discontinuous} extend this result
to the case of a general nonlinear Hamilton--Jacobi--Bellman
equations, including that of (\ref{eq:nondivergence0-inhomogeneous})
with nonzero $c$ and $\vec b$ under
condition~(\ref{def:general-Cordes-condition}).
\subsection{Dirichlet boundary conditions}
\label{sec:Dirichlet-boundary-conditions}
We assume $r\in\sobh{3/2}(\boundary\W)$, i.e., $r$ is the
restriction (boundary trace) of a function, also denoted $r$, in
$\sobh2{(\W)}$ satisfying
\begin{equation}
  \label{eqn:trace-inequality}
  \inf
  \setofsuch{\Norm\phi_{\sobh2(\W)}}{\phi\in\sobh2(\W)\tand\phi-r\in\sobhz1(\W)}
  =:
  \Norm{r}_{\sobh{3/2}(\boundary\W)}
  \leq
  \constref{eqn:trace-inequality}
  \Norm{r}_{\sobh2(\W)}
  ,
\end{equation}
for some $\constref{eqn:trace-inequality}>0$ depending only on
$\W$. The function $v= u-r$ satisfies the problem
\begin{equation}
  \label{eq:homogeneous}
  \linop L v=f-\linop Lr
  \tand
  \boundarytraceof\W v=0
  .
\end{equation}
We will assume, except in the numerical experiments, that
$r=0$ in order to focus on the homogeneous boundary value problem
\begin{equation}
  \label{eq:nondivergence0}
  \linop L%
  u
  =
  f
  \tand
  \restriction u{\boundary\W}=0.
\end{equation} 
\subsection{A least-squares problem}
\label{sec:least-squares-problem}
We propose to formulate a least-squares alternative to
(\ref{eq:nondivergence0}) which allows for weaker solutions. Consider
 $0\leq\theta\leq1$\index{$\theta$} and start by introducing the linear operator
 \index{$\linoptheta$}
\begin{equation}
  \label{eq:def:linoptheta}
  \begin{gathered}
    \dfunkmapsto[.]{\linoptheta}{(\varphi, \vec\psi, \mat\Xi)}{\fatlinspace Y
    }{
      \mat A \frobinner \mat\Xi 
      + 
      \vec b
      \inner
      (
      \theta\vec \psi +(1-\theta) \grad \varphi
      )
      -
      c \varphi
      =:{\linoptheta}(\varphi,\vec\psi,\mat\Xi)
    }{\leb2(\W)}
  \end{gathered}
\end{equation}
The role of $\linoptheta$ is to approach the operator $\linop L$
from a mixed view point via
\begin{equation}
  \linop L\varphi = \linoptheta\triple\varphi{\grad\varphi}{\D^2\varphi}
  \text{ for }
  \varphi
  \text{ twice differentiable}.
\end{equation}
Although the aforementioned problem of finding a strong solution $u$
of (\ref{eq:nondivergence0}) in {$\sobh2(\W)\meet\sobhz1(\W)$} is
well-posed, working with such a high regularity assumption has
undesirable effects such as additional computational difficulties. As
we aim to a numerical scheme, to circumvent too stringent regularity
assumptions on $u$, we reformulate (\ref{eq:nondivergence0}) to an
appropriate alternative in $\sobh1(\W)$. The idea behind the
reformulation and the theory that follows is, similar to mixed
formulation, considering $u \in \sobh2(\W) \cap \sobhz1(\W)$ as $u \in
\sobhz1(\W)$ which also $\grad u \in\sobh1(\W; \R d) $. Motivated by
this reasoning, we introduce the following quadratic functional on
$\fatlinspace w$ \index{$E_\theta$}
\begin{equation}
  \begin{split}
    \label{functional:minimize}
    E_\theta(\varphi,\vec\psi,\mat\Xi)
    :=
    &
    \Norm{ \grad\varphi-\vec\psi}_{\leb2(\W)} ^2
    +
    \Norm{ \D\vec\psi-\mat\Xi}_{\leb2(\W)} ^2
    \\
    &
    +
    \Norm{ \rot\vec\psi }_{\leb2(\W)}^2
    +
    \Norm{ \linoptheta(\varphi, \vec\psi,\mat\Xi)  -f } _{\leb2(\W)}^2
  \end{split}
\end{equation} 
and consider the convex minimization problem of finding
\begin{equation}
  \label{eq:minimization}
  (u, \vec g, \mat  H)
  =\underset
  {\substack{
   \triple\varphi{\vec\psi}{\mat\Xi}\in\fatlinspace w
  }}
  \argmin%
  E_\theta(\varphi ,\vec\psi , \mat\Xi)
  .
\end{equation}
We recall that the
\indexemph{rotational} or \indexemph{curl operator}
\index{$\curl$}\index{$\rot$}
\begin{equation}
  \dfunkmapsto{\rot{}}{\vec\psi}{\sobh1(\W; \R d)}{\rot\vec\psi}{\leb2(\W)^{\hat d}} 
  \text{ for }
  \hat d
  :=
  \binom d2
  =
  \begin{cases}
    1
    &
    \text{ if }
    d=2
    ,
    \\
    3
    &
    \text{ if }
    d=3
  \end{cases}
\end{equation}
is such that in Cartesian coordinates one has
\begin{equation}
  \label{eq:def:curl/rot:coordinate-wise}
  \rot
  \begin{bmatrix}
  \psi_1\\
  \psi_2
  \end{bmatrix}
  =
  \partial_1{\psi_2} - \partial_2{\psi_1}, \qquad 
  \rot
  \begin{bmatrix}
  \psi_1\\
  \psi_2 \\
  \psi_3
  \end{bmatrix}
  =
  \begin{bmatrix}
  \partial_2{\psi_3}-\partial_3{\psi_2}\\
  \partial_3{\psi_1}-\partial_1{\psi_3} \\
  \partial_1{\psi_2}-\partial_2{\psi_1}
  \end{bmatrix}.
\end{equation}
{More generally, a coordinate and dimension $d$-independent
  definition of curl is the doubled
  skew-symmetric part of the Jacobian,
  \begin{equation}
    \label{eq:curl-and-skew-symmetric-Jacobian}
    \D\times\vec\psi:=\D\vec\psi-\D\vec\psi\transposed,
    \text{ for }\vec\psi\in\sobh1(\W;\R d)
  \end{equation}
  whereby when $d=2$ or $3$
  the usual curl is characterized by
  \begin{equation}
    \qp{\rot\vec\psi}\times\vec x=\qp{\D\times\vec\psi}\vec x
    \Foreach\vec x\in\R d,
  \end{equation}
  where, for (columns) $\vec x,\vec y\in\R{d}$, $\vec x\times\vec y$
  is the usual vector (external) product in $d=3$, and is
  $\det\disrowvectwo{\vec x}{\vec y}$ in $d=2$.  In terms of
  exterior algebra (and calculus) we are simply identifying
  elements of $\Lambda^2(\R d)$ (the alternating $2$-forms) (or
  skew-symmetric $d\times d$ matrices if preferred) with elements
  of $\R{\hat d}$, through the map $\mat J$
  \begin{equation}
    a\mapsto\dismattwo0{-a}a0=:\mat Ja,
    \tand
    \vec v=
    \discolvecithree v
    \mapsto
    \dismatthree0{-v_3}{v_2}{v_3}0{-v_1}{-v_2}{v_1}0
    =:
    \mat J\vec v,
  \end{equation}
  for $d=2$ and $d=3$ respectively.
  A useful consequence of this is that
  \begin{equation}
    \label{eq:inner-equals-frobinner-half}
    \vec v\inner\vec w
    =
    \frac12
    \mat J\vec v\frobinner\mat J\vec w.
  \end{equation}
}
\begin{Obs}[equivalence of (\ref{eq:nondivergence0}) and (\ref{eq:minimization})]
\label{rem:equivalence}
  If $u$ is a strong solution to (\ref{eq:nondivergence0}), then
  $(u,\grad u, \D^2u)$ minimizes the non-negative convex functional
  $E_\theta$. Since (\ref{eq:nondivergence0}) has a strong solution,
  the minimum value of $E_\theta$ is zero. Conversely, if $E_\theta$
  takes a minimum value at $(u, \vec g, \mat H)$, then $u$ is also a
  strong solution to (\ref{eq:nondivergence0}) and $\grad u= \vec g$,
  $\D^2u = \mat H$ in $\leb2(\W)$. Therefore the problem of finding
  strong solution to (\ref{eq:nondivergence0}) and problem
  (\ref{eq:minimization}) are equivalent. In the rest of the paper
  $\vec g$ and $\mat H$ will be synonymous with $\grad u$ and $\D^2 u$.
\end{Obs}
\subsection{Euler--Lagrange equations}
\label{sec:Euler-Lagrange-equations}
The Euler--Lagrange equation of the minimization problem 
(\ref{eq:minimization}) consist in finding\index{$\vec g$}\index{$\mat h$}
\({
  (u, \vec g, \mat H) \in \fatlinspace w%
}\)
such that
\begin{multline}
  \label{eq:Euler-Lagrange}
  \qa{ \grad u-\vec g, \grad \varphi - \vec \psi }
  +
  \qa{ \D\vec g - \mat H, \D\vec \psi - \mat \Xi } 
  \\
  +
  \qa{\rot\vec g,\rot\vec\psi}
  +
  \qa{ \linoptheta(u, \vec g,\mat H) , \linoptheta(\varphi, \vec \psi,\mat \Xi) }
  \\
  = 
  \qa{ f, \linoptheta(\varphi, \vec \psi,\mat \Xi) }
  ~\Foreach (\varphi, \vec\psi, \mat\Xi) \in \fatlinspace w.
\end{multline}
For numerical purposes it will be useful to rewrite
the Euler--Lagrange equation (\ref{eq:Euler-Lagrange})
in the following equivalent system-form
\begin{equation}
 \label{eq:E-L-system-3}
 \begin{aligned}
   \langle
   \grad u - \vec g
   &
   +
   (1-\theta)
   \linoptheta(u, \vec g,\mat H) 
   \vec b
   ,{
     \grad\varphi
   }
   \rangle
   -
   \ltwop{
     \linoptheta(u, \vec g,\mat H)c
   }{
     \varphi
   }
    \\
    &
    =
    (1-\theta)
    \ltwop{
      f\vec b}{
      \grad \varphi
    }
    -
    \ltwop{f c}{\varphi }
    \qquad\Foreach
    \phi\in\sobhz1(\W)
    ,
    \\
    \ltwop{
      \grad u - \vec g
    }{
      -\vec \psi }
    &
    +
    \ltwop{
      \D\vec g - \mat H
    }{
      \D\vec \psi }
    +
    \ltwop{
    \rot \vec g
    }{
    \rot \vec \psi}
        +
    \ltwop{
      \theta\linoptheta(u, \vec g,\mat H) \vec b
    }{
      \vec \psi
    }
    \\
    &
    =
    \theta
    \ltwop{
      f
      \vec b
    }{
      \vec \psi}
    \qquad\Foreach\vec\psi\in\linspace v
    ,
    \\
    \big\langle
      \linoptheta(u, \vec g,\mat H)
      \mat A
      &
      -
      \qp{\D\vec g-\mat H}
    ,{
      \mat\Xi
    }\big\rangle
    \\
    &
    =
    \ltwop{
      f \mat A
    }{
      \mat \Xi
    }
    \qquad\Foreach\mat\Xi\in\leb2(\W;\Symmatrices d).
  \end{aligned}
\end{equation}
\begin{Obs}[a Hessian-less approach]
  \label{obs:Hessian-less-approach}
  We may consider the \indexemph{Hessian-less objective functional}%
  \begin{equation}
    \label{functional:minimizing-2}
    (\varphi,\vec\psi)
    \mapsto
    E_\theta(\varphi,\vec\psi, \D\vec\psi),
  \end{equation} 
  the corresponding Euler-Lagrange equation be turned to finding
  $(u, \vec g) \in\sobhz1(\W)\times \linspace v $ such that
  \begin{multline}
    \label{eq:E-L-2}
    \qa{
      \grad u-\vec g
      ,
      \grad \varphi - \vec\psi }
    +
    \qa{
      \rot \vec g
      ,
      \rot \vec \psi }
    +
    \ltwop{
      \linoptheta(u , \vec g,\D\vec g)
    }{
      \linoptheta(\varphi, \vec \psi,\D\vec \psi)
    }
    \\
    = 
    \qa{ f,\linoptheta(\varphi, \vec \psi,\D\vec \psi) }
    \Foreach (\varphi, \vec\psi) \in \sobhz1(\W) \times \linspace v,
  \end{multline}
  or in equivalent system-form
  \begin{equation}
    \label{eq:E-L-system-2}
    \begin{aligned}
      \ltwop{
        \grad u-\vec g + (1-\theta) \linoptheta(u , \vec g,\D\vec g) \vec b
      }{
        \grad \varphi
      }
      -
      \ltwop{
        \linoptheta(u , \vec g,\D\vec g) c
      }{
        \varphi
      }
      &
      \\
      =
      (1-\theta)
      \ltwop{
        f\vec b}{
        \grad \varphi
      }
      -
      \ltwop{f c}{\varphi }
      \qquad\Foreach\varphi\in\sobhz1(\W)
      ,
      \\
      \ltwop{
        \grad u  - \vec g
      }{
        -\vec \psi
      }
      +
      \ltwop{
        \rot \vec g
      }{
        \rot \vec \psi
      }			
      +
      \ltwop{
        \theta\linoptheta(u , \vec g,\D\vec g)\vec b
      }{
        \vec \psi
      }
      +
      \ltwop{
        \linoptheta(u , \vec g,\D\vec g)\mat A
      }{
        \D\vec \psi
      }
      &
      \\
      =
      \theta
      \ltwop{
        f \vec b
      }{
        \vec \psi
      }
      +
      \ltwop{
        f \mat A
      }{
        \D\vec \psi
      }
      \qquad\Foreach\vec\psi\in \linspace v
      .
    \end{aligned}
  \end{equation}
\end{Obs}
\subsection{Bilinear forms}
\label{sec:bilinear-forms}
In keeping with (\ref{eq:Euler-Lagrange}) and (\ref{eq:E-L-2}),
we define the symmetric bilinear forms
\index{$a_\theta$}\index{$\hat a_\theta$}
\begin{equation}
  \funk
  {a_\theta}{
    \ppow{
      \fatlinspace{Y}
    }2
  }
  \reals
  \tand
  \funk{
   \hat{a}_\theta}{
    \qppow{
      \sobh1(\W)
      \times
      \sobh1\qp{\W;\R d}
    }2}
  \reals
 \end{equation}
by the expressions
\begin{gather}
  \label{eq:def:bilinear-a}
  \begin{split}
    \abil[a_\theta]{
      \varphi ,\vec\psi, \mat\Xi
    }{
      \varphi' , \vec\psi', \mat\Xi'
    }
    &
    := 
    \ltwop{
      \grad \varphi-\vec\psi
    }{
      \grad \varphi'-\vec\psi'
    }
    +
    \ltwop{
      \D\vec\psi - \mat\Xi
    }{
      \D \vec\psi' - \mat\Xi'
    }
    \\
    &
    \phantom{:=}
    +
    \ltwop{
      \rot\vec\psi
    }{
      \rot\vec\psi'
    }
    +
    \ltwop{
      \linoptheta(\varphi ,\vec\psi ,\mat\Xi )
    }{
      \linoptheta(\varphi',\vec\psi',\mat\Xi')
    }
  \end{split}
  \\
  \intertext{and}
  \label{eq:def:bilinear-a-hat}
  \abil[\hat{a}_\theta]{
    \varphi,
    \vec\psi
  }{
    \varphi',
    \vec\psi'
  }
  := 
  \abil[a_\theta]{
    \varphi ,
    \vec\psi,
    \D\vec\psi
  }{
    \varphi',
    \vec\psi',
    \D\vec\psi'
  }
\end{gather}
respectively for all $\triple{\varphi}{\vec\psi}{\vec\Xi}$ and
$\triple{\varphi'}{\vec\psi'}{\vec\Xi'}$ in the appropriate spaces.

Note that for any $v \in \sobh 2(\W) \cap \sobhz 1(\W)$ we have
$\grad v \in \linspace v$. In the analysis of the problem
(\ref{eq:minimization}) we need an estimate that is more general
than the classical \indexemph{Miranda--Talenti} estimate,
\begin{equation}
  \Norm{ \D^2 v}_{\leb2(\W)} \leq \Norm{ \lap v} _{\leb2(\W)}
  \Foreach
  v \in \sobh 2(\W) \cap \sobhz 1(\W).
\end{equation}
Indeed we need to bound $\Norm{ \div \vec \psi }_{\leb2(\W)}^2 + \Norm{\rot \vec \psi}_{\leb2(\W)}^2$ 
from below by $\Norm{ \D\vec \psi }_{\leb2(\W)}^2$.
\subsection{The role of the curl and Maxwell's estimate} 
\label{lem:Maxwell-estimate}
A motivation for considering the $\Norm{ \rot \vec \psi }_{\leb2(\W)}^2$ 
in the functional $E_\theta$ lies in the fact, known as
\indexemph{Maxwell estimate},
that since $\W$ is a convex domain, for any $\vec\psi \in\linspace v$,
we have
\begin{equation}
  \label{eqn:Maxwell-estimate}
  \Norm{ \D\vec\psi }_{\leb2(\W)}^2
  \leq
  \Norm{ \div \vec\psi }_{\leb2(\W)}^2
  +
  \Norm{\rot\vec\psi}_{\leb2(\W)}^2
  .
\end{equation}
We refer to \citet{CostabelDauge:99:article:Maxwell} for more details.
\section{Coercivity and continuity of the cost functional}
\label{sec:coercivity-continuity-of-cost-functional}
We now show that problem (\ref{eq:Euler-Lagrange}) is well-posed via a
Lax--Milgram approach. To effect this it is sufficient to show that
the bilinear form $a_\theta$, defined in \S~\ref{sec:bilinear-forms}),
is coercive and continuous. After discussing our main strategy in
\S~\ref{sec:our-approach}, and giving some preliminaries, including a
Miranda--Talenti type consequence of the Cordes condition in
Lemma~\ref{lem:Miranda-Talenti-special-Cordes}. This is further developed
into Theorem~\ref{the:miranda-talenti}, which for $\theta = 0$ is proved by
\citet[Lem.~2.1]{GallistlSuli:19:article:Mixed} and we extend it for any
$0 \leq \theta \leq 1$.

Based on these results we then prove the main results of this
section, namely, that $\hat a_\theta$ and $a_\theta$ are coercive in
theorems \ref{the:coercivity-u,g} and \ref{the:coercivity-u,g,H},
respectively and continuity is shown in \S~\ref{sec:continuity-a_theta}.
Finally, in \S~\ref{relaxing-problem} and
\S~\ref{sec:nonzero-boundary-values} we show the necessity of the zero
tangential-trace condition and adapt the minimization problem to the
case of nonzero boundary values problem.
\subsection{Key ideas of our least-squares approach}
\label{sec:our-approach}
We develop the proof of $a_\theta$'s coercivity in two steps. First,
we prove that $\hat{a}_\theta$ is coercive on $\sobhz1(\W)
\times \linspace v$; the key of the proof is considering an
appropriate operator on $\sobhz1(\W) \times \sobh1\qp{\W;\R d}$ say
$\linop{D}$ which for any $(\varphi , \vec\psi) \in \sobhz1(\W) \times
\linspace v$ is close to $\linoptheta(\varphi , \vec\psi, \D\vec\psi)$ 
and for some constant $C>0$
\begin{equation}
  \Norm{\grad \varphi - \vec\psi}_{\leb2(\W)}^2
  +
  \Norm{\rot \vec \psi}_{\leb2(\W)}^2
  +
  \Norm { 
    \linop{D}(\varphi , \vec\psi)
  }_{\leb2(\W)}^2 
  \geq
  C\qp{
  \Norm{\varphi }_{\sobh1(\W)}^2
  +
  \Norm{\vec\psi}_{\sobh1(\W)}^2
  }.
\end{equation}
Then, by comparing $\D\vec\psi$ with $\mat\Xi$ we get the coercivity
of $a_\theta$ on $\fatlinspace w$. 

Recalling the notation from (\ref{def:general-Cordes-condition}) and
(\ref{def:special-Cordes-condition}) introduce the \indexemph{scaling
function}\index{$\gamma$}
\begin{equation} 
  \label{def:gamma}
  \gamma
  :=
  \begin{cases}
    \frac{\trace\mat A}{\norm{ \mat A }^2}
    &\text{ if }\lambda = 0,
    \\
    \frac{\trace\mat A
      + c/\lambda}{\norm{ \mat A} ^2
      + \norm{ \vec b } ^2/2\lambda
      + (c/\lambda)^2}
    &\text{ if }\lambda> 0,
  \end{cases}
\end{equation}
which was used in \citet{SmearsSuli:13:article:Discontinuous}.
Uniform ellipticity (\ref{def:uniformly-elliptic}), non-negativity of $c$ and 
uniform boundedness of the coefficients of $\linop L$ imply that $\inf_\W\gamma>0$ 
and %
\begin{equation}
  \infty>
  \Norm{\gamma}_{\leb\infty(\W)}
  =:
  \constdef[\linop L]{const:sup-gamma}.
\end{equation}
\begin{Lem}[a Miranda--Talenti estimate]
  \label{lem:Miranda-Talenti-special-Cordes}
  If $\mat A$ satisfies the Cordes condition with $\lambda=0$
  (\ref{def:special-Cordes-condition}), then for any
  $\vec\psi\in\linspace v$
  \begin{equation}
    \label{eqn:Miranda-Talenti-special-Cordes}
    \Norm{ \rot \vec \psi }_{\leb2(\W)}^2
    +
    \Norm{ \mat A\frobinner \D\vec \psi }_{\leb2(\W)}^2
    \geq
    \constref{const:special-Cordes}
    \Norm{ \D\vec \psi }_{\leb2(\W)}^2
  \end{equation}
  where
  \begin{equation}
    \constdef[\linop L]{const:special-Cordes}
    :=
    \frac{
      (1-\sqrt{1-\varepsilon})^2
    }{\qgroup{
        \maxi{\constref{const:sup-gamma}^2}{1}
    }}.
  \end{equation}
\end{Lem}
\begin{Proof}
  The definition of $\gamma$ in (\ref{def:gamma})
  and the Cordes condition (\ref{def:special-Cordes-condition}) imply
  that
  \begin{equation}
    \norm{ \gamma \mat A - \eye }^2
    = d - \frac{\norm{\mat A}^2}{(\trace \mat A)^2}
    \leq 1-\varepsilon.
  \end{equation}
  Hence we have
  \begin{equation}
    \label{eqn:triangle-ineq-cordes}
    \Norm{ (\gamma \mat A- \eye)\frobinner \D \vec \psi }_{\leb2(\W)}
    \leq
    \sqrt{1-\varepsilon} \Norm{ \D \vec \psi }_{\leb2(\W)}.
  \end{equation}
  Adding and subtracting $\eye \frobinner \D \vec \psi $ and then using 
  (\ref{eqn:Maxwell-estimate}) and (\ref{eqn:triangle-ineq-cordes}) lead to
  \begin{multline}
    \Norm{\rot \vec \psi}_{\leb2(\W)}^2
    +
    \Norm{\gamma \mat A \frobinner \D \vec \psi }_{\leb2(\W)}^2
    =
    \Norm{\rot \vec \psi}_{\leb2(\W)}^2
    +
    \Norm{(\gamma \mat A- \eye + \eye )\frobinner \D \vec \psi }_{\leb2(\W)}^2
    \\
    \geq
    \qp{
      \sqrt{\Norm{\rot \vec \psi}_{\leb2(\W)}^2 +\Norm{\div \vec \psi}_{\leb2(\W)}^2 }
      - \Norm{(\gamma \mat A - \eye)\frobinner \D \vec \psi  }_{\leb2(\W)}
    }^2
    \\
    \geq 
    (1-\sqrt{1-\varepsilon})^2 \Norm{\D \vec \psi}_{\leb2(\W)}^2,
  \end{multline}
  from which we conclude.
\end{Proof}
\begin{Def}[an auxiliary perturbed mixed Laplace operator]
  \label{def:auxilliary-Laplace-operator}
  Recalling the parameter $\lambda$ entering the Cordes
  condition~(\ref{def:general-Cordes-condition})
  we define the \indexemph{perturbed mixed Laplace operator}
  \index{$\linop D_\lambda$}
  $\funk{\linop{D}_\lambda}{\sobhz1(\W)\times\sobh1(\W;\R d)}{\leb2(\W)}$
  as
  \begin{equation}
    \label{def:functional-operator-L_lambda}
    \linop{D}_\lambda(\varphi, \vec \psi )
    := \div \vec \psi - \lambda \varphi
    .
  \end{equation}
  The name of this operator, which we need for our proof, rests on the
  fact that our intention behind the variable $(\varphi,\vec\psi)$ is for
  it to equate $(u,\grad u)$ and obtain the characteristic operator
  \begin{equation}
    \label{eqn:auxilliary-characteristic-operator-ellipitic}
    \linop{D}_\lambda(u,\grad u)=\lap u-\lambda u.
  \end{equation}
  A similar idea
  of using this operator can be found in
  \citet[eq. (2.12)]{SmearsSuli:14:article:Discontinuous}.
\end{Def}
\begin{Def}[an auxiliary parameter-dependent norm]
 Given two parameters $0\leq\theta\leq1$ and $\lambda>0$, as introduced
 before, define the following norm for 
 $(\varphi,\vec\psi)\in\sobhz1(\W)\times\linspace v$
 \index{$\Norm{ \qp{\varphi , \vec\psi} }_{\lambda,\theta}^2$}
 \begin{equation}
   \label{eq:def:parameter-dependent-seminorm}
   \Norm{ \qp{\varphi , \vec\psi} }_{\lambda,\theta}^2 
   :=
   \Norm{ \D \vec\psi }_{\leb2(\W)}^2
   +
   2\lambda \Norm{ \theta \vec\psi + (1-\theta)\grad \varphi }_{\leb2(\W)}^2
   +
   \lambda^2 \Norm{ \varphi }_{\leb2(\W)}^2 
   .
 \end{equation}
\end{Def}
\begin{Obs}[Poincaré's inequality]
  Let $\W$ be a bounded domain, then for any
  $(\varphi,\vec\psi)\in\sobhz1(\W)\times\linspace{v}$ there
  corresponds $\constref[\W]{eqn:Poincare}>0$ such that
  \begin{equation}
    \label{eqn:Poincare}
    \Norm{\D\vec\psi}_{\leb2(\W)}^2
    \geq
    \constref[\W]{eqn:Poincare}
    \Norm{\vec\psi }_{\sobh1(\W)}^2, 
    \tand
    \Norm{\grad \varphi}_{\leb2(\W)}^2
    \geq
    \constref[\W]{eqn:Poincare}
    \Norm {\varphi}_{\sobh1(\W)}^2.
\end{equation}
\end{Obs}
\begin{The}[a modified Miranda--Talenti estimate]
 \label{the:miranda-talenti}
 If $\W$ is a bounded open convex subset of $\mathbb{R}^d$, $0<\rho<2$ and 
 $0 \leq \theta \leq 1$ then for any $(\varphi,\vec\psi)\in\sobhz1(\W) \times\linspace v$ 
 we have
 \begin{multline}
   \label{eqn:Miranda-Talenti} 
   \qp{1- \fracl{\rho}{2}}
   \Norm{\qp{\varphi , \vec\psi} }_{\lambda,\theta}^2
   \\
   \leq
   \Norm{\rot \vec \psi }_{\leb2(\W)}^2
   +
   \Norm{ \linop{D}_\lambda(\varphi, \vec \psi ) }_{\leb2(\W)}^2
   +
   \qp{\theta^2 + (1-\theta)^2} 
   \fracl{\lambda}{\rho}
   \Norm {\grad \varphi - \vec\psi }_{\leb2(\W)}^2 
   .
 \end{multline}
\end{The}
\begin{Proof}
 We start the proof by noting that thanks to 
    $\varphi\in\sobhz1(\W)$ and $\vec\psi\in\sobh1(\W)^d$ we have
    \begin{equation}
      \ltwop{\div\vec\psi}\varphi
      =
      -\ltwop{\vec\psi}{\grad\varphi}.
    \end{equation}
 Using the Maxwell estimate \eqref{eqn:Maxwell-estimate} and expanding 
 $\Norm{\theta \vec\psi + (1-\theta)\grad \varphi }_{\leb2(\W)}^2$
  imply
  \begin{multline}
    \Norm{ \qp{\varphi , \vec\psi}}_{\lambda,\theta}^2
    \leq
    \Norm{\rot \vec \psi}_{\leb2(\W)}^2
    +
    \Norm{ \linop{D}_\lambda(\varphi, \vec \psi ) }_{\leb2(\W)}^2
    -
    2\lambda \langle \vec \psi, \grad \varphi\rangle
    \\
    +
    2\lambda \theta^2 \Norm{\vec \psi}_{\leb2(\W)}^2
    +
    2\lambda (1-\theta)^2 \Norm{\grad \varphi}_{\leb2(\W)}^2
    +
    4\lambda \theta (1-\theta)\langle \vec\psi, \grad\varphi \rangle
    \\
    =
    \Norm{\rot \vec \psi}_{\leb2(\W)}^2
    +
    \Norm{ \linop{D}_\lambda(\varphi, \vec \psi ) }_{\leb2(\W)}^2
    +
    2\lambda \theta^2 \langle \vec\psi, \vec\psi - \grad \varphi \rangle
    +
    2\lambda (1-\theta)^2 \langle \grad\varphi,   \grad \varphi - \vec\psi \rangle  
    .
   \end{multline} 
  Applying a weighted Young's inequality leads to
  \begin{multline}
    \Norm{ \qp{\varphi , \vec\psi} }_{\lambda,\theta}^2
    \leq
    \Norm{\rot \vec \psi}_{\leb2(\W)}^2
    +
    \Norm{ \linop{D}_\lambda(\varphi, \vec \psi ) }_{\leb2(\W)}^2
    +
    \lambda \theta^2 \rho \Norm{\vec\psi}_{\leb2(\W)}^2
    \\
    +
    \frac{\lambda \theta^2}{\rho} \Norm{\grad\varphi - \vec\psi}_{\leb2(\W)}^2
    +
    \lambda (1-\theta)^2 \rho \Norm{\grad\varphi}_{\leb2(\W)}^2
    +
    \frac{\lambda (1-\theta)^2}{\rho} \Norm{\grad\varphi - \vec\psi}_{\leb2(\W)}^2
    .
  \end{multline} 
  By subtracting $\lambda \theta^2 \rho \Norm{\vec\psi}_{\leb2(\W)}^2 +  
  \lambda (1-\theta)^2 \rho \Norm{\grad\varphi}_{\leb2(\W)}^2$ from 
  both sides and reversing the inequality we get
  \begin{equation}
    \begin{split} 
      \Norm{ \rot \vec \psi }_{\leb2(\W)}^2
      &
      +
      \Norm{ \linop{D}_\lambda(\varphi, \vec \psi ) }_{\leb2(\W)}^2
      +
      \frac{\lambda \theta^2}{\rho} \Norm{\grad\varphi - \vec\psi}_{\leb2(\W)}^2
      \\
      &
      +
      \frac{\lambda (1-\theta)^2}{\rho}
      \Norm{\grad\varphi - \vec\psi}_{\leb2(\W)}^2
      \\
      \geq
      &
      \Norm{ \qp{\varphi , \vec\psi} }_{\lambda,\theta}^2 
      -
      \lambda \theta^2 \rho \Norm{\vec\psi}_{\leb2(\W)}^2
      -
      \lambda (1-\theta)^2 \rho \Norm{\grad\varphi}_{\leb2(\W)}^2
      \\
      =
      &
      \Norm{\D \vec\psi}_{\leb2(\W)}^2 
      + 
      2\lambda(1-\fracl{\rho}{2}) 
      \Norm{\theta \vec\psi + (1-\theta)\grad\varphi}_{\leb2(\W)}^2
      \\
      &
      +
      \lambda^2 \Norm{\varphi}_{\leb2(\W)}^2
      -
      2\lambda\theta(1-\theta)\rho\ltwop{\vec\psi}{\grad\varphi}
      \\
      =
      &
      \Norm{\D \vec\psi}_{\leb2(\W)}^2 
      + 
      2\lambda(1-\fracl{\rho}{2}) 
      \Norm{\theta \vec\psi + (1-\theta)\grad\varphi}_{\leb2(\W)}^2
      \\
      &
      +
      \lambda^2 \Norm{\varphi}_{\leb2(\W)}^2
      +
      2\theta(1-\theta)\rho\ltwop{\div\vec\psi}{\lambda\varphi}
      \\
      \geq
      &
      2\lambda(1-\fracl{\rho}{2}) 
      \Norm{\theta \vec\psi + (1-\theta)\grad\varphi}_{\leb2(\W)}^2
      +
      \Norm{\D \vec\psi}_{\leb2(\W)}^2 
      \\
      &
      +
      \lambda^2 \Norm{\varphi}_{\leb2(\W)}^2
      -
      \theta(1-\theta)\rho\Norm{\div\vec\psi}_{\leb2(\W)}^2
      -
      \theta(1-\theta)\rho\lambda^2\Norm{\varphi}_{\leb2(\W)}^2
      \\
      \geq
      &
      2\lambda(1-\fracl{\rho}{2}) 
      \Norm{\theta \vec\psi + (1-\theta)\grad\varphi}_{\leb2(\W)}^2
      +
      (1-\fracl{\rho}{4})\Norm{\D\vec\psi}_{\leb2(\W)}^2
      \\
      &
      +
      (1-\fracl{\rho}{4})\lambda^2\Norm{\varphi}_{\leb2(\W)}^2
      \\
      \geq
      &
      (1-\fracl{\rho}{2}) \Norm{ (\varphi , \vec\psi) }_{\lambda,\theta}^2
      ,
    \end{split}
  \end{equation}
  as claimed. 
\end{Proof}  
\begin{The}[coercivity of $\hat{a}_\theta$]
  \label{the:coercivity-u,g}
  Let $\W$ be a bounded convex open subset of $\mathbb{R}^d$
  and the coefficients $\mat A, \vec b, c $ satisfy the Cordes condition (either
  (\ref{def:general-Cordes-condition}) with $\lambda >0$ or
  (\ref{def:special-Cordes-condition}) with $\vec b = 0, c = 0$
  and $\lambda=0$). Then the restricted bilinear form $\hat a_\theta$ defined in
  (\ref{eq:def:bilinear-a-hat}) satisfies
  \begin{equation}
    \label{eqn:coercivity-theta-ug}
    \begin{split}
      &
      \abil[\hat a_\theta]{\varphi,\vec\psi}{\varphi,\vec\psi}
      \geq 
      \constref{const:coercivity-theta-ug}
      \qp{
        \Norm{\varphi }_{\sobh1(\W)}^2
        +
        \Norm{ \vec\psi }_{\sobh1(\W)}^2
      }
    \end{split}
  \end{equation} 
  for all $(\varphi,\vec\psi)\in\sobhz1(\W) \times\linspace v$
  where
  \begin{equation}
    \constdef[\W,\theta,\lambda,\varepsilon,{\linop L}]{const:coercivity-theta-ug}
    :=
    \left\{
    \begin{matrix}
      \frac{\constref{eqn:Poincare}}{2}
      \qp{\mini{1}{
          \frac{\qp{1-\sqrt{1-\varepsilon}}^2
            \constref{eqn:Poincare}
          }{
            \qp{\maxi{\constref{const:sup-gamma}^2}{1}}}}
      }
      &
      \text{ if }\lambda=0
      ,
      \\[12pt]
      \frac{\qp{\sqrt[4]{1-\varepsilon}-\sqrt{1-\varepsilon}}^2
        \qp{ \mini{{\constref{eqn:Poincare}}}{{4}\lambda^2}}
      }{
        \maxi{
       	  {2\qp{\sqrt[4]{1-\varepsilon}-\sqrt{1-\varepsilon}}^2
      	    \constref{eqn:Poincare}
       	  }
       	  + 
       	  \frac{2\lambda\qp{\theta^2+\qp{1-\theta}^2}
       	  }{1 - \sqrt{1-\varepsilon}}
       	}{
     	  \maxi{1}{\constref{const:sup-gamma}^2}
     	} 
      } 
      &
      \text{ if }\lambda>0
      .
    \end{matrix}
    \right.
  \end{equation}
\end{The}
\begin{Proof}
  We distinguish two cases according to whether $\lambda=0$
  or $\lambda>0$.
  \begin{CasesList}
  \item
    Consider $\lambda=0$, then Lemma
    \ref{lem:Miranda-Talenti-special-Cordes} leads to
    \begin{multline}
      \label{eqn:Gallistl-to-L}
      \Norm{ \grad \varphi - \vec\psi } _{\leb2(\W)}^2
      +
      \Norm{\rot \vec \psi}_{\leb2(\W)}^2
      +
      \Norm{\mat A\frobinner \D\vec \psi } _{\leb2(\W)}^2 
      \\
      \geq
      \Norm{ \grad \varphi - \vec\psi } _{\leb2(\W)}^2
      +
      \frac{\qp{1-\sqrt{1-\varepsilon}}^2}{
        \qp{ \maxi{\constref{const:sup-gamma}^2} {1}
      }}
      \Norm{\D\vec\psi }_{\leb2(\W)}^2
      .
    \end{multline}
    Putting
    \begin{equation}
      \constdef{const:hatC}
      :=\qpmini
          {1}
          {\frac{ 
              \qp{1-\sqrt{1-\varepsilon}}^2}{ \qpmaxi{
                \constref{const:sup-gamma}^2
              }{1}
            }
            \constref{eqn:Poincare}
          },
    \end{equation}
    and using Young's and Poincaré's inequality we arrive at
    \begin{equation}
      \label{eqn:op-L-2}
      \begin{aligned}
        \Norm{ \grad \varphi -\vec\psi } _{\leb2(\W)}^2
        &
        +
        \Norm{\rot \vec \psi}_{\leb2(\W)}^2
        +
        \Norm{ \mat A\frobinner \D\vec \psi} _{\leb2(\W)}^2
        \\
        &
        \geq
        \frac{
          \constref{const:hatC}
        }{2}
        \Norm{\grad\varphi}_{\leb2(\W)}^2 
        -
        \constref{const:hatC} \Norm{  \vec\psi  }_{\leb2(\W)}^2 
        \\
        &
        \phantom\geq
        +
        \frac{\qp{1-\sqrt{1-\varepsilon}}^2}{\qp{\maxi{\constref{const:sup-gamma}^2}{1}}}
        \constref{eqn:Poincare}\qp{
          \Norm{ \vec\psi  }_{\leb2(\W)}^2
          +
          \Norm{ \D\vec\psi  }_{\leb2(\W)}^2 }
        \\
        &
        \geq
        \frac{
          \constref{const:hatC}
          \constref{eqn:Poincare}
        }2
        \qp{
          \Norm{ \varphi }_{\sobh1(\W)}^2 
          +
          \Norm{ \vec\psi  }_{\sobh1(\W)}^2
        },
      \end{aligned}
    \end{equation}
    which establishes the result for zero $\lambda$.
  \item
    Suppose {$\lambda>0$},
    let $\rho=2-2\sqrt{1-\varepsilon}$ and define
    \begin{equation}
      \constdef[\lambda,\theta,\varepsilon]{const:mu}:=
      \frac{\lambda(\theta^2+(1-\theta)^2)}{2-2\sqrt{1-\varepsilon}}, 
    \end{equation}
    then from the Miranda--Talenti estimate, Theorem~\ref{the:miranda-talenti},
    we first note that
    \begin{equation}  
      \label{eq:proof:a-hat-coercive:Miranda-Talenti-for-op-L_lambda}
      \constref{const:mu} \Norm{ \grad\varphi - \vec\psi } _{\leb2(\W)}^2
      +
      \Norm{ \linop{D}_\lambda (\varphi, \vec\psi  ) } _{ \leb2(\W)} ^2 
      +
      \Norm{\rot \vec \psi}_{\leb2(\W)}^2
      \geq
      \sqrt{1-\varepsilon} \Norm{ (\varphi , \vec\psi) }_{\lambda,\theta}^2
      .
    \end{equation}
    On the other hand, the \CBS implies
    \begin{multline}
      \label{eq:op-(gamma.L-L_lambda)-0}
      \Norm{ 
        \gamma \linoptheta(\varphi,  \vec\psi , \D\vec\psi )
        - 
        \linop{D}_\lambda(\varphi, \vec\psi )  
      }_ {\leb2(\W)}^2
      \\
      =
      \Norm{
        (\gamma \mat A-\eye):\D\vec\psi 
        +
        \gamma
        \vec b\inner
        \qp{\theta\vec\psi  + (1-\theta)\grad \varphi}
        +
        \qp{\lambda -\gamma c}\varphi}_{\leb2(\W)}^2        
      \\
      \leq
      \Normbig{
        \norm{ \gamma \mat A-\eye }^2
        +
        \fraclff{\normreg\gamma^2 \norm{\vec b}^2}{2 \lambda}
        +
        \fraclff{\normreg{\lambda-\gamma c}^2}{\lambda^2}
      }_{\leb\infty(\W)}
      \\
      \phantom\leq
      \qp{
        \Norm{ \D\vec\psi }_{\leb2(\W)}^2
        +
        2\lambda
        \Norm{
          \theta\vec\psi   
          +
          (1-\theta)\grad \varphi
        }_{\leb2(\W)}^2
        +
        \lambda^2 \Norm{\varphi}_{\leb2(\W)}^2
      }
      .
    \end{multline}
    Rearranging the first factor in the right-hand side of
    \eqref{eq:op-(gamma.L-L_lambda)-0} and recalling the definition of
    the scaling function $\gamma$ (\ref{def:gamma}), as well as the
    the Cordes condition (\ref{def:general-Cordes-condition}) yield
    \begin{multline}
      \norm{\gamma \mat A-\eye }^2
      +
      \frac{\norm{ \gamma } ^2 \norm{ \vec{b} }^2}{2 \lambda}
      +
      \frac{\norm{ \lambda -\gamma c } ^2}{\lambda^2}
      \\
      =
      d+1
      -
      2\gamma\qp{\trace\mat A+\frac{c}{\lambda}}
      +
      \norm{ \gamma }^2(\norm{ \mat A }^2
      +
      \frac{\norm{ \vec b }^2}{2\lambda} 
      +
      \frac{\norm{ c }^2}{\lambda^2})
      \leq
      1-\varepsilon
      .
    \end{multline}
    {Owing to definition \eqref{eq:def:parameter-dependent-seminorm} we have}
    \begin{equation}
      \label{eqn:op-(gamma.L-L_lambda)}
      \Norm{  
        \gamma 
        \linoptheta(\varphi, \vec\psi ,\D\vec\psi )
        -
        \linop{D}_\lambda(\varphi, \vec\psi ) 
      }_ {\leb2(\W)}^2 
      \leq 
      (1-\varepsilon) \Norm{(\varphi , \vec\psi)}_{\lambda,\theta}^2 
      .
    \end{equation}
     Adding--subtracting
    $\linop{D}_\lambda(\varphi, \vec\psi )$, some manipulations,
    \eqref{eq:proof:a-hat-coercive:Miranda-Talenti-for-op-L_lambda}
    and
    (\ref{eqn:op-(gamma.L-L_lambda)})
    lead us to
    \begin{equation}
      \label{eqn:proof:coercivity-of-a-hat:coercivity-in-theta-lambda-norm}
      \begin{aligned}
        &
        \constref{const:mu}
        \Norm{\grad\varphi-\vec\psi}_{\leb2(\W)}^2
        +
        \Norm{\rot\psi}_{\leb2(\W)}^2
        +
        \Norm{
          \gamma\linoptheta(\varphi, \vec\psi , \D\vec\psi )
        }_{\leb2(\W)}^2 
        \\
        &=
        \constref{const:mu} \Norm{ \grad \varphi - \vec\psi  } _{\leb2(\W)}^2
        +
        \Norm{\rot \psi}_{\leb2(\W)}^2
        \\
        &\phantom=
        +
        \Norm{ 
          \qgroup{
            \gamma
            \linoptheta(\varphi, \vec\psi , \D\vec\psi )
            -
            \linop{D}_\lambda(\varphi, \vec\psi )
            +
            \linop{D}_\lambda(\varphi, \vec\psi )} 
        }_{\leb2(\W)}^2
        \\ 
        &\geq
        \bigg(
          \qpsqrt{
            \constref{const:mu}
            \Norm{ \grad \varphi - \vec\psi  } _{\leb2(\W)}^2
            +
            \Norm{\rot \psi}_{\leb2(\W)}^2
            +
            \Norm{ \linop{D}_\lambda(\varphi,\vec\psi ) } _{\leb2(\W)}^2
          }
          \\
          &\phantom{\geq\Big(}
          -
          \Norm{
            \gamma(\linoptheta(\varphi, \vec\psi , \D\vec\psi )
            -
            \linop{D}_\lambda(\varphi, \vec\psi ) 
          }_{\leb2(\W)}
        \bigg)^2
        \\
        &\geq
        \qp{ \sqrt[4]{1-\varepsilon}-\sqrt{1-\varepsilon}}^2 
        \Norm{(\varphi , \vec\psi)}_{\lambda,\theta}^2
        \mathchanges, 
      \end{aligned}
    \end{equation}
    where in the last step we use the $\lambda,\theta$-norm defined in
    (\ref{eq:def:parameter-dependent-seminorm}).
    
    Young's inequality and Poincaré's inequality (\ref{eqn:Poincare}) combined with
    (\ref{eqn:proof:coercivity-of-a-hat:coercivity-in-theta-lambda-norm})
    imply that
    \begin{multline}%
      \label{eqn:proof:coercivity-of-a-hat:from-theta-lambda-to-energy-norm}
      \qpmaxi{\frac12{
          \qp{\sqrt[4]{1-\varepsilon}-\sqrt{1-\varepsilon}}^2
      	  \constref{eqn:Poincare}
       	}
       	+
        \constref{const:mu}		
      }{\maxi{\constref{const:sup-gamma}^2}{1}}
      \aqua[\hat a_\theta]{\varphi,\vec\psi}
      \\
      \begin{split}
        &
        \geq
        \frac{\qp{\sqrt[4]{1-\varepsilon}-\sqrt{1-\varepsilon}}^2\constref{eqn:Poincare}}{2} 
        \Norm{ \grad \varphi - \vec\psi   }_{\leb2(\W)}^2
        \\
        &\phantom\geq
        +
        \constref{const:mu} \Norm{ \grad \varphi - \vec\psi   }_{\leb2(\W)}^2
        +
        \Norm{\rot \psi}_{\leb2(\W)}^2
        +
        \Norm{
          \gamma\linoptheta(\varphi, \vec\psi , \D\vec\psi )
        }_{\leb2(\W)}^2 
        \\
        &
        \geq
        \frac{
          \qp{\sqrt[4]{1-\varepsilon}-\sqrt{1-\varepsilon}}^2\constref{eqn:Poincare}}{4}
        \Big(
        \Norm{\grad \varphi}_{\leb2(\W)}^2
        -
        2\Norm{\vec \psi}_{\leb2(\W)}^2
        \\
        &\phantom\geq
        +
        \frac{4}{\constref{eqn:Poincare}}
        \qp{
          \Norm{\D\vec\psi}_{\leb2(\W)}^2
          +
          \lambda^2\Norm{\varphi}_{\leb2(\W)}^2
        }
        \Big)
        \\
        &
        \geq
        \frac{\qp{\sqrt[4]{1-\varepsilon}-\sqrt{1-\varepsilon}}^2\constref{eqn:Poincare}}{4}
        \\
        &\phantom\geq
        \times
        \qp{
          2\Norm{\vec\psi}_{\sobh1(\W)} ^2
          +
          \Norm{\grad \varphi}_{\leb2(\W)}^2
          +
          \frac{4\lambda^2}{\constref{eqn:Poincare}}
          \Norm{\varphi}_{\leb2(\W)}^2 
        }
        .
      \end{split}
    \end{multline}
    We then deduce the coercivity
    \begin{equation}
      \label{eqn:proof:coercivity-of-a-hat:conclusion}
      {\aqua[\hat a_\theta]{\varphi,\vec\psi}}
      \geq 
      \constref{const:coercivity-theta-ug}
      (\Norm{ \varphi }_{\sobh1(\W)}^2 
      +
      \Norm{ \vec\psi  }_{\sobh1(\W)}^2),
    \end{equation}
  \end{CasesList}
  which is the claim for $\lambda$ strictly positive.
\end{Proof}
\begin{The}[coercivity of $a_\theta$]%
  \label{the:coercivity-u,g,H}
  Under the same assumptions of Theorem~\ref{the:coercivity-u,g}
  we have
  \begin{equation}
    \label{eqn:coercivity-u,g,H}
    \begin{split}
      \abil[a_\theta]{\varphi,\vec\psi, \mat \Xi}{\varphi,\vec\psi, \mat \Xi}
      \geq
      \constref%
          {const:coercivity-theta-ugH}\qp{
            \Norm{ \varphi }_{\sobh1(\W)}^2
            +
            \Norm{ \vec  \psi }_{\sobh1(\W)}^2
            +
            \Norm{ \mat \Xi }_{\leb2(\W)}^2
          }
    \end{split}
  \end{equation}
  for all $(\varphi,\vec\psi, \mat \Xi)\in \fatlinspace w$ where 
  {%
    \begin{equation}
      \constdef[\W,\lambda,\theta,\varepsilon,{\linop L}]{const:coercivity-theta-ugH}
      :=
      \frac{
        \mini{
          {\constref{const:coercivity-theta-ug}}
        }{
          4\Norm{\mat A}_{\leb{\infty}(\W)}^2
        }
      }{
        \maxi8{16 \Norm{\mat A}_{\leb{\infty}(\W)}^2}
      }
    \end{equation}
  }
\end{The}
\begin{Proof}
    Posing
    \begin{equation}
      \label{eq:replace}
      \mat M := \D\vec\psi - \mat\Xi,
    \end{equation}
    maximum property, {some algebraic manipulations and 
    Young's inequality together with 
    Theorem \ref{the:coercivity-u,g} respectively imply the first, 
    second and third inequalities of the following:
    \begin{multline}
      \qpmaxi{1}{2 \Norm{\mat A}_{\leb{\infty}(\W)}^2}
      \abil[a_\theta]{\varphi,\vec\psi, \mat \Xi}{\varphi,\vec\psi, \mat \Xi} 
      \\
      \begin{aligned}
        & \geq
        \Norm{ \grad \varphi - \vec \psi }_{\leb2(\W)} ^2
        +
        \qp{\maxi{1}{2 \Norm{\mat A}_{\leb{\infty}(\W)}^2}} \Norm{ \mat M }_{\leb2(\W)} ^2
        +
        \Norm{\rot \vec \psi}_{\leb2(\W)}^2
        \\
        & \quad +
        \Norm{
          \linoptheta(\varphi, \vec \psi, \D\vec \psi - \mat M)
        }_{\leb2(\W)}^2 
        \\
        & \geq
        \bigg(
        \qpsqrt{
          \Norm{ \grad \varphi - \vec \psi}_{\leb2(\W)} ^2 
          +
          \Norm{\rot \vec \psi}_{\leb2(\W)}^2
          +
          \Norm{ \linoptheta(\varphi, \vec \psi,\D\vec \psi) }_{\leb2(\W)} ^2
        }
        \\
        &\phantom{\geq\Big(}
        -\Norm{ \mat A \frobinner \mat M }_{\leb2(\W)} 
        \bigg)^2 
        +
        \qp{\maxi{1}{2 \Norm{\mat A}_{\leb{\infty}(\W)}^2}} \Norm{ \mat M }_{\leb2(\W)} ^2
        \\
        & 
        \geq
        \frac{ \constref{const:coercivity-theta-ug}}{2}
        \qp{
          \Norm{ \varphi }_{\sobh1(\W)}^2
          +
          \Norm{ \vec \psi }_{\sobh1(\W)}^2} 
        +
        \Norm{ \mat A }_{\leb{\infty}(\W)}^2
        \Norm{ \mat M }_{\leb2(\W)} ^2
        \\
        &
        \geq
        \frac{ \constref{const:coercivity-theta-ug}}{2}
        \qp{
          \Norm{ \varphi }_{\sobh1(\W)}^2
          +
          \Norm{ \vec \psi }_{\sobh1(\W)}^2} 
        \\
        &\phantom\geq
        +
        \qpmini{1}{
            \fracl{\constref{const:coercivity-theta-ug}}{4 \Norm{\mat A}_{\leb{\infty}(\W)}^2}
        }
        \Norm{ \mat A }_{\leb{\infty}(\W)}^2
        \Norm{ \mat M }_{\leb2(\W)} ^2
        .
      \end{aligned}
    \end{multline}}%

    By replacing (\ref{eq:replace}) and using Young's inequality, we infer that 
    \begin{multline}
      \label{ineq:using-young}
      \qp{\maxi{1}{2 \Norm{\mat A}_{\leb{\infty}(\W)}^2}} 
      \abil[a_\theta]{\varphi,\vec\psi, \mat \Xi}{\varphi,\vec\psi, \mat \Xi}
      \\
      \begin{split}
        &\geq
        \frac{ \constref{const:coercivity-theta-ug}}{2}
        \qp{
          \Norm{ \varphi }_{\sobh1(\W)}^2
          +
          \Norm{ \vec \psi }_{\sobh1(\W)}^2} 
        \\
        &\phantom\geq
        +
        \qp{
          \mini{1}{
            \frac{\constref{const:coercivity-theta-ug}}{4 \Norm{\mat A}_{\leb{\infty}(\W)}^2}
        }}
        \Norm{\mat A}_{\leb{\infty}(\W)}^2
        \qp{ \frac{1}{2}\Norm{ \mat \Xi }_{\leb2(\W)} ^2
          -
          \Norm{ \D\vec \psi}_{\leb2(\W)} ^2}
        \\
        &\geq
        \qp{
          \mini{ \frac{\constref{const:coercivity-theta-ug}}{8}
          }{\frac{ \Norm{\mat A}_{\leb{\infty}(\W)}^2}{2}}
        }
        \qp{
          \Norm{ \varphi }_{\sobh1(\W)}^2
          +
          \Norm{ \vec \psi }_{\sobh1(\W)}^2
          +
          \Norm{ \mat \Xi }_{\leb2(\W)} ^2
        }.     
      \end{split}
    \end{multline}
    Dividing both sides of (\ref{ineq:using-young}) by
    $\maxi{1}{2\Norm{\mat A}_{\leb{\infty}(\W)}^2}$ establishes the
    claim.%
\end{Proof}
\subsection{Continuity of $a_\theta$}
\label{sec:continuity-a_theta}
We now look at the continuity of $a_\theta$ on $\fatlinspace Y$,
 which includes $\fatlinspace w$.

{
Following \citet{CostabelDauge:99:article:Maxwell}, but for any $d$,
any $\vec\psi,\vec\psi'\in\sobh1(\W; \R d)$ we have
as revealed from
\eqref{eq:curl-and-skew-symmetric-Jacobian},
\eqref{eq:inner-equals-frobinner-half}
and basic Frobenius inner product algebra
that
\begin{equation}
  \begin{split}
    \qp{\rot \vec\psi}\inner\qp{\rot \vec\psi'}
    &
    =
    \frac12{\D\times\vec\psi}\frobinner{\D\times\vec\psi'}
    =
    \frac12\qp{
      \D\vec\psi-\D\vec\psi\transposed
    }\frobinner\qp{
      \D\vec\psi'-\D\vec\psi'\transposed}
    \\
    &
    =
    \D\vec\psi \frobinner \D\vec\psi'
    -
    \D\vec\psi \frobinner\transposeof{\D\vec\psi'}
  \end{split}
\end{equation}
}%
The following inequality follows
\begin{equation}
  \label{eq:curl-inner-Frobenius-identity:Costabel-Dauge}
  \ltwop{\rot \vec\psi}{\rot \vec\psi'}
  \leq
  2 \Norm{\D\vec\psi}_{\leb2(\W)} \Norm{\D\vec\psi'}_{\leb2(\W)}.
\end{equation}
By using \CBS, we realize that
\begin{equation}
  \begin{aligned}
    \big|%
    &
    \ltwop{
      \grad \varphi  - \vec\psi}{
      \grad \varphi'-\vec\psi' }
    +
    \ltwop{
      \D\vec\psi - \mat\Xi}{
      \D \vec\psi'-\mat\Xi' }
    +
    \ltwop{ \rot \vec \psi}{ \rot \vec\psi' }
    \\
    &
    +
    \ltwop{ \linoptheta(\varphi , \vec\psi,\mat\Xi)
    }{
      \linoptheta(\varphi', \vec\psi',\mat\Xi')
    }
    \big|
    \\
    & 
    \leq
    \Norm{\grad\varphi-\vec\psi}_{\leb2(\W)}
    \Norm{\grad\varphi'-\vec\psi'}_{\leb2(\W)}
    +
    \Norm{\D\vec\psi-\mat\Xi}_{\leb2(\W)}
    \Norm{\D\vec\psi'-\mat\Xi'}_{\leb2(\W)}
    \\
    &
    +
    \Norm{\rot \vec \psi}_{\leb2(\W)} \Norm{\rot \vec\psi'}_{\leb2(\W)}
    +
    \Norm{\linoptheta(\varphi , \vec\psi,\mat\Xi)}_{\leb2(\W)}
    \Norm{\linoptheta(\varphi', \vec\psi',\mat\Xi') }_{\leb2(\W)}
    \\
    &
    \leq
    \qp{
      \Norm{\grad\varphi-\vec\psi}_{\leb2(\W)} 
      +
      \Norm{\D\vec\psi-\mat\Xi}_{\leb2(\W)}
      +
      \Norm{\rot\vec\psi}_{\leb2(\W)}
      +
      \Norm{\linoptheta(\varphi,\vec\psi,\mat\Xi)}_{\leb2(\W)}
    }
    \\
    &  
    \times\qp{
      \Norm{\grad \varphi' - \vec\psi'}_{\leb2(\W)}
      +
      \Norm{\D \vec\psi' - \mat\Xi' }_{\leb2(\W)}
      +
      \Norm{\rot \vec\psi'}_{\leb2(\W)}
      +
      \Norm{\linoptheta(\varphi', \vec\psi',\mat\Xi') }_{\leb2(\W)}
    }
    \\ 
    & 
    \leq
    \Big(
    \Norm{\grad \varphi}_{\leb2(\W)} + \Norm{\vec \psi}_{\leb2(\W)}
    + 
    \Norm{\D \vec \psi }_{\leb2(\W)} + \Norm{\mat \Xi}_{\leb2(\W)}
    + 
    \Norm{\rot \vec \psi}_{\leb2(\W)}
    \\
    &
    +
    \Norm{\mat A \frobinner \mat \Xi}_{\leb2(\W)}
    +
    \theta \Norm{\vec b\inner\vec \psi}_{\leb2(\W)}
    +
    (1-\theta) \Norm{\vec b\inner\grad \varphi}_{\leb2(\W)}
    +
    \Norm{c\varphi}_{\leb2(\W)}
    \Big)
    \\
    &
    \times
    \Big(
    \Norm{\grad \varphi'}_{\leb2(\W)} + \Norm{\vec\psi'}_{\leb2(\W)}
    + 
    \Norm{\D \vec\psi' }_{\leb2(\W)} + \Norm{\mat\Xi'}_{\leb2(\W)}
    + 
    \Norm{\rot \vec\psi'}_{\leb2(\W)}
    \\
    &
    +
    \Norm{\mat A \frobinner \mat\Xi'}_{\leb2(\W)}  
    +
    \theta \Norm{\vec b\inner\vec\psi'}_{\leb2(\W)}
    +
    (1-\theta) \Norm{\vec b\inner\grad \varphi'}_{\leb2(\W)}
    +
    \Norm{c\varphi'}_{\leb2(\W)}
    \Big) 
    \\
    &  
    \leq
    \constref{const:continuity-of-atheta}
    \Norm{(\varphi,\vec \psi,\mat \Xi)}_{\fatlinspace{Y}} 
    \Norm{(\varphi',\vec\psi',\mat\Xi')}_{\fatlinspace{Y}}.
  \end{aligned}
\end{equation}
where we introduce the \indexemph{continuity constant}
\begin{multline}
  \constdef[\W,\linop L,\theta]%
      {const:continuity-of-atheta}
    :=
    5\Big(
    \maxi{
      \Norm{c}_{\leb\infty(\W)}
    }
    {
      \maxi{
        \qp{1+d(1-\theta)\Norm{\vec b}_{\leb\infty(\W)}}
      }
      \\
      {
        \maxi{
          \qp{1+d \theta \Norm{\vec b}_{\leb\infty(\W)}}
        } 
      {
        \maxi{
         \qp{1+ \sqrt{2}}
        }{
          \qp{1+d^2\Norm{\mat A}_{\leb\infty (\W)}} 
        }
      }
    } 
    }
\Big)^2.
\end{multline}
We have thus established that 
\begin{equation}
  \label{eq:continuity-a_theta}
  \norm{\abil[a_\theta]{
      \varphi ,\vec\psi, \mat\Xi
    }{
      \varphi', \vec\psi', \mat\Xi'
    }}
  \leq
  \constref{const:continuity-of-atheta}
  \Norm{ (\varphi,\vec \psi, \mat \Xi)}_{\fatlinspace Y} 
  \Norm{ (\varphi',\vec\psi', \mat\Xi')}_{\fatlinspace Y}.
\end{equation}
By the same argument, we can also show the continuity of
$\hat{a}_\theta$ on $\sobhz1(\W) \times \sobh1(\W; \R d)$, 
which includes $\sobhz1(\W) \times \linspace v$. The continuity
of $a_\theta$ on $\fatlinspace w$ and Theorem~\ref{the:coercivity-u,g,H} 
imply that the problem~(\ref{eq:Euler-Lagrange}) is well-posed
and also, the continuity of $\hat{a}_\theta$ on $\sobhz1(\W)
\times \linspace v$ and Theorem \ref{the:coercivity-u,g} imply that the 
problem~(\ref{eq:E-L-2}) is well-posed.
\subsection{Necessity of the zero tangential trace condition}
\label{relaxing-problem}
If we define the functional $ \tilde{E}_\theta$ on $\fatlinspace y$ by
\begin{equation}
  \label{functional:minimize:relaxed}
    \tilde E_\theta(\varphi,\vec\psi,\mat\Xi)
    :=
    \Norm{ \grad\varphi-\vec\psi}_{\leb2(\W)} ^2
    +
    \Norm{ \D\vec\psi-\mat\Xi}_{\leb2(\W)} ^2
    + 
    \Norm{ \linoptheta(\varphi, \vec\psi,\mat\Xi)  -f } _{\leb2(\W)}^2,
\end{equation}
as a straightforward alternative to $E_\theta$,
it still provides equivalence between the minimization 
problem and the strong solution of (\ref{eq:nondivergence0}).
Nonetheless additional conditions on the space, e.g.,
zero-tangential-trace assumption for the field-space (containing
$\vec{g}$ and $\vec\psi$) and the functional, e.g.,
the extra term $\Norm{\rot\vec\psi}_{\leb2(\W)}^2$
in (\ref{functional:minimize}) provide coercivity for $E_\theta$ which 
may fail for $\tilde E_\theta$.
\par
To illustrate how
$E_\theta(\phi,\vec\psi,\mat\Xi)$'s coercivity may fail when its
second argument $\vec\psi$ is a generic element of $\sobh1(\W)^d$
with nonzero tangential trace, take $\mat A = \eye$, $\vec
b = 0$, $c = 0$ and consider
$\varphi = 0$, $\mat\Xi=\D\vec\psi$. Let us show that
\begin{equation}
  \label{ineq:contradict-claim}
  \Norm{\vec\psi}_{\leb2(\W)} ^2
  + 
  \Norm{ \rot\vec\psi }_{\leb2(\W)}^2
  +
  \Norm{\div \vec \psi}_{\leb2(\W)}^2
  \geq
  C \Norm{\vec \psi}_{\sobh1(\W)}^2
\end{equation}
is not always satisfied on $\sobh1(\W)^d$.

In this regard, let $\seqsinat q n$ be a sequence 
in $\sobh{1/2}(\boundary\W)$ with $\int_{\boundary\W} q_n =0$, 
satisfying $\lim_{n\to\infty}{\Norm{q_n}_{\sobh{1/2}(\boundary\W)}}=\infty$ and
$\Norm{q_n}_{\sobh{-1/2}(\boundary\W)}\leq{C}$ bounded, uniformly in
$n$. Obviously, for each $n \in \naturals$, problem of finding $v_n \in \sobh2(\W)$ with 
$\int_{\W}v_n=0$ such that
\begin{equation}
  \label{eq:counterexample}
  \lap v_n = 0
  \tand
  \boundarytraceof\W{\normalderto\W{v_n}}%
  =q_n,
\end{equation}
is well-posed. Stability of $v_n$ and the trace theorem imply that there exist
constants $\constref[1]{eq:trace-constants:counterexample}$ and
$\constref[2]{eq:trace-constants:counterexample}$ such that
\begin{equation}
  \label{eq:trace-constants:counterexample}
  \Norm{v_n}_{\sobh1(\W)} 
  \leq 
  \constref[1]{eq:trace-constants:counterexample}
  \Norm{q_n}_{\sobh{-1/2}(\boundary\W)}
  \tand
  \Norm{q_n}_{\sobh{1/2}(\boundary\W)}
  \leq
  \constref[2]{eq:trace-constants:counterexample}
  \Norm{v_n}_{\sobh2(\W)}.
\end{equation}
Our assumptions on $\seqsinat qn$ thus imply that 
\begin{equation}
  \label{ineq:contradiction}
  \Norm{v_n} _{\sobh1(\W)} \leq\constref[1]{eq:trace-constants:counterexample}C
  \tand
  \lim_{n \rightarrow \infty} \Norm{\D^2 v_n}_{\leb2(\W)} = \infty.
\end{equation}
By setting $\vec \psi_n = \grad v_n$, it is clear that
\begin{equation}
  \vec \psi_n \in \sobh1\qp{\W;\R d}, 
  \quad
  \div \vec \psi_n = 0,
  \quad 
  \rot \vec \psi_n = 0.
\end{equation}
Now by replacing $ \vec \psi_n$ in (\ref{ineq:contradict-claim}) and
taking the limit $n \rightarrow \infty$ of both sides,
(\ref{ineq:contradiction}) makes a contradiction.

This example shows also that:
\begin{itemize}
\item
  Lemma~\ref{lem:Miranda-Talenti-special-Cordes} and consequently
  Theorem~\ref{the:coercivity-u,g} and \ref{the:coercivity-u,g,H}
  are not valid without the zero tangential-trace condition on $\vec\psi$;
\item 
  coercivity is merely sufficient, not necessary, for the unique
  minimization of $E_\theta$ and the solvability of
  (\ref{eq:nondivergence0}), because $\tilde E_\theta$ 
  also takes the minimum value at $(u, \grad u , \D^2u)$. %
\end{itemize}
\subsection{Nonzero boundary values}
\label{sec:nonzero-boundary-values}
Since in problem~(\ref{eq:nondivergence0-inhomogeneous}), when
heterogeneous, i.e., $\restriction u{\boundary\W}=r \neq 0$, a full
extension of $r$ to all of $\closure\W$ may not be explicitly
available while its approximation must be sought numerically or built
into the discrete solution space. In this case, a reasonable solution is 
to use the following extension of the functional $E_\theta$ (which
we call the same) on $\fatlinspace Y$
\begin{equation}
  \begin{aligned}
    E_\theta (\varphi,\vec\psi,\mat\Xi)
    := &
    \Norm{ \grad\varphi-\vec\psi}_{\leb2(\W)} ^2
    +
    \Norm{ \D\vec\psi-\mat\Xi}_{\leb2(\W)} ^2
    +
    \Norm{ \rot \vec \psi }_{\leb2(\W)} ^2
    \\
    & 
    +
    \Norm{ \linoptheta(\varphi, \vec\psi,\mat\Xi)  -f } _{\leb2(\W)}^2
    +
    \Norm{ \varphi-r }_{\leb2(\boundary\W)}^2,
  \end{aligned}
\end{equation} 
and then considering the Euler--Lagrange equation of the minimization problem 
\begin{equation}
  \label{eq:minimization-nonhomogeneous}
  (u, \vec g, \mat  H)
  =\underset
  {\substack{
      \triple\varphi{\vec\psi}{\mat\Xi}\in\fatlinspace Y
  }}
  \argmin E_\theta(\varphi ,\vec\psi , \mat\Xi) .
\end{equation}
It is easy to check that (\ref{eq:minimization-nonhomogeneous}) and
the problem of finding strong solution to
(\ref{eq:nondivergence0-inhomogeneous}) are equivalent. Although the
setting of proving coercivity of the bilinear form corresponding to
(\ref{eq:minimization-nonhomogeneous}) is no longer provided, we would
like to point out that coercivity of the bilinear form is not
necessary to establish that the problem is well-posed.
\section{A conforming Galerkin finite element method}
\label{sec:a-conforming-Galerkin-FEM}
In this section, we derive via a Galerkin approach, discrete
counterparts of the infinite dimensional problems of
\S~\ref{sec:coercivity-continuity-of-cost-functional}; we specifically
use conforming Galerkin finite elements where the finite dimensional
subspace of the functional spaces $\fatlinspace{w}$ or
$\fatlinspace{y}$. Using first an abstract choice of Galerkin subspaces
and the coercivity of the exact problem we derive abstract a priori
error estimates in Theorem~\ref{pro:error-estimates}.

{We analyze the method and the well-posed nature of the problem
with zero boundary condition, i.e., problem
(\ref{eq:discrete-system-finite-dim}), but we will use a nonhomogenous
boundary value problem (\ref{eq:discrete-system-nonhomogeneous}) in
the numerical tests of \S~\ref{test-problem-nonzero-boundary}, the
numerical results is as good as zero boundary problem. Since
coercivity on a normed space is inherited by its subspaces, thanks to
Theorem~\ref{the:coercivity-u,g,H}, the resulting discrete problems
(\ref{eq:discrete-system-finite-dim}) are automatically well posed.}

We realize the abstract results into concrete theorems by
introducing a conforming finite element discretization and discuss
about how well a solution may be approximated by proposed method. We
provide an \aposteriori error estimate, with fully computable
estimators, via the plain residual provided by the least-squares
functional in Theorem \ref{the:residual-error-bound}, as well as an
\apriori error bound in Theorem \ref{the:convergence-rate}. Finally
we use the \aposteriori error indicators to design Algorithm
\ref{alg:adaptive} for adaptive mesh refinement based on the by-now
classical loop of the form
  $\text{solve}\to\text{estimate}\to\text{mark}\to\text{refine}.$

{\emph{We like to remind the reader of Remark~\ref{rem:equivalence}
 implying we \emph{always} have $\vec{g}=\grad{u}$ and $\mat{H}=\D\vec{g}=\D^2u$}.
}
\subsection{An abstract discrete problem}
\label{def:abstract-discrete-problem}
Consider finite dimensional subspaces (to be specified later)
satisfying
\begin{equation}
  \tilde{\fespace{U}}\subset\sobh1(\W),
  \:
  \tilde{\fespace{G}}\subset\sobh1(\W;\R{d})
  \tand
  \fespace{H}\subset\leb2(\W;\Symmatrices{d}).
\end{equation}
\index{\(\tilde{\fespace{U}}\)}
\index{\(\tilde{\fespace{G}}\)}
Set
\begin{equation}
  \fespace{U}:=\tilde{\fespace U}\meet\sobhz1(\W)
  \tand
  \fespace{G}:=\tilde{\fespace G}\meet\linspace{v},
\end{equation}
\index{\(\fespace u\)}
\index{\(\fespace g\)}
and define the Galerkin spaces
\begin{equation}
  \fespace V:=\fespace{U}\times\fespace{G} \times \fespace{H},
  \:
  \fespace{x}
  :=\fespace{u}\times\tilde{\fespace{g}}\times\fespace{h}
  \tand
  \fespace Y:=\tilde{\fespace{U}}\times\tilde{\fespace{G}} \times \fespace{H}.
\end{equation}
\index{\(\fespace v\)}\index{\(\fespace x\)}\index{\(\fespace y\)}
We consider the discrete counterpart of (\ref{eq:Euler-Lagrange})
consisting in finding $(\fespacefun vu,\vecfespacefun vg,\matfespacefun vh)\in\fespace{V}$
such that
\begin{equation}
  \label{eq:discrete-system-finite-dim}
  \abil[a_\theta]{
    \fespacefun vu,\vecfespacefun vg, \matfespacefun vH
  }{
    \fe\varphi , \vecfe\psi, \matfe\Xi
  }
  =
  \ltwop{
    f
  }{
    \linoptheta(\fe\varphi,\vecfe\psi,\matfe\Xi)
  }
  \Foreach (\fe\varphi, \vecfe\psi, \matfe\Xi) \in {\fespace[]V},
\end{equation}
which we will analyze in this section; the analogue on the space
$\fespace x$ replacing $\fespace v$ denoted $\triple{\fespacefun
  xu}{\vecfespacefun xg}{\matfespacefun xh}$ will be used in
\S~\ref{sec:numerical-experiments}.
\par
To treat possible nonzero boundary values $r$ we also consider the
discrete problem of finding
$(\fespacefun{y}u,\vecfespacefun{y}g,\matfespacefun{y}h)\in\fespace{Y}$
such that
{
  \begin{multline}
    \label{eq:discrete-system-nonhomogeneous}
    \abil[a_\theta]{
      \fespacefun yu ,\vecfespacefun{y}g, \matfespacefun yh
    }{
      \fe\varphi , \vecfe\psi, \matfe\Xi
    }
    +
    \qa{\fespacefun yu,\fe\varphi}_{\boundary\W}
    =
    \qa{r,\fe\varphi}_{\boundary\W}
    +
    \ltwop{
      f
    }{
      \linoptheta(\fe\varphi,\vecfe\psi,\matfe\Xi)
    }
    \\
    \Foreach (\fe\varphi, \vecfe\psi, \matfe\Xi) \in {\fespace[]Y}.
  \end{multline}
  \begin{Obs}[our approach vs. standard FEM]
    \label{obs:our-approach-vs-standard-FEM}
    Strictly speaking our approach here does not extend the classical
    finite element approach but should be viewed as a variant.  We only
    test the boundary value $r$ with $\fe\varphi$ while the rest of the
    equation is tested with
    $\linoptheta(\fe\varphi,\vecfe\psi,\matfe\Xi)$.  Thus even letting
    $\mat A=\eye, \vec b=0, c=0$ we do not get the standard Poisson
    solver arising from its weak formulation, since we work with strong
    formulation and do not integrate by parts the Laplacian term.
  \end{Obs}
}
\begin{The}[quasi-optimality]
\label{pro:error-estimates}
Consider
$(\fespacefun{v}u,\vecfespacefun{v}g,\matfespacefun{v}H)  \in \fespace[]V$ 
is the unique solution of discrete problem (\ref{eq:discrete-system-finite-dim}). 
It satisfies the error estimate 
\begin{equation}
  \label{eqn:apriori-error-bound}
  \Norm{ (u,\grad u, \D^2 u)-(\fespacefun{v}u , \vecfespacefun{v}g , \matfespacefun{v}H) }_{\fatlinspace Y}
  \leq
  \frac{\constref{const:continuity-of-atheta}}{\constref%
    {const:coercivity-theta-ugH}}
  \inf_{({\fe\varphi},{\vecfe\psi},\matfe\Xi) \in {\fespace[]V}}
  \Norm{ (u,\grad u, \D^2 u)-({\fe\varphi},{\vecfe\psi}, \matfe\Xi)}_{\fatlinspace Y}.
\end{equation}
where $\constref{const:coercivity-theta-ugH}$ and
$\constref{const:continuity-of-atheta}$ respectively are the coercivity
and the continuity constants of $a_\theta$ relative to $\fatlinspace
V$.
\end{The}
\begin{Proof}
  It is easy to check that the following Galerkin orthogonality relation holds
  \begin{equation}
    a_\theta((u,\vec g, \mat H)
    -
    (\fespacefun{v}u,\vecfespacefun{v}g,\matfespacefun{v}H)
    , 
    (\fe\varphi,{\vecfe\psi}, {\matfe\Xi})
    -
    (\fespacefun{v}u ,\vecfespacefun{v}g , \matfespacefun{v}H))
    = 0
    ~ \Foreach (\fe\varphi,{\vecfe\psi}, {\matfe\Xi}) \in {\fespace[]V}.
  \end{equation}
 Therefore, for any $(\fe\varphi,{\vecfe\psi}, {\matfe\Xi}) \in {\fespace[]V}$,
 we get
  \begin{equation}
    \begin{aligned}
      &a_\theta((u,\vec g, \mat H)-(\fespacefun{v}u ,\vecfespacefun{v}g , \matfespacefun{v}H)
      , 
      (u,\vec g, \mat H)-(\fespacefun{v}u ,\vecfespacefun{v}g , \matfespacefun{v}H))
      \\
      &\qquad = 
      a_\theta(
      (u,\vec g, \mat H)-(\fespacefun{v}u ,\vecfespacefun{v}g , \matfespacefun{v}H)
      , 
      (u,\vec g, \mat H) -
      (\fe\varphi,{\vecfe\psi}, {\matfe\Xi})
      +
      (\fe\varphi,{\vecfe\psi}, {\matfe\Xi})
      -
      (\fespacefun{v}u ,\vecfespacefun{v}g , \matfespacefun{v}H))
      \\
      &\qquad =
      a_\theta((u,\vec g, \mat H)-(\fespacefun{v}u ,\vecfespacefun{v}g , \matfespacefun{v}H)
      , 
      (u,\vec g, \mat H)-(\fe\varphi,{\vecfe\psi}, {\matfe\Xi}))
      \\
      & \qquad 
      +
      a_\theta((u,\vec g, \mat H)-(\fespacefun{v}u ,\vecfespacefun{v}g , \matfespacefun{v}H)
      , 
      (\fe\varphi,{\vecfe\psi}, {\matfe\Xi})
       - 
      (\fespacefun{v}u ,\vecfespacefun{v}g , \matfespacefun{v}H))
      \\
      &\qquad =
      a_\theta((u,\vec g, \mat H)-(\fespacefun{v}u ,\vecfespacefun{v}g , \matfespacefun{v}H)
      ,
      (u,\vec g, \mat H) - (\fe\varphi,{\vecfe\psi}, {\matfe\Xi}) )
      .
    \end{aligned}
  \end{equation}
 Coercivity (\ref{eqn:coercivity-u,g,H}) and continuity (\ref{eq:continuity-a_theta}) 
 imply that %
 we have
  \begin{equation}
    \begin{aligned}
      \big\|%
        (u,\vec g,& \mat H)-(\fespacefun{v}u ,\vecfespacefun{v}g , \matfespacefun{v}H)
      \big\|%
      _{\fatlinspace Y}^2
      \\
      &
      \leq \constref%
           {const:coercivity-theta-ugH}^{-1} 
      a_\theta((u,\vec g, \mat H)-(\fespacefun{v}u ,\vecfespacefun{v}g , \matfespacefun{v}H)
      , 
      (u,\vec g, \mat H)-(\fespacefun{v}u ,\vecfespacefun{v}g , \matfespacefun{v}H))
      \\
      &%
      =
      \constref%
               {const:coercivity-theta-ugH}^{-1} 
      a_\theta((u,\vec g, \mat H)-(\fespacefun{v}u ,\vecfespacefun{v}g , \matfespacefun{v}H)
      , 
      (u,\vec g, \mat H)-(\fe\varphi,{\vecfe\psi}, {\matfe\Xi}))
      \\
      &%
      \leq
      \frac{\constref{const:continuity-of-atheta}}{\constref%
        {const:coercivity-theta-ugH}}
      \Norm{(u,\vec g, \mat H)-(\fespacefun{v}u ,\vecfespacefun{v}g , \matfespacefun{v}H)}_{\fatlinspace Y}
      \Norm{(u,\vec g, \mat H)-(\fe\varphi,{\vecfe\psi}, {\matfe\Xi})}_{\fatlinspace Y}.
    \end{aligned}
  \end{equation}
  Dividing both sides by 
  $\Norm{(u,\vec g, \mat H)-(\fespacefun{v}u ,\vecfespacefun{v}g , \matfespacefun{v}H)}_{\fatlinspace Y}$
  yields the assertion.
\end{Proof}
\begin{The}[error-residual \aposteriori estimates]
\label{the:residual-error-bound}
Suppose that\\$(\fespacefun{v}u ,\vecfespacefun{v}g , \matfespacefun{v}H) \in \fespace[]V$ 
is the unique solution of the discrete problem (\ref{eq:discrete-system-finite-dim}). 
  \begin{enumerate}[(i)]
  \item
    \label{item:aposteriori-residual:upper-bound}
    The following \indexen{a posteriori residual upper bound} %
    holds 
    \begin{multline}
      \label{eqn:aposteriori-residual:upper-bound}
        \Norm{ (u, \grad u, \D^2 u)-(\fespacefun{v}u ,\vecfespacefun{v}g , \matfespacefun{v}H) }_{\fatlinspace Y}^2
        \leq
        \constref%
            {const:coercivity-theta-ugH}^{-1}
        \Big(%
          \Norm{ \grad \fespacefun{v}u - \vecfespacefun{v}g  }_{\leb2(\W)}^2
          +
          \Norm{ \D\vecfespacefun{v}g  - \matfespacefun{v}H }_{\leb2(\W)}^2
          \\
          \phantom{%
            \leq
            \constref%
                     {const:coercivity-theta-ugH}^{-1}
            \Big(%
          }
          +
          \Norm{ \rot \vecfespacefun{v}g  }_{\leb2(\W)}^2
          +
          \Norm{ \linoptheta(\fespacefun{v}u , \vecfespacefun{v}g , \matfespacefun{v}H) -f }_{\leb2(\W)}^2 
        \Big)%
        .
    \end{multline}
  \item
    \label{item:aposteriori-residual:lower-bound}
    For each open subdomain $\w \subseteq \Omega$ we have 
    \begin{equation}
      \label{eqn:aposteriori-residual:lower-bound}
      \begin{aligned}
        \|
        \grad \fespacefun{v}u 
        &
        - \vecfespacefun{v}g 
        \|%
        _{\leb2(\w)}^2
        +%
        \Norm{ \D\vecfespacefun{v}g  - \matfespacefun{v}H }_{\leb2(\w)}^2
        \\
        &
        +
        \Norm{ \rot \vecfespacefun{v}g  }_{\leb2(\w)}^2
        +
        \Norm{ \linoptheta(\fespacefun{v}u , \vecfespacefun{v}g , \matfespacefun{v}H) -f }_{\leb2(\w)}^2  
        \\
        &%
        \leq
        \mathchanges{\constref[\w]{const:continuity-of-atheta:local}}
        \Big(%
          \Norm{u-\fespacefun{v}u}_{\sobh1{(\w)}}^2
	  +
	  \Norm{\grad u-\vecfespacefun{v}g }_{\sobh1{(\w)}}^2
	  +
	  \Norm{\D^2 u-\matfespacefun{v}H}_{\leb2{(\w)}}^2
	  \Big)%
          ,
      \end{aligned}
    \end{equation}
    {where
      \begin{equation}
        \constdef[\w,\linop L,\theta]{const:continuity-of-atheta:local}
        :=\constref[\w,\linop L,\theta]{const:continuity-of-atheta}
      \end{equation}
    }%
    is the continuity constant of the analogue of $a_\theta$ on
    the space $\fatlinspace Y$ (without boundary values)
    albeit over $\w$ instead of $\Omega$ defined in
    \eqref{const:continuity-of-atheta}. 
  \end{enumerate}
\end{The}
\begin{Proof}
  The coercivity of $a_\theta$ from Theorem~\ref{the:coercivity-u,g,H}
  immediately implies the \aposteriori residual--error upper bound
  (\ref{eqn:aposteriori-residual:upper-bound}). The continuity of
  $a_\theta$, in view of (\ref{eq:continuity-a_theta}) on
  $\fatlinspace{y}$ albeit with $\W$ replaced by a subset $\w$,
  implies (\ref{eqn:aposteriori-residual:lower-bound}).
\end{Proof}
\subsection{Triangulations and finite element spaces}
Let $\mathfrak T$ \index{$\mathfrak T$} be a 
collection of \indexen{conforming simplicial partitions}, also known as
\indexen{meshes}.
For each mesh $\mesh t$ in $\mathfrak T$ the domain $\W\subseteq\R{d}$ such that
\begin{equation}
  \overline{\W}=\overline{\W}_{\mesh t}
  :=
 \unions K{\mesh t}K
  ,
\end{equation}
which requires $\W$ to be a polyhedral domain. If $\W$ is not
polyhedral, it is necessary to approximate pieces of $\boundary\W$ by
(possibly curved) simplex sides, which can give rise to simplices
having curved sides and isoparametric elements; for simplicity, we do
not treat the details of this more general case in this work, although
many parts can be modified to include it.

For each element $K\in\mesht\in\mathfrak T$, denote
$\meshsize[K]:=\diam{K}$, $\rho_K$ be the lowest upper bound on the
radius of a ball contained in $K$, and $\sigma(K):=h_K/\rho_K$ its
(inverse) \indexemph{shape-regularity} or \indexemph{chunkiness
  parameter} as in \citet{BrennerScott:08:book:The-mathematical},
which we follow for many notations and results herein. We define
$\sigma(\mesht):=\max_{K\in\mesht}\sigma(K)$ and
$\sigma(\mathfrak{T}):=\sup_{\mesht\in\mathfrak T}\sigma(\mesht)$ and
we assume that this is a strictly positive finite real number.
Finally denote by
$\meshsize:=\meshsize[\mesht]:=\max_{{K}\in\mesht}h_{K}$ the
\indexemph{mesh-size function} defined on all of $\W$ (although the
meshsize $\meshsize[\mesht]$ depends on $\mesht$ we drop this
dependence and use $\meshsize$ to lighten notation).
Consider the following concrete realization of the Galerkin finite element
spaces defined in \S~\ref{def:abstract-discrete-problem}
\begin{gather}
  \label{def:concrete-Galerkin-spaces}
  \begin{aligned}
    \tilde{\fespace u}
    :=
    \poly{k}\qp{\mesht}\meet\sobh1(\W)
    ,
    \qquad
    \fespace u
    &
    :=\tilde{\fespace u}\meet\sobhz1(\W)
    ,
    \\
    \:
    \tilde{\fespace{g}}
    :=
    \poly{k}\qp{\mesht;\R d}\meet\sobh1(\W;\R{d})
    ,
    \qquad
    \fespace g
    &
    :=
    \tilde{\fespace g}
    \meet
    \linspace v
    ,
  \end{aligned}
  \intertext{%
    and }
  \fespace h
  :=
  \poly{k-1}\qp{\mesht;\Symmatrices d}.%
\end{gather}
Denote by $\pinteron u$ and $\pinteron g$ %
  a corresponding \indexemph{nodal interpolators}.
\begin{Lem}[intepolation error estimates]
  \label{Lem:interpolation-error}
Let $\mesht$ be in a collection $\mathfrak T$ of
  shape-regular conforming simplicial meshes on the polyhedral domain
  $\W\subseteq\R d$. For each of $X=\reals$ or $\R d$, %
  consider the space
 \begin{equation}
   \label{def:interpolation-space}
    \begin{split}
      \fespace W:=\poly{k}\qp{\mesht;X}\meet \sobh1(\W;X).
    \end{split}
  \end{equation}
  For any $\varphi \in \sobh{s}(\W;X)$ with $1 
 \leq s \leq k+1$, suppose that $\pinteron w \varphi$ denotes
 nodal interpolation of $\varphi$ in $\fespace W$. Then
 there exists $\constref{eqn:interpolation-error}>0$, which depends
 on the shape-regularity of $\mesht$, such that
 \begin{equation}
   \label{eqn:interpolation-error}
   \Norm{\varphi - \mathcal{I}_{\fespace W}\varphi}_{\sobh1(\W)}
   \leq\constref{eqn:interpolation-error}h^{s-1} 
   \Norm{\varphi}_{\sobh{s}(\W)}
           \quad \text{for }  %
   0<h\leq1
   .
 \end{equation}
\end{Lem}
\begin{Proof}
  This is a standard result \citep[Th.4.4.20]{BrennerScott:08:book:The-mathematical}.
\end{Proof}
%
\begin{comment}
  With respect to the introduced Galerkin finite element spaces, we 
  consider the discrete problem corresponding to zero boundary condition, as finding 
  $(u_h, \vec g_h, \mat H_h) \in \fespace u \times \fespace g \times \fespace h$ 
  such that
  \begin{equation}
    \label{eq:discrete-system-Galerkin}
    \abil[a_\theta]{
      u_h ,\vec g_h, \mat H_h
    }{
      \varphi_h , \vec\psi_h, \mat\Xi_h
    }
    =
    \ltwop{
      f
    }{
      \linoptheta(\varphi_h, \vec\psi_h ,\mat\Xi_h)
    }
    ~
    \Foreach (\varphi_h, \vec\psi_h, \mat\Xi_h) \in \fespace u \times \fespace g \times \fespace h.
  \end{equation}
  The following theorem states the main result of this section.
\end{comment}
%
\begin{The}[a priori error estimate]
  \label{the:convergence-rate}
  Suppose the collection of meshes $\mathfrak T$ satisfies the assumptions
  of Lemma~\ref{Lem:interpolation-error}, that the strong solution $u$
  of (\ref{eq:nondivergence0}) satisfies $u \in \sobh{\alpha+2}(\W)$,
  for some real $0< \alpha\leq k$ and let $(\fespacefun vu,\vecfespacefun
  vg,\matfespacefun vh)\in\fespace v=\fespace u\times\fespace
  g\times\fespace h$ be the finite element solution of
  (\ref{eq:discrete-system-finite-dim})
  relative to the mesh $\mesht$, i.e., with the choice of spaces~(\ref{def:concrete-Galerkin-spaces}) .
  Then for some $\constref{eqn:convergence-rate}>0$
  independent of $u$ and $h$ we have
  \begin{equation}
    \label{eqn:convergence-rate}
    \Norm{(u,\grad u, \D^2 u)-(\fespacefun vu,\vecfespacefun vg, \matfespacefun vH)}_{\fatlinspace Y} 
    \leq
    \constref{eqn:convergence-rate}
    h^{\alpha}
    \Norm{u}_{\sobh{\alpha+2}(\W)}.
  \end{equation}
\end{The}
\begin{Proof}
  From Lemma \ref{Lem:interpolation-error} 
  we know the interpolation inequalities
  \begin{gather}
    \Norm{u - \pinteron u u}_{\sobh1(\W)}
    \leq
    \constref{eqn:interpolation-error}
    h^\alpha \Norm{u}_{\sobh{\alpha+1}(\W)}
    \leq
    \constref{eqn:interpolation-error}
    h^\alpha \Norm{u}_{\sobh{\alpha+2}(\W)},
    \\
    \Norm{\grad u - \pinteron g \grad u }_{\sobh1(\W)}
    \leq
    \constref{eqn:interpolation-error}
    h^\alpha \Norm{\grad u}_{\sobh{\alpha+1}(\W)}
    \leq
    \constref{eqn:interpolation-error}
    h^\alpha\Norm{u}_{\sobh{\alpha+2}(\W)},
    \\
    \Norm{\D^2 u - \D(\pinteron g \grad u)}_{\leb2(\W)} 
    \leq
    \constref{eqn:interpolation-error}
    h^\alpha
    \Norm{\grad u}_{\sobh{\alpha+1}(\W)}
    \leq 
    \constref{eqn:interpolation-error}
    h^\alpha\Norm{u}_{\sobh{\alpha+2}(\W)},
  \end{gather}
  hence
  \begin{multline}
    {\Norm{ u - \pinteron u  u}_{\sobh1(\W)}^2
      +
      \Norm{\grad u - \pinteron g \grad u}_{\sobh1(\W)}^2
      +
      \Norm{\D^2 u - \D(\pinteron g \grad u)}_{\leb2(\W)}^2}
    \\
    \leq
    3\constref{eqn:interpolation-error}^2
    h^{2\alpha}
    \Norm{u}_{\sobh{\alpha+2}(\W)}^2.
  \end{multline}
  The assertion now follows from Theorem \ref{pro:error-estimates}.
\end{Proof}
\begin{Obs}[curved domain]
In Theorem~\ref{the:convergence-rate}, the domain is assumed polyhedral, so 
that it can be triangulated exactly. If $\W$ has a curved boundary, isoparametric 
finite elements may be used. In isoparametric method, a smooth or piecewise smooth 
boundary, $\boundary\W$, guarantees that the elements with curved boundary 
are not too distorted from triangles. Consequently, an error bound similar to 
that of Lemma~\ref{Lem:interpolation-error} can be established. The final result is that the 
error using isoparametric finite element goes to zero at the same rate as if 
ordinary Lagrange triangles were used on polyhedral domain.
This claim can be found in \citet[\S\S~4.3--4]{Ciarlet:02:book:The-finite}.
\end{Obs}
%
\begin{comment}
\subsection{Conforming finite element discretization for nonzero boundary}
Corresponding to the problem with nonzero boundary condition, define 
also the following Galerkin finite element spaces
\begin{gather}
  \label{def:interpolation-space}
  \tilde{\fespace u}
  :=
  \poly{k}\qp{\mesht}\meet\sobh1(\W),
  \quad
  \tilde{\fespace g}
  :=
  \poly{k}\qp{\mesht;\R d} \meet  \sobh1(\W; \R d).
  \end{gather}
Then we consider the discrete problem of finding 
$(u_h, \vec g_h, \mat H_h) \in \tilde{\fespace u} \times \tilde{\fespace G} \times \fespace h $ 
such that
%
\begin{multline}
  \label{eq:discrete-system-Galerkin-heterogen}
  \abil[a_\theta]{
    u_h ,\vec g_h, \mat H_h
  }{
    \varphi_h , \vec\psi_h, \mat\Xi_h
  } 
  +
  \qa{u_h - r,\varphi_h}_{\boundary\W}
  =
  \ltwop{
    f
  }{
    \linoptheta(\varphi_h, \vec\psi_h ,\mat\Xi_h)
  }
  \\ 
  \Foreach (\varphi_h, \vec\psi_h, \mat\Xi_h) \in \tilde{\fespace u} \times \tilde{\fespace g} \times \fespace h.
\end{multline}
%
%
\end{comment}
%
\subsection{Adaptive mesh refinement strategy}
We close this section by proposing an \indexen{adaptive algorithm} based on the
\aposteriori residual error bounds, Theorem~\ref{the:residual-error-bound}.
Controlling the error of a numerical approximation is prerequisite for
more reliable simulations, while adapting the discretization to local
features of problem can be lead to more efficient simulations. In this
regard, the \aposteriori residual error estimate of Theorem~\ref{the:residual-error-bound}
paves a way to use adaptive refinement approach. By considering
the \indexemph{local error indicator} for each ${K}\in\mesht$
\begin{equation}
  \label{def:error-estimator}
  \begin{split}
    \eta({K})^2
    &
    :=
    \Norm{\grad\fespacefun vu-\vecfespacefun vg}_{\leb2({K})}^2
    +
    \Norm{\D\vecfespacefun vg-\matfespacefun vH}_{\leb2({K})}^2
    \\
    &\phantom{:=}
    +
    \Norm{\rot\vecfespacefun vg}_{\leb2({K})}^2
    +
    \Norm{
      \linoptheta(\fespacefun vu,\vecfespacefun vg,\matfespacefun vH)-f}_{\leb2({K})}^2,             
  \end{split}
\end{equation}
and
\begin{equation}
  \eta^2:={\sum_{{K}\in\mesht}{\eta({K})^2}},
\end{equation}
we track the following adaptive algorithm which we shall test in \S~\ref{sec:numerical-experiments}.
\begin{Alg}[adaptive least squares nondivergence Galerkin solver]
  \label{alg:adaptive}
  { Following is an adaptive mesh refinement algorithm, based
    on the \aposteriori error indicator algorithm pioneered by
    \citet{Dorfler:96:article:A-convergent} and subsequently developed
    into variants by many authors \citep[and references
      therein]{Verfurth:13:book:A-posteriori}.  We use a
    \indexemph{bulk-chasing} \akaindexemph{Dörfler's marking} strategy
    modified as follows: we use sorting and based on a fixed ratio
    $\beta$ of triangles (instead of the fixed ration $\theta$ of
    indicator).  Namely, at each adaptive level $l$ we mark for
    refinement those elements $K$, forming a subset $\mesh m$ of the
    domain's partition $\mesh T_l$, with the highest $\eta(K)$s and of
    cardinality $\card\mesh m=\ceil{\beta\card\mesh T_l}$ (the
    smallest integer bigger than $\beta$ times the cardinality of
    $\mesh T_l$) for some fixed ``element-fraction'', whereas
    \citet{Dorfler:96:article:A-convergent} uses a
    ``indicator-fraction'' (called $\theta$ therein) corresponding to
    a subset $\mesh m\subseteq\mesh T_l$ such that
    $\sum_{K\in\mesh{m}}\eta(K)^2\approx\theta\sum_{K\in\mesh
      T_l}\eta(K)^2$.
  }%
  \begin{algorithmic}[1]
    \REQUIRE{%
    data of Problem~(\ref{eq:nondivergence0}),
    refinement fraction
    $\beta\in\opinter01$, tolerance $\tol$ and maximum number of
    iterations $\maxiter$}
    \ENSURE{sequences $\listidotsfromto {\fe u}0L$, $\listidotsfromto {\vecfe g}0L$, $\listidotsfromto {\matfe h}0L$ of discrete solutions
      of~(\ref{eq:discrete-system-finite-dim}) either with
      $$\Norm{(u,\grad u,\D^2u)-(\fe[L]u,\vecfe[L]g,\matfe[L]h)}_{\fatlinspace y}\leq{\constref{const:coercivity-theta-ugH}^{-1}} \tol$$ 
      or after $\maxiter$ iterations}
    \PROCEDURE{{Adaptive-Least-Squares-Solver}}{$\qp{\W,\mat A,\vec b,c,f,r,\beta,\tol,\maxiter}$}
    \STATE{construct an initial admissible partition $\mesh[0]t$}
    \STATE{$l\mapsfrom 0$}
    \STATE{$\eta^2\mapsfrom \tol +1$}
    \WHILE{$l\leq\maxiter$ and $\eta^2>\tol$}
    \STATE{%
      \textsf{\textbf{solve}}
      for $(\fe[l]u,\vecfe[l]g,\matfe[l]h)\mapsfrom(\fespacefun vu,\vecfespacefun vg,\matfespacefun vh)$
      problem~(\ref{eq:discrete-system-finite-dim})
      with $\mesht\mapsfrom\mesh[l]t$
    }
    \FOR{${K} \in \mesh[l]t$}
    \STATE{
      compute $\eta({K})^2$ via (\ref{def:error-estimator})
    }
    \ENDFOR
    \STATE{%
      \textsf{\textbf{estimate}} by computing
      $\eta^2\mapsfrom\sums K{\mesh[l]t}\eta({K})^2$
    }%
    \STATE{%
      sort array
      $\seqs{\eta({K})^2}{K\in\mesh[l]t}$
      in decreasing order
    }
    \STATE{%
      \textsf{\textbf{mark}} the first $\ceil{\beta\card\mesh[l]t}$
      elements $K$ with the highest $\eta(K)^2$
    }
    \STATE{\textsf{\textbf{refine}}
      $\mesh[l]t\mapsto\mesh[{l+1}]t$ ensuring
      split of all marked elements and $l\mapsfrom l+1$
    } 
    \ENDWHILE
    \ENDPROCEDURE
  \end{algorithmic}
\end{Alg}
\section{Numerical experiments}
\label{sec:numerical-experiments}
This section reports on the numerical performance of the schemes
described in \S~\ref{sec:a-conforming-Galerkin-FEM}. We first
describe our numerical treatment of the zero tangential-trace condition
and introduce the intermediate finite element space $\fespace
x=\fespace u\times\tilde{\fespace g}\times\fespace h$ in
\S~\ref{sec:numerical-treatment-of-zero-tangential-trace}. We then
study four $\R2$-based experiments aimed at demonstrating the
robustness and testing the convergence rates of our method. In all
experiments the solution is known and computations are performed using
the FEniCS/Dolfin package~\citep{LoggMardalWells:12:book:Automated}. The
various error measures, include $\Norm{u-\fespacefun
  xu}_{\sobh1{(\W)}}$,
$\Norm{\grad{u}-\vecfespacefun{x}g}_{\sobh1{(\W)}}$,
$\Norm{\D^2u-\matfespacefun{x}H}_{\leb2{(\W)}}$ and
\({
  \Norm{(u,\grad{u},\D^2u)
    -
    (\fespacefun{x}u,\vecfespacefun{x}g,\matfespacefun{x}H)
  }_{\fatlinspace{Y}}
}\)
are plotted in logarithmic scale against the number of degrees of
freedom, ndof, that is the number of locations needed to store the
information on the computer. In test
problems~\ref{test-problem-nonzero-boundary},
\ref{test-problem-full-lower-order} and
\ref{test-problem-disk-domain}, the solution is chosen smooth
enough. The numerical results confirm the convergence analysis of
Theorem~\ref{the:convergence-rate}.
To benchmark our tests, we use the \indexemph{experimental orders of
  convergence} (\indexen{EOC}) associated with a numerical experiment
with errors $e_i$ and (uniform) mesh-sizes $h_i$, $i\integerbetween0I$,
which is defined by
\begin{equation}
  \EOC := \frac{\log(e_{i+1}/e_i)}{\log(h_{i+1}/h_i)}.
\end{equation}
\par
{%
  We also test the performance of the adaptive
  algorithm~\ref{alg:adaptive} in examples where the exact solution
  exhibits features such as rapid changes in localized parts of
  the domain and including a singularity as well. 
  For this we consider the test problem~\ref{test-problem-sharp-peak} as a problem 
  with a \indexemph{sharp peak} in the interior of domain and 
  test problem~\ref{test-problem-singularity} 
  as a problem with singular solution. 
  In these two cases, the convergence rate of the adaptive approach 
  with that of the uniform mesh refinement are compared.
}
\subsection{Numerical treatment of the zero tangential-trace}
\label{sec:numerical-treatment-of-zero-tangential-trace}
In proving the coercivity of $(u,\vec g,\mat H)\mapsto a_\theta(u,\vec
g,\mat H)$, and thus the error estimates, we took $\vec g
\in \linspace v$ (i.e., $\vec g$ is a Sobolev fields with vanishing
tangential-trace) and consequently $\vecfespacefun vg \in \fespace
g$. However, enforcing a zero tangential-trace condition onto the
finite element space is not trivial. One way to effect such a
boundary condition is to consider the appropriate constraint on the
space and introduce a Lagrange multiplier variable; in this case, we
must determine subspaces that satisfy the corresponding inf-sup
condition and this may limit the choice of finite element spaces. To
circumvent this limitation, based on the discussion in
\S~\ref{relaxing-problem}, we replace the zero-tangential-trace space
$\fespace g$, with the wider space $\tilde{\fespace g}$ in the
implementation and monitor the tangential-trace. Specifically, we
consider the discrete problem of finding
\({
  (\fespacefun{x}u,\vecfespacefun{x}g,\matfespacefun{x}H)
  \in
  \fespace{x}
  :=\fespace{u}\times\tilde{\fespace{g}}\times\fespace{h}
}\)
satisfying
\begin{equation}
  \label{eq:discrete-system-expr}
  \abil[a_\theta]{
    \fespacefun xu ,\vecfespacefun xg, \matfespacefun xH
  }{
    \fe\varphi , \vecfe\psi, \matfe\Xi
  }
  =
  \ltwop{
    f
  }{
    \linoptheta(\fe\varphi, \vecfe\psi, \matfe\Xi)
  }
  \Foreach (\fe\varphi, \vecfe\psi, \matfe\Xi) \in
  \fespace x
\end{equation}
corresponding to a zero boundary problem and
(\ref{eq:discrete-system-nonhomogeneous}) corresponding to a
nonzero boundary value.
\subsection{Test problem with nonzero boundary condition}
\label{test-problem-nonzero-boundary}
The first test problem considered by
\citet{LakkisPryer:11:article:A-finite}. Let
$\W=(-1,1)\times(-1,1)$ and
\begin{equation}
  \label{coef:test-nonzero}
\mat{A}(\vec x) = 
\begin{bmatrix}
  1 & 0 \\
  0 & a(\vec x)
\end{bmatrix},
\quad 
\vec b=[0,0],
\quad
c= 0,
\end{equation}
where $a(\vec x) = \arctan(5000(\norm{\vec x }^2 -1 ))+2 $. $\mat A(\vec x)$
satisfies the Cordes condition~(\ref{def:special-Cordes-condition}) with $\varepsilon = 0.37$. 
We choose right hand side $f$ and nonzero boundary condition $r$ such 
that the exact solution is
\begin{equation}
  u(\vec x) = \sin(\pic x_1) \sin(\pic x_2) + \sin (\pic(x_1+x_2)).
\end{equation}
We test the discrete problem (\ref{eq:discrete-system-nonhomogeneous})
for polynomial degree $k=1,2$ in uniform mesh.
Figure~\ref{fig:convergence-order-nonzero} bears results of the EOC. It clearly shows
that the method in used norms performs with optimal convergence rates.
\subsection{Test problem with full lower order terms}
\label{test-problem-full-lower-order}
In this test problem, let $\W = (-1, 1) \times(-1,1)$ and 
\begin{equation}
  \label{coef:test-full-lower-order}
  \mat{A}(\vec x) = 
  \begin{bmatrix}
    2 & \sign(x_1x_2)
    \\
    \sign(x_1x_2) & 2
  \end{bmatrix},
  \quad
  \vec b = [0.5, 0.5],
  \quad
  c = 1.
\end{equation}
We consider data $f$ such that the exact solution is
\begin{equation}
  u(\vec x) = x_1x_2(1-\exp(1-\norm{ x_1}))(1-\exp(1-\norm{ x_2 })).
\end{equation}
Although the secondary diagonal elements of $\mat A(\vec x)$ are
discontinuous on the axes, for $\lambda=1$, $\mat A(\vec x)$, $\vec b$
and $c$ satisfy the Cordes
condition~(\ref{def:general-Cordes-condition}) with $\varepsilon =0.22
$. We test the discrete problem (\ref{eq:discrete-system-expr}) for
$\theta = 0,0.5,1$ and polynomial degree $k=1,2$ in uniform
mesh.
\changefromto{
Fig.~\ref{fig:convergence-order-full-lower-theta-zero},
Fig.~\ref{fig:convergence-order-full-lower-theta-half} and
Fig.~\ref{fig:convergence-order-full-lower-theta-one} 
}{
Figs.~\ref{fig:convergence-order-full-lower-theta-zero}--\ref{fig:convergence-order-full-lower-theta-one}
}
show the optimal
convergence rates of the method through results of the EOC,
corresponding to $\theta=0,0.5,1$ respectively.
\subsection{Test problem in disk-domain}
\label{test-problem-disk-domain}
\begin{figure}[tbh]
  \begin{center}
  \caption[Test~\ref{test-problem-disk-domain} numerical solution]{
    \label{fig:disk-domain:solution}
    Test problem~\ref{test-problem-disk-domain}.  Numerically computed
    solution via discrete problem~(\ref{eq:discrete-system-expr}) with
    $\theta=0.5$ in the unit disk domain with
    coefficients~(\ref{coef:test-disk-domain}) and choosing the forcing
    $f$ such that
    \\
    $u(\vec x) = \sin(\pic(x_1^2+x_2^2))\cos(\pic(x_1-x_2)),$
    \\ by isoparametric $\poly{2}$-element,
    $k=2$, and $605973$ degrees of freedom.  }
  \includegraphics[scale=0.25]{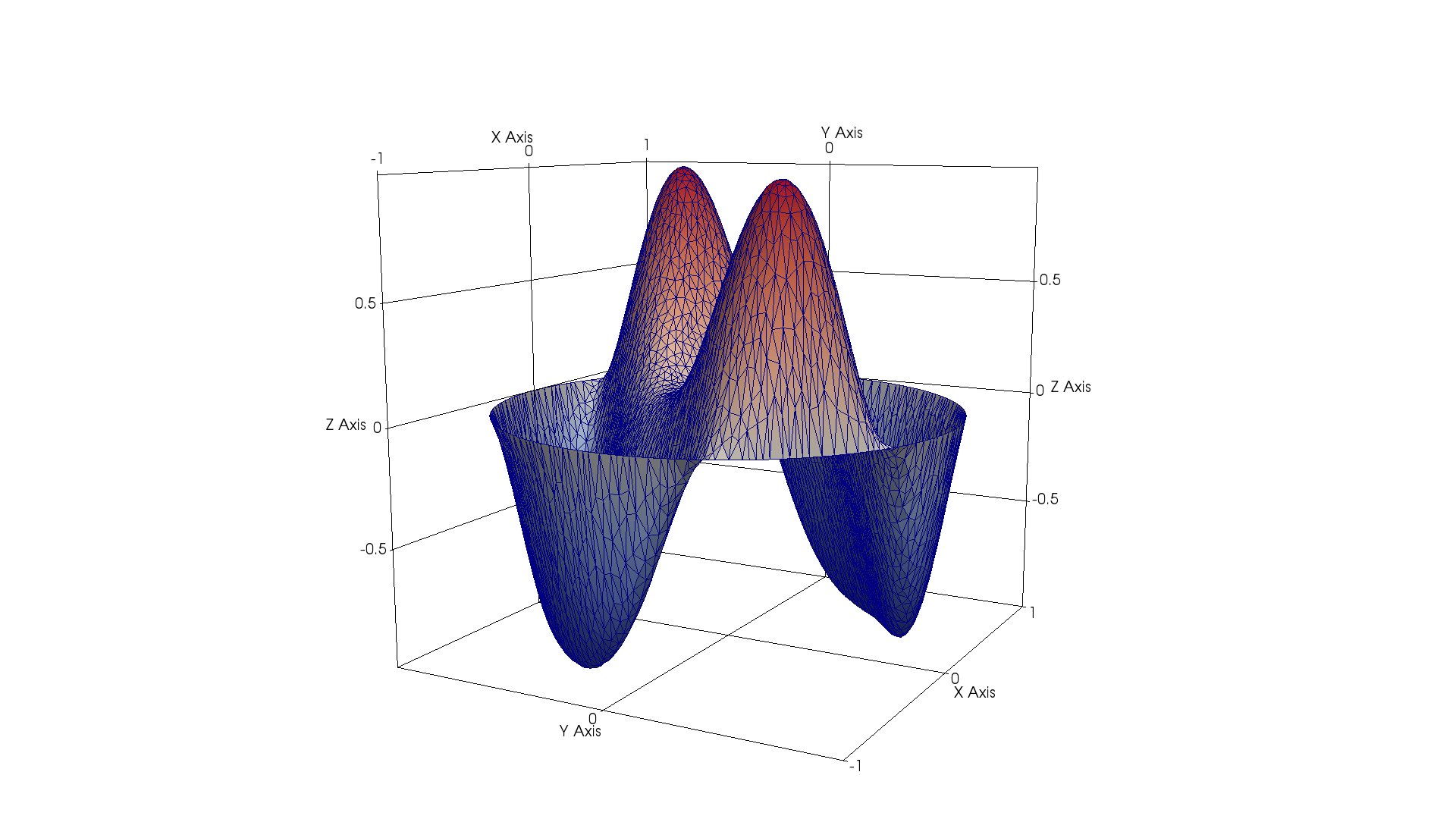}%
  \end{center}
\end{figure}

\begin{figure}[h]
  \centering
  \caption[Test~\ref{test-problem-disk-domain} error and rates]{
    \label{fig:convergence-order-disk-domain}
    Test problem~\ref{test-problem-disk-domain}.  We report the
    (log--log) error vs. degrees of freedom and the convergence rates
    for the discrete problem (\ref{eq:discrete-system-expr}) with
    $\theta=0.5$, applied to a nondivergence form problem in the unit
    disk domain, with coefficients (\ref{coef:test-disk-domain}) and
    choosing the forcing $f$ such that
    \\ $u(\vec{x})=\sin(\pic(x_1^2+x_2^2)) \cos(\pic(x_1-x_2))$.
    \\ For $\poly k$ elements with both $k=1$ and $2$, we observe
    optimal convergence rates, that is
    \({
      \Norm{u-\fespacefun{x}u}_{\sobh1(\W)} =
      \Norm{\grad{u}-\vecfespacefun{x}g}_{\sobh1(\W)} =
      \Norm{\D^2{u}-\matfespacefun{x}H}_{\leb2(\W)} =\Oh(h^k).
    }\)
    For $k=2$, the isoparametric finite element is used.
  }
  \subfloat[{\(\poly{1}\) elements}]{
    \includegraphics[scale=.37]{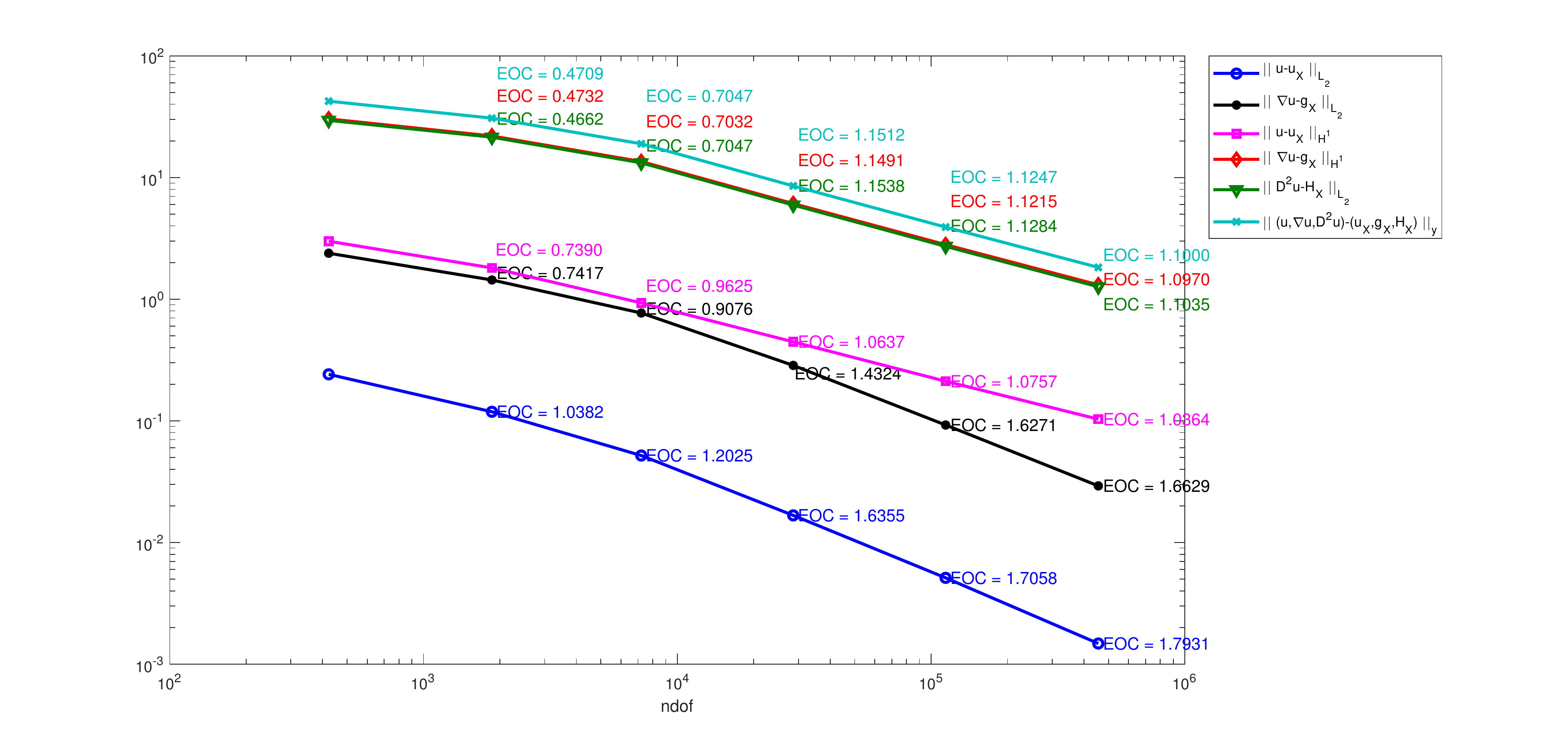}}
  \\
  \subfloat[{\(\poly{2}\) elements}]{
    \includegraphics[scale=.37]{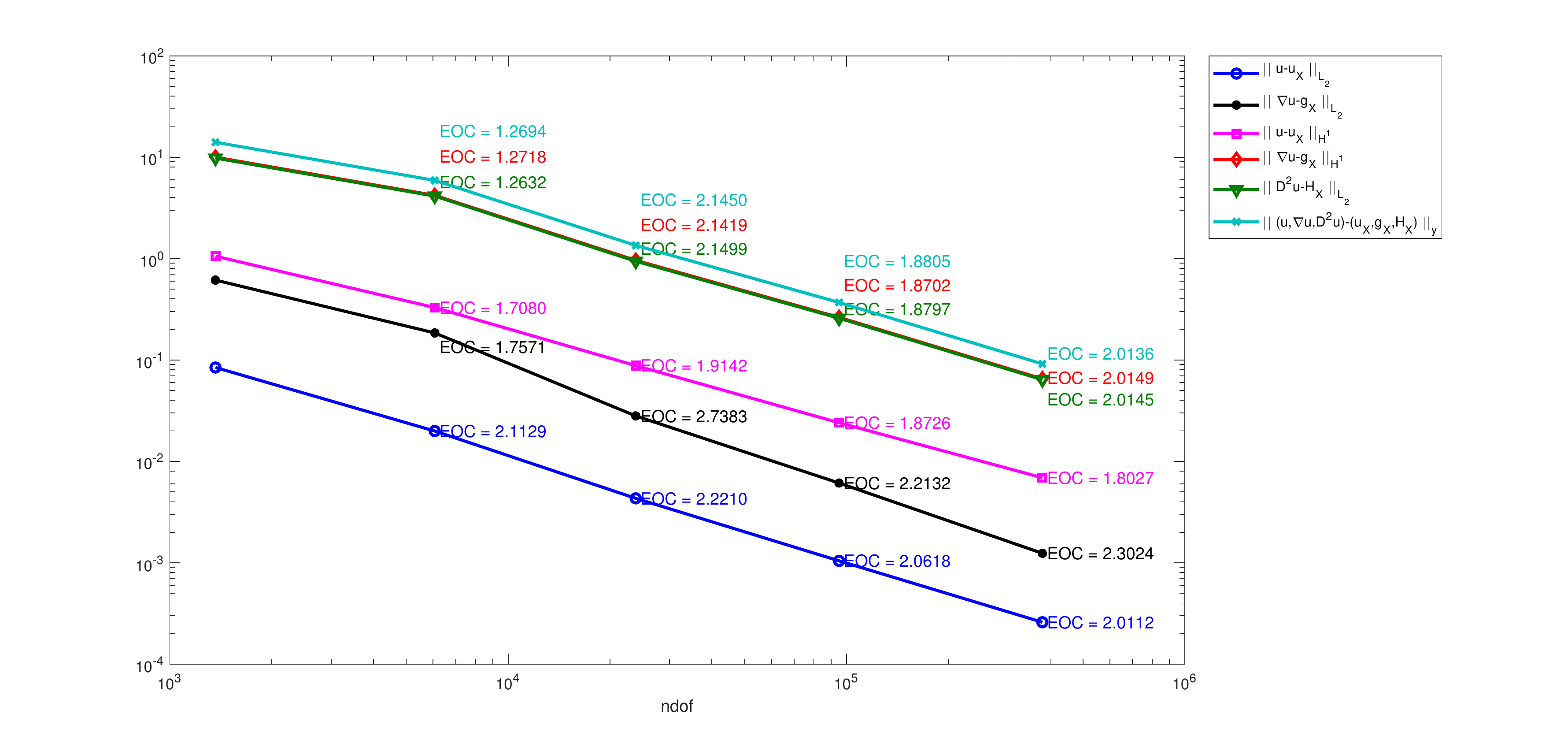}}
\end{figure}

\begin{figure}[tbh]
  \begin{center}
  \caption[Adaptive test~\ref{test-problem-sharp-peak} mesh.]{%
      \label{fig:adaptive-peak-mesh}
      Test problem~\ref{test-problem-sharp-peak}. 
      Adaptively refined mesh, generated by Algorithm~\ref{alg:adaptive}
      with $\beta=0.3$ and after $8$ iterations, for polynomial degree 
      $k=2$ ( and $122598$ degrees of freedom).
  }
    \includegraphics[scale=.5]{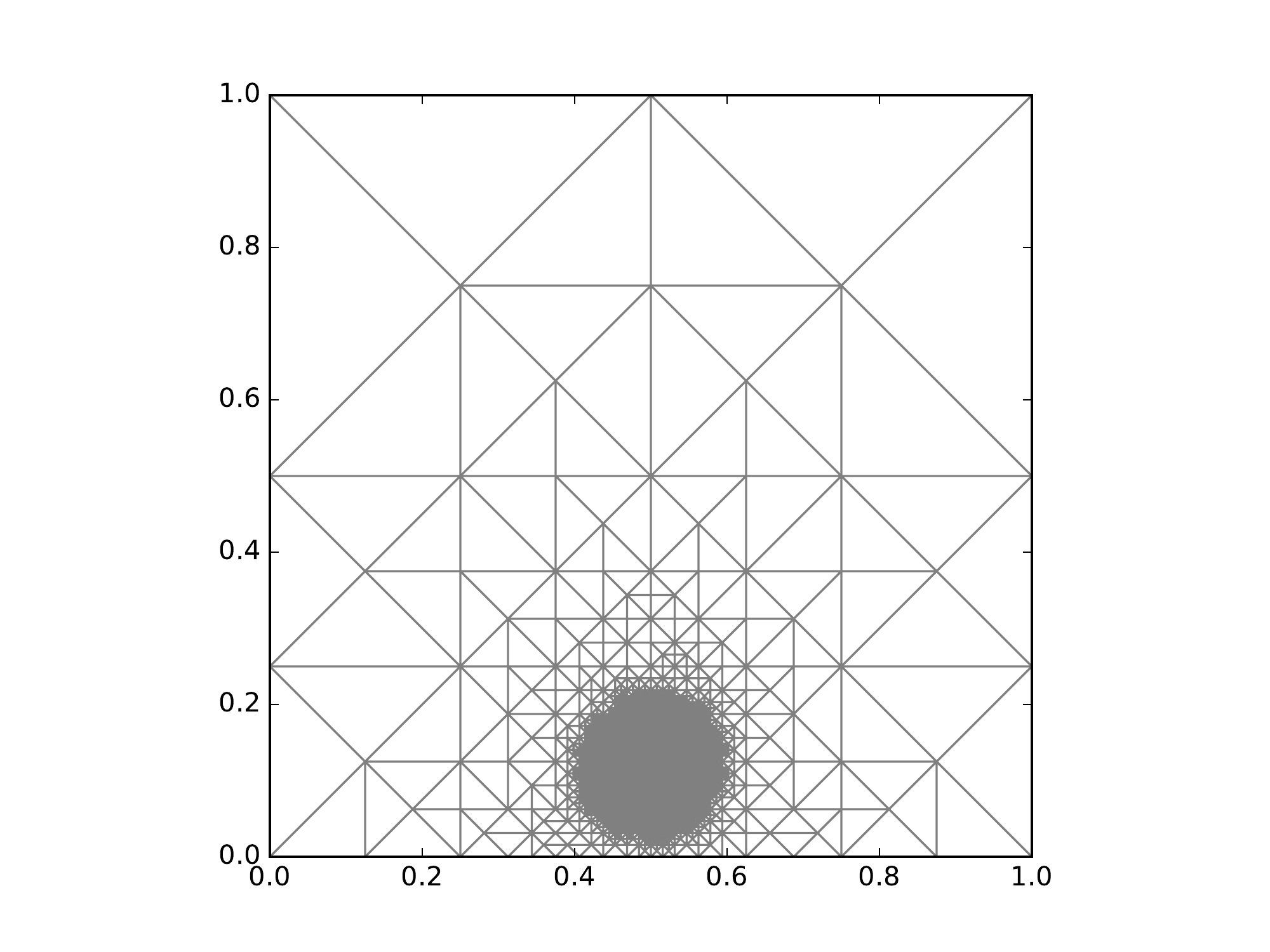}
  \end{center}
\end{figure}

\begin{figure}[ht]
  \begin{center}
  \caption[Adaptive test~\ref{test-problem-singularity}: mesh.]{%
    {%
      \label{fig:adaptive-singular-mesh}
      Test problem~\ref{test-problem-singularity}. 
      Adaptively refined mesh, generated by Algorithm~\ref{alg:adaptive}
      with $\beta=0.3$ and after $8$ iterations, for polynomial degree 
      $k=2$ ( and $98679$ degrees of freedom).
    }%
  }
  \includegraphics[scale=.5]{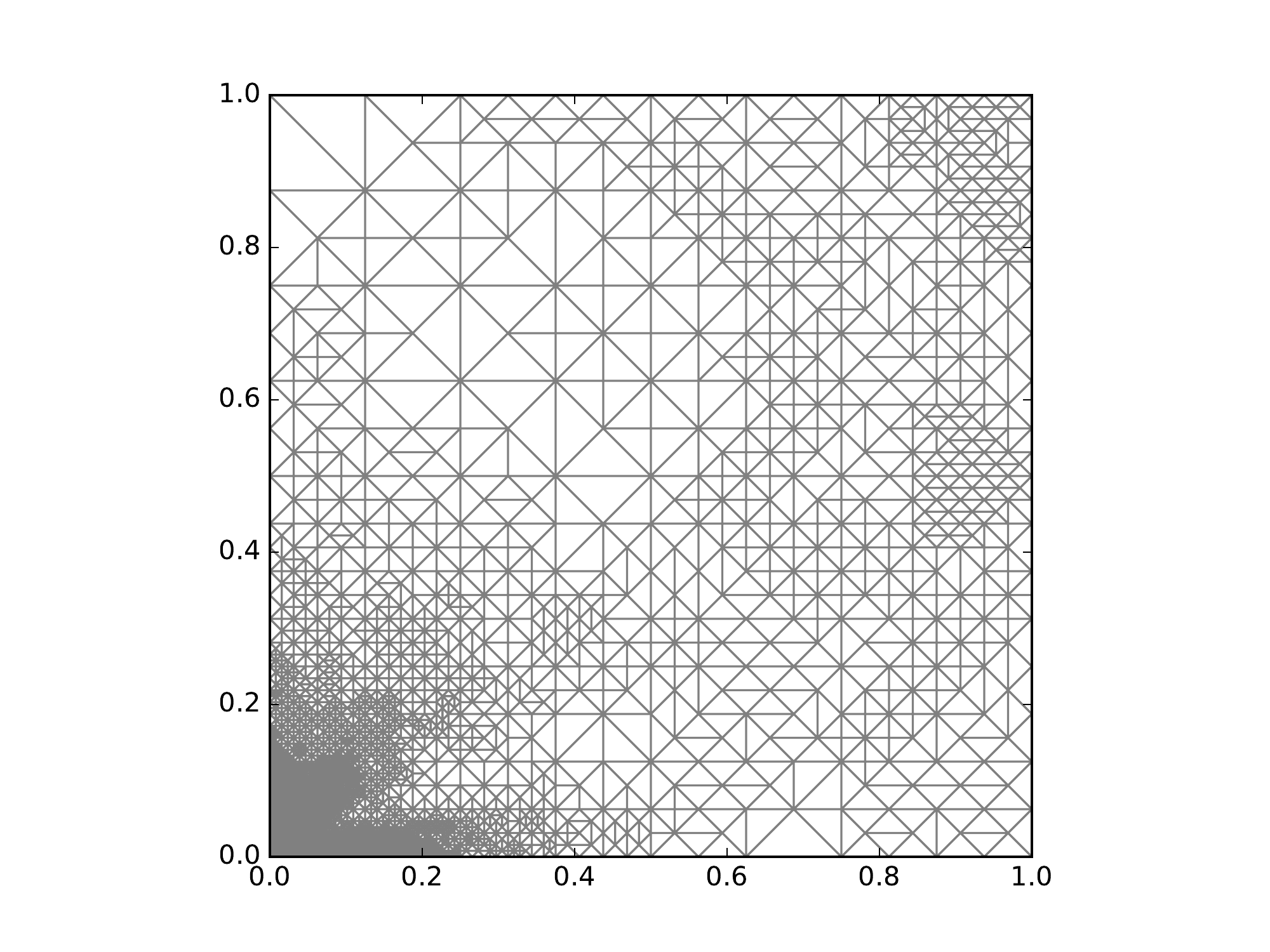}
  \end{center}
\end{figure}

In this test problem, let $\W $ be the unit disk domain and 
\begin{equation}
\label{coef:test-disk-domain}
\mat{A}  = 
\begin{bmatrix}
2 & 1\\
1 & 1
\end{bmatrix},
\quad
\vec b(\vec x)=[x_1x_2,0],
\quad
c = 2.
\end{equation}
For $\lambda =1$, these data satisfy the Cordes 
condition~(\ref{def:general-Cordes-condition}) 
with $\varepsilon =0.17 $.
 We choose data $f$ such that the exact solution is
\begin{equation}
u(\vec x) = \sin(\pic(x_1^2+x_2^2)) \cos(\pic(x_1-x_2)).
\end{equation}
We test the discrete problem (\ref{eq:discrete-system-expr}) for $\theta = 0.5$ and 
polynomial degree $k=1,2$ in {unstructured} quasi-uniform mesh. Since the 
domain includes curved boundary, for $k = 2$, we use isoparametric finite element. 
The approximate solution is shown in Fig.~\ref{fig:disk-domain:solution} and results 
of the EOC for the approximation can be found in Fig.~\ref{fig:convergence-order-disk-domain}, 
which demonstrates the optimal convergence rates of the method.
{%
  \subsection{Numerical results of the adaptive refinement}
  In the following examples, we test the performance of the adaptive refinement 
  based on Algorithm~\ref{alg:adaptive}. We set refinement 
  fraction $\beta = 0.3$, tolerance $\tol = 10^{-6}$ and maximum number of iteration 
  $\maxiter =12$ such that on each $\mesh[k]t$, the discrete problem 
  (\ref{eq:discrete-system-expr}) with $\theta = 0.5$ and polynomial degree $k=1,2$ is 
  applied. In the all following test problems, we also set the coefficients as
  \begin{equation}
    \label{coef:tests-adaptive}
    \mat{A}(\vec x) = 
    \begin{bmatrix}
      1 & (x_1x_2)^{2/3}\\
      (x_1x_2)^{2/3} & 4
    \end{bmatrix}, 
    \quad 
    \vec b(\vec x)= \left[ (x_1x_2)^{1/3}, (x_1x_2)^{1/3}\right] ,
    \quad
    c = 2.
  \end{equation}
  For $\lambda =1$, these data satisfy the Cordes 
  condition~(\ref{def:general-Cordes-condition}) 
  with $\varepsilon =0.04 $, in the considered domains.
  \subsection{Test problem with sharp peak}
  \label{test-problem-sharp-peak}
  In this test problem let $\W = (0, 1) \times(0,1)$ and choose 
  data $f$ such that the exact solution is
  \begin{equation}
    u(\vec x) = x_1x_2(x_1-1)(x_2-1) \exp(-1000((x_1-0.5)^2+(x_2-0.117)^2)).
  \end{equation}
  The solution includes sharp peak at $(x_1, x_2)= (0.5, 0.117)$. 
  An obvious remedy to deal with this difficulty is to refine the discretization near the 
  critical regions. The adaptive refined mesh is shown in Fig.~\ref{fig:adaptive-peak-mesh}. 
  To demonstrate the performance of the adaptive refinement, we compare the error of the 
  method in uniform with adaptive mesh for polynomial degree $k =1,2$ in Fig.~\ref{fig:adaptive-peak-1} 
  and Fig.~\ref{fig:adaptive-peak-2} respectively.
  \subsection{Test problem with a salient corner singularity}
  \label{test-problem-singularity}
  In this test problem let $\W = (0, 1) \times(0,1)$ and choose 
  data $f$ such that the exact solution is
  \begin{equation}
    u(\vec x) =2(x_1-x_1^2)(x_2-x_2^2)\pownorm{-\fracl12}{\vec x}
  \end{equation}
  and has thus a singularity at $(0,0)$.
  One should note that 
  $u \in \sobh{s}(\W)$ for $s <1+3/2$. As we see in Fig.~\ref{fig:adaptive-singularity-1} 
  and Fig.~\ref{fig:adaptive-singularity-2}, singularity of solution 
  $u(\vec x)$ at $(0,0)$ leads to lack of optimal convergence rate on uniform mesh. 
  Through the adaptive approach, we expect an improvement of the convergence rates
  (at least) for 
  $\Norm{\grad{u}-\vecfespacefun{x}g}_{\sobh1(\W)}$, 
  $\Norm{\D^2 u-\matfespacefun xH}_{\leb2(\W)}$ and
  $\Norm{(u,\grad{u},\D^2u)-(\fespacefun{x}u,\vecfespacefun{x}g,\matfespacefun{x}H)}_{\fatlinspace Y}$. 
  The adaptive refined mesh is shown in 
  Fig.~\ref{fig:adaptive-singular-mesh}. We compare the error of the method in uniform 
  with adaptive mesh for polynomial degree $k =1,2$ in Fig.~\ref{fig:adaptive-singularity-1} 
  and Fig.~\ref{fig:adaptive-singularity-2} respectively.
}
\section{Conclusions and outlook}
\label{sec:conclusion}
The least-squares based gradient or Hessian
recovery method presented is a practical and effective method for the
numerical approximation of solutions to linear elliptic equations in
nondivergence form. The advantages of the method herewith proposed
are:
\begin{enumerate}[(a)\ ]
\item
  Method~(\ref{eq:discrete-system-finite-dim}) allows the use of a
  wide choice of finite elements, including all standard conforming.
  With the appropriate modifications one could envisage extending our
  method to nonconforming elements as well, e.g.,
  \citet{SmearsSuli:13:article:Discontinuous}.
\item
  \label{item:comparison-with-Gallistl-Suli:conclusion}
  Our least squares Lax--Milgram-based approach circumvents the need
  for {Lagrange multipliers or a curl-penalty stabilization}
  in inf-sup stable combinations 
  for $(\fe u,\vecfe g)$ (let alone
  $(\fe u,\vecfe g,\matfe h)$ when the Hessian is needed)
  as in \citet{Gallistl:17:article:Variational},
  {\cite{Gallistl:19:article:Numerical}
  and \citet{GallistlSuli:19:article:Mixed}. 
  We also can use a Céa quasi-optimality in the error analysis.}
  {
  \item
    Through the least-squares approach, we are capable of considering
    constraints that are assumed on the function spaces (to ensure well-posedness of the problem) 
    as square terms of the quadratic cost functional and then working in 
    general function spaces.  }
\item
  We are able to derive straightforward \aposteriori error bounds,
  with easily implemented estimators and indicators for which
  the adaptive method shows convergence.
\item
  We can choose between a gradient-and-Hessian and gradient-only
  recovery as observed in \S~\ref{obs:Hessian-less-approach}. The
  Hessian is useful when our method is applied as the linear look
  within a Newton or fixed-point method to a nonlinear elliptic
  equation as in \citet{LakkisPryer:13:article:A-finite},
  \cite{Neilan:14:article:Finite},
  \cite{LakkisPryer:15:inproceedings:An-adaptive} and
  \cite{KaweckiLakkisPryer:18:article:A-finite}.
\item\label{item:conclusions:using-dg}
  {
    An interesting issue, which we did not have room to address in this
    paper, is the use of discontinuous Galerkin piecewise polynomial spaces
    for the approximation of the gradient or the Hessian.  Our method,
    at least from the computational side can be easily adapted to
    use such spaces, but the outcomes and gains are not clear, in that
    the analysis would need serious reworking and the penalization
    parameters required to get coercivity going might just give an unexpected
    sting in the tail.}
\end{enumerate}
Our method is not without drawbacks of which we note the lack of
optimal convergence rate for the function value error
$\Norm{u-\fespacefun{v}u}_{\leb2(\W)}$ {and the slow convergence
  for viscosity solutions (which we have not included in this work)}. We are aiming to
address thes issues in forthcoming work {announced by
\cite{LakkisMousavi:20:techreport:A-least-squares}}.

Our FEniCS-based implementation is available on request for testing
and further research.

\providecommand{\figuniformtest}[9]{
  \begin{figure}[h]
    \centering
    \caption[Test~\ref{#1}, $\theta=#3$]{%
      \label{#2}
      Test problem~\ref{#1}. We report the (log--log) error
      vs. degrees of freedom and the convergence rates for the
      discrete problem (\ref{#8})%
      , applied to a
      nondivergence form problem
      (\ref{eq:nondivergence0-inhomogeneous}) with domain $\W=#4$,
      coefficients (\ref{coef:#5}) and choosing right hand side $f$
      such that
      \\
      $u(\vec{x})=#6$.
      \\
      For $\poly k$ elements with
      both $k=1$ and $2$, we observe optimal convergence rates, that
      is
      \\
      \({
        \Norm{u-\fespacefun{#9}u}_{\sobh1(\W)} =
        \Norm{\grad{u}-\vecfespacefun{#9}g}_{\sobh1(\W)} =
        \Norm{\D^2{u}-\matfespacefun{#9}H}_{\leb2(\W)} =\Oh(h^k).
      }\)%
    }
    \subfloat[{$\poly{1}$ elements}]{
      \includegraphics[scale=.37]{Picture/#7_1.pdf}}
    \\
    \subfloat[{$\poly{2}$ elements}]{
      \includegraphics[scale=.37]{Picture/#7_2.pdf}}
  \end{figure}
}
\figuniformtest{%
  test-problem-nonzero-boundary}{%
  fig:convergence-order-nonzero}{%
  0.5}{%
  \opinter{-1}1^2}{%
  test-nonzero}{%
  \sin(\pic x_1)\sin(\pic x_2)+\sin(\pic(x_1+x_2))}{%
  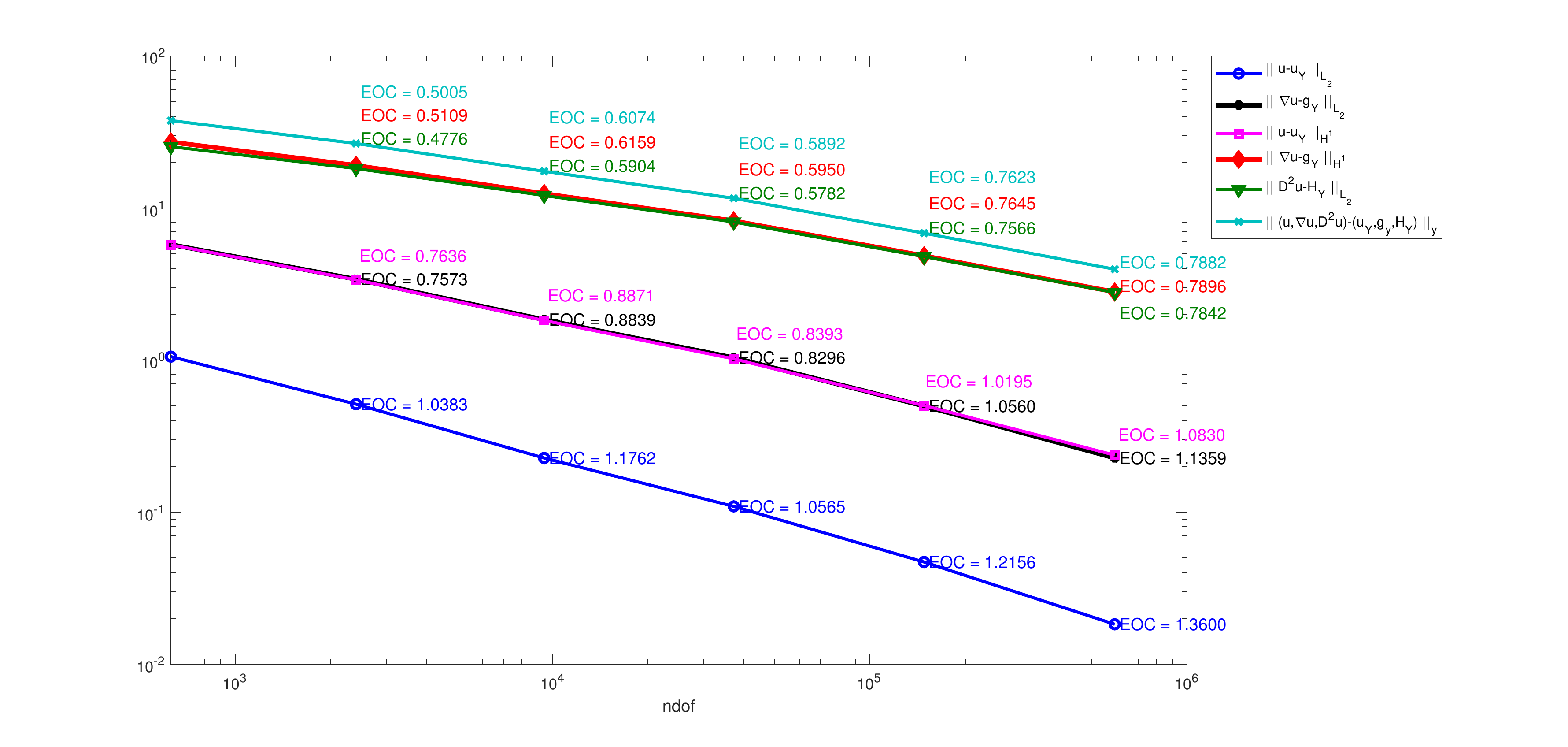}{%
  eq:discrete-system-nonhomogeneous}{y}
\figuniformtest{%
  test-problem-full-lower-order}{fig:convergence-order-full-lower-theta-zero}{%
  0}{%
  \opinter{-1}1^2}{%
  test-full-lower-order}{%
  x_1x_2\qp{1-\expp{1-\norm{x_1}}}\qp{1-\expp{1-\norm{x_2}}}}{%
  test2_theta_zero_k}{%
  eq:discrete-system-expr}{x}
\figuniformtest{%
  test-problem-full-lower-order}{fig:convergence-order-full-lower-theta-half}{%
  0.5}{%
  \opinter{-1}1^2}{%
  test-full-lower-order}{%
  x_1x_2(1-\exp(1-\norm{x_1}))(1-\exp(1-\norm{ x_2 }))}{%
  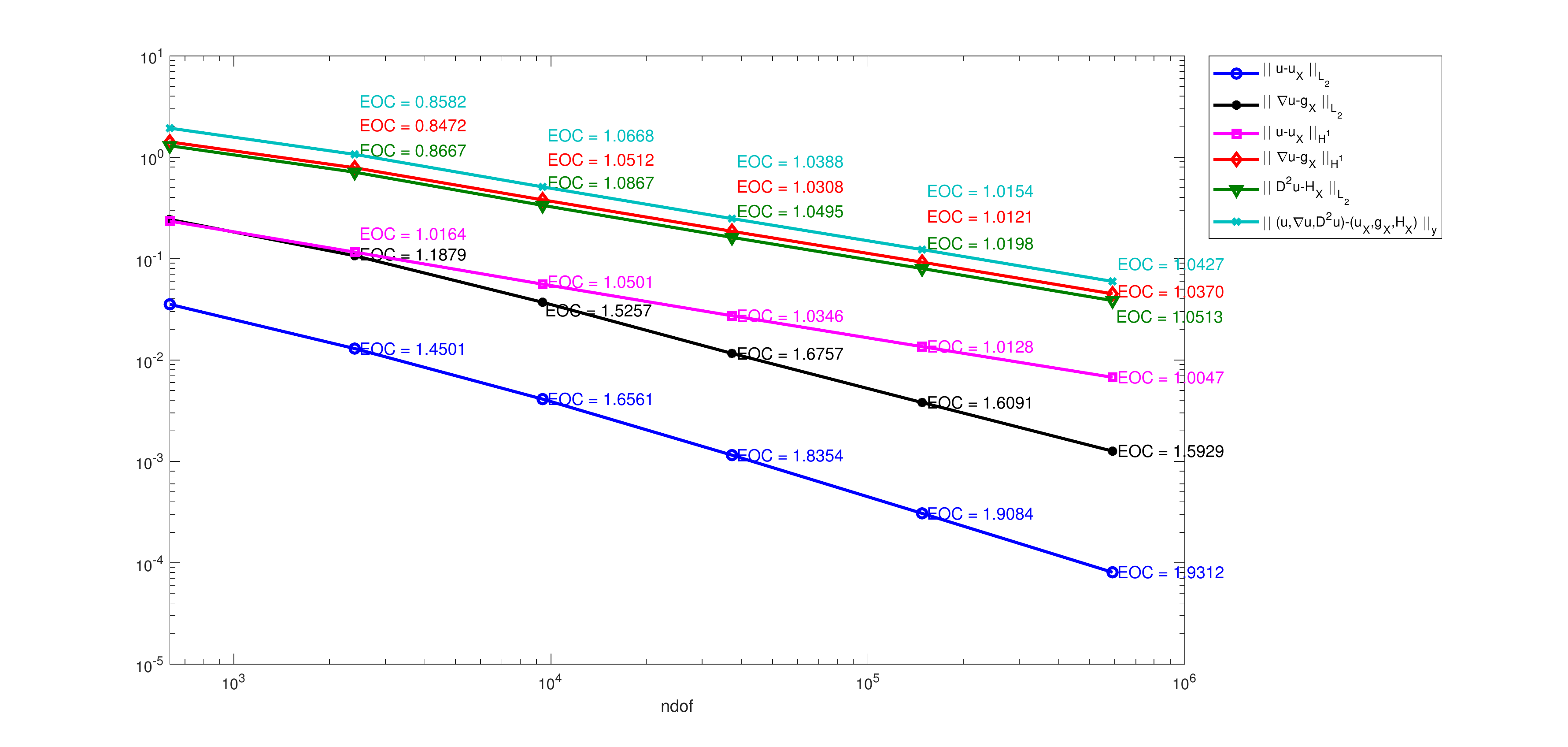}{%
  eq:discrete-system-expr}{x}
\figuniformtest{%
  test-problem-full-lower-order}{fig:convergence-order-full-lower-theta-one}{%
  1}{%
  \opinter{-1}1^2}{%
  test-full-lower-order}{%
  x_1x_2(1-\exp(1-\norm{x_1}))(1-\exp(1-\norm{ x_2 }))}{%
  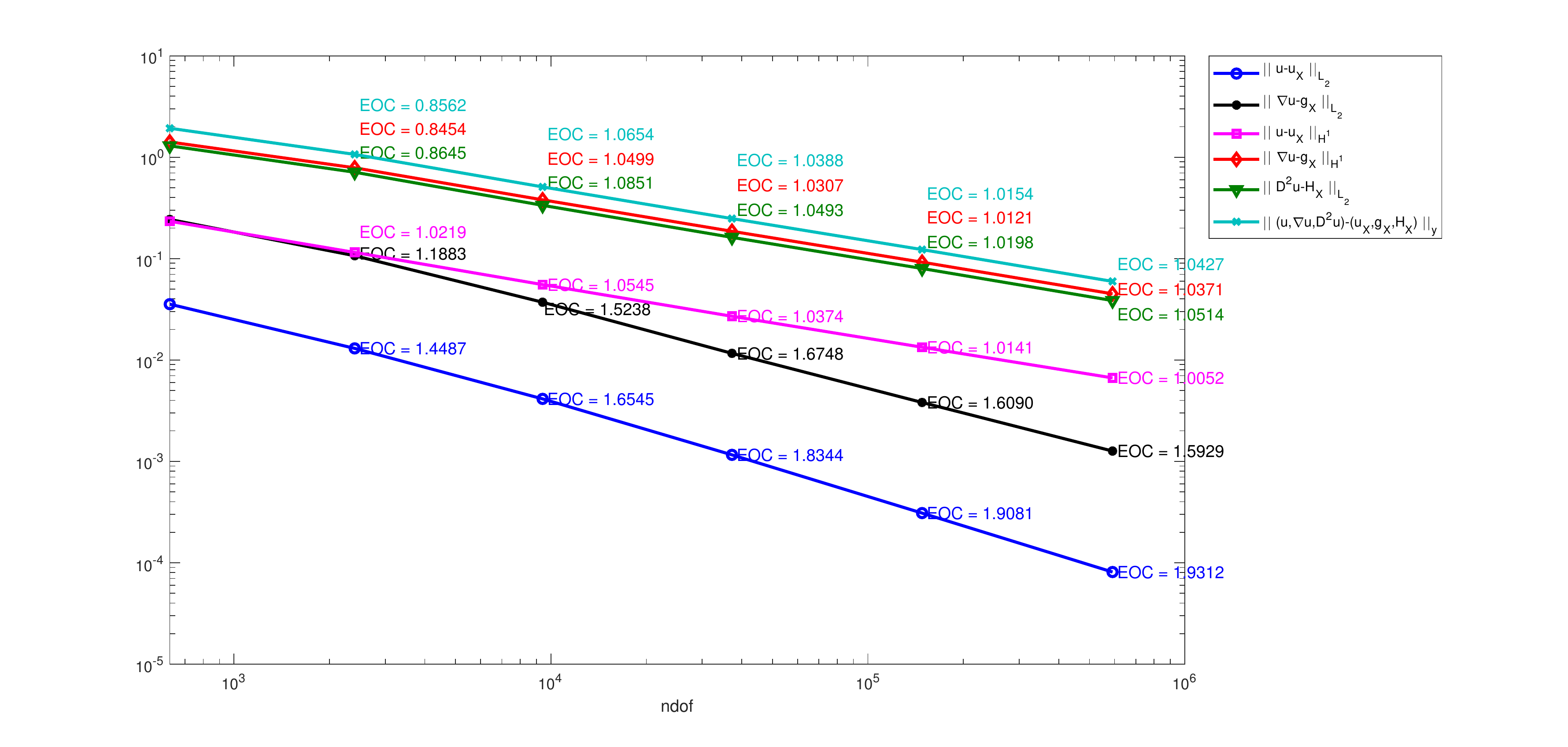}{%
  eq:discrete-system-expr}{x}
%
%
%
%
%
%
%
%
%
%
%
%
%
%
%
%
%
%
%
%
%
%
%
%
%
%
%
%
%
%
%
%
%
%
%
%
%
%
%
%
%
%
%
%
%
%
%
%
%
%
%
%
%
%
%
%
%
%
%
%
%
%
%\subfloat[]{\includegraphics[scale=.5]{test4_k_2_g.pdf}}\\%{test4_m_2_g.pdf}}\\
%
%
%
%
\providecommand{\figpeaksolution}[2]{}
\renewcommand{\figpeaksolution}[2]{%
  \begin{figure}[ht]
    \centering
    \caption[Adaptive test~\ref{test-problem-sharp-peak}: rates with $\poly{#1}$ elements]{%
      \label{fig:adaptive-peak-#1}
      Adaptive mesh refinement Algorithm~\ref{alg:adaptive} on
      problem~\ref{test-problem-sharp-peak} with $\poly{#1}$ elements.
      We plot the errors in various norms of the discrete problem 
      (\ref{eq:discrete-system-expr}) with $\theta =0.5$ for
      {\color{blue}uniform} and {\color{red}adaptive} mesh, on the
      domain $\W=(0,1)\times(0,1)$ with coefficients
      (\ref{coef:tests-adaptive}) and exact solution
      \\
      $u(\vec x) = x_1x_2(x_1-1)(x_2-1) \exp(-1000((x_1-0.5)^2+(x_2-0.117)^2))$.
      #2%
    }
    \subfloat[full norm][%
      Full error-norm:
      {
        \tiny
        \({\Norm{(u,\grad{u},\D^2u)
            -
            (\fespacefun{x}u,\vecfespacefun{x}g,\matfespacefun{x}H)}_{\fatlinspace Y}
        }\)%
      }%
    ]{
      \includegraphics[scale=.35]{Picture/test4_#1_Y.pdf}}%
    \subfloat[L2 norm][$\Norm{u-\fespacefun xu}_{\leb2(\W)}$]
             {%
               \includegraphics[scale=.35]{Picture/test4_#1_L2u.pdf}}%
    \\
    \subfloat[H1 seminorm][$\norm{u-\fespacefun xu}_{\sobh1(\W)}$]
        {%
          \includegraphics[scale=.35]{Picture/test4_#1_SH1u.pdf}}%
    \subfloat[L2 of recovered gradient norm][$\Norm{\grad{u}-\fespacefun xg}_{\leb2(\W)}$]
        {%
          \includegraphics[scale=.35]{Picture/test4_#1_L2g.pdf}}%
    \\
    \subfloat[H1 of recovered gradient seminorm][$\norm{\grad{u}-\fespacefun xg}_{\sobh1(\W)}$]
        {%
          \includegraphics[scale=.35]{Picture/test4_#1_SH1g.pdf}}%
    \subfloat[L2 of recovered Hessian norm][$\Norm{\D^2 u-\matfespacefun xH}_{\leb2(\W)}$]
        {%
          \includegraphics[scale=.35]{Picture/test4_#1_H.pdf}}%
  \end{figure}
}
\figpeaksolution1{%
  Although uniform and adaptive errors seem asymptotically equivalent
  (because the solution is not really singular), the
  adaptive error is an order of magnitude smaller.%
}
\figpeaksolution2{%
  Compared to \ref{fig:adaptive-peak-1}, also in this
  case we see that despite their asymptotic equivalence, the adaptive
  error in all norms becomes an order of magnitude smaller than the
  uniform error after $8$ iterations. The higher polynomial degree
  makes this shift more pronounced.
}
\renewcommand{\figpeaksolution}[2]{
  \begin{figure}[ht]
    \centering
    \caption[Adaptive test~\ref{test-problem-singularity}: rates with $\poly{#1}$ elements]{
      \label{fig:adaptive-singularity-#1}
      Adaptive mesh refinement Algorithm~\ref{alg:adaptive} on
      problem~\ref{test-problem-singularity} with $\poly{#1}$ elements.
      We plot the errors in various norms of the discrete problem 
      (\ref{eq:discrete-system-expr}) with $\theta =0.5$ for
      {\color{blue}uniform} and {\color{red}adaptive} mesh, on the
      domain $\W=(0,1)\times(0,1)$ with coefficients
      (\ref{coef:tests-adaptive}) and exact solution
      \ifx#2{}\else\\{#2}\fi%
    }
    \subfloat[full error norm][{%
        Full error-norm:
        \\{
          \tiny
          \({\Norm{(u,\grad{u},\D^2u)
              -
              (\fespacefun{x}u,\vecfespacefun{x}g,\matfespacefun{x}H)}_{\fatlinspace Y}
          }\)%
        }
    }]{\includegraphics[scale=.35]{Picture/test5_#1_Y.pdf}}%
    \subfloat[L2 error norm][{%
        $\Norm{u-\fespacefun xu}_{\leb2(\W)}$%
    }]{\includegraphics[scale=.35]{Picture/test5_#1_L2u.pdf}}%
    \\
    \subfloat[H1 error seminorm][{%
        $\norm{u-\fespacefun xu}_{\sobh1(\W)}$%
    }]{\includegraphics[scale=.35]{Picture/test5_#1_SH1u.pdf}}%
    \subfloat[L2 of recovered gradient norm][{%
        $\Norm{\grad{u}-\fespacefun xg}_{\leb2(\W)}$%
    }]{\includegraphics[scale=.35]{Picture/test5_#1_L2g.pdf}}%
    \\
    \subfloat[H1 or recovered gradient norm][{%
        $\norm{\grad{u}-\fespacefun xg}_{\sobh1(\W)}$%
    }]{\includegraphics[scale=.35]{Picture/test5_#1_SH1g.pdf}}%
    \subfloat[L2 of recovered Hessian norm][{%
        $\Norm{\D^2 u-\matfespacefun xH}_{\leb2(\W)}$%
    }]{\includegraphics[scale=.35]{Picture/test5_#1_H.pdf}}%
  \end{figure}
}
\figpeaksolution 1{Although the performance of $\poly1$ elements is not the best,
  this example shows that the gradient is better approximated in the $\sobh1(\W)^2$ norm.}
\figpeaksolution 2{The superiority of the $\poly2$ elements in combination
  with the adaptive method versus the uniform $\poly2$ elements is clearly exhibited here,
  especially in the approximation of the gradient and the Hessian.}
\bibliographystyle{abbrvnat}
\ifthenelse{\boolean{shownotes}}{%
  \clearpage%
  \printindex%
}{}%
\end{document}